\documentclass[11 pt, reqno]{article}

\usepackage{amsmath,amsfonts,amssymb,amsthm}
\usepackage[colorlinks=true,hyperindex=true]{hyperref}
\usepackage{cancel,bbm}
\usepackage{mathabx} 
\usepackage{comment}
\usepackage{enumitem}
\usepackage{cite}

\usepackage{chngcntr}
\counterwithin*{equation}{section}

\usepackage[left=3cm, right=3cm, top=3 cm, bottom= 3 cm]{geometry}
\usepackage{color}
\usepackage{fancyhdr}
\usepackage{latexsym}

\usepackage{bm}

\newtheorem{ccounter}{ccounter}[section]
\newtheorem{thm}[ccounter]{Theorem}
\newtheorem{lem}[ccounter]{Lemma}
\newtheorem{cor}[ccounter]{Corollary}
\newtheorem{defn}[ccounter]{Definition}
\newtheorem{prop}[ccounter]{Proposition}
\newtheorem{ass}[ccounter]{Assumption}
\newtheorem{ex}[ccounter]{Example}

\newtheorem{rmk}[ccounter]{Remark}

\def\bet{\begin{thm}}
\def\eet{\end{thm}}
\def\bel{\begin{lem}}
\def\eel{\end{lem}}
\def\bas{\begin{ass}}
\def\eas{\end{ass}}
\def\bec{\begin{cor}}
\def\eec{\end{cor}}
\def\bed{\begin{defn}}
\def\eed{\end{defn}}
\def\bep{\begin{prop}}
\def\eep{\end{prop}}
\def\beq{\begin{equation}}
\def\eeq{\end{equation}}
\def\proof{\noindent {\bf Proof.}\ \ }
\def\bea{\begin{equation*}}
\def\eea{\end{equation*}}
\def\tr{\mathrm{tr}}
\def\bex{\begin{ex}}
\def\eex{\end{ex}}
\def\remark{\noindent{\bf Remark. }}

\def\rr{\mathbb{R}}

\def\cc{\mathbb{C}}
\def\1{\boldsymbol{1}}
\def\Im{\mathrm{Im}}
\def\Re{\mathrm{Re}}
\def\e{\mathrm{e}}
\def\i{\mathrm{i}}
\def\del{\partial}
\def\d{\mathrm{d}}
\def\eps{\varepsilon}
\renewcommand\leq\varleq
\renewcommand\geq\vargeq
\def\ee{\mathrm{E}}

\def\O{\mathcal{O}}

\def\ee{\mathbb{E}}

\def\pp{\mathbb{P}}

\def\msc{m_{\mathrm{sc}}}
\def\rhosc{\rho_{\mathrm{sc}}}
\def\D{\mathcal{D}}
\def\I{\mathcal{I}}
\def\mfa{\mathfrak{m}}
\def\A{\mathcal{A}}

\def\J{\mathcal{J}}

\def\tilz{\tilde{z}}
\def\Var{\mathrm{Var}}
\def\bx{\mathbf{x}}
\def\by{\mathbf{y}}
\def\S{\mathcal{S}}
\def\dto{\downarrow}
\def\D{\mathcal{D}}
\def\C{\mathcal{C}}
\def\Cov{\mathrm{Cov}}
\def\be{\mathbf{e}}
\def\tilphi{\tilde{\varphi}}
\def\mfa{\mathfrak{a}}
\def\Oma{\Omega_{\mathfrak{a}}}
\def\tilf{\tilde{f}}
\def\ea{e_{\mathfrak{a}}}
\def\psia{\psi_{\mathfrak{a}}}
\def\mfb{\mathfrak{b}}
\def\mfs{\mathfrak{s}}
\def\tilZ{\tilde{Z}}
\def\hatmu{\hat{\mu}}
\def\hatx{\hat{x}}
\def\haty{\hat{y}}
\def\tilz{\tilde{z}}
\def\tilw{\tilde{w}}
\def\muij{\hat{\mu}^{(ij)} }
\def\lamij{\hat{\lambda}^{(ij)}}
\def\Ai{\mathrm{Ai}}
\def\phiN{\varphi_N^{(N)}}
\def\phiNm{\varphi_{N-1}^{(N)}}
\def\sg{\mathrm{sg}}
\def\mfc{\mathfrak{c}}
\def\N{\mathcal{N}}

\begin{document}
\title{Deformed GOE}

\begin{table}
\centering

\begin{tabular}{c}

\multicolumn{1}{c}{\parbox{12cm}{\begin{center}\Large{\bf Almost-optimal bulk regularity conditions in the CLT for Wigner matrices  }\end{center}}}\\
\\
\end{tabular}
\begin{tabular}{ c c c  }
Benjamin Landon
& \phantom{blah} & 
Philippe Sosoe
 \\
 & & \\  
 \small{University of Toronto} & & \small{ Cornell University } \\
 \small{Department of Mathematics} & & \small{Department of Mathematics} \\
 \small{\texttt{blandon@math.toronto.edu}} & & \small{\texttt{ps934@cornell.edu}} \\
  & & \\
\end{tabular}
\\
\begin{tabular}{c}
\multicolumn{1}{c}{\today}\\
\\
\end{tabular}

\begin{tabular}{p{15 cm}}
\small{{\bf Abstract:}  We consider linear spectral statistics of the form $\tr ( \varphi (H))$ for test functions $\varphi$ of low regularity and Wigner matrices $H$ with smooth entry distribution. We show that for functions $\varphi$ in the Sobolev space $H^{1/2+\eps}$ or the space $C^{1/2+\eps}$,  that are supported within the spectral bulk of the semicircle distribution, these linear spectral statistics have asymptotic Gaussian fluctuations with the same variance as in the CLT for functions of higher regularity, for any $\eps >0$.}
\end{tabular}
\end{table}

\tableofcontents

\section{Introduction}

This work is concerned with the fluctuations of linear spectral statistics of Wigner matrices $H$ for low regularity functions. Denote the eigenvalues of $H$ in increasing order by 
\[\lambda_1 \leq \lambda_2 \leq \dots \leq \lambda_N.\] 
Wigner's semicircle law states that for smooth $\varphi: \rr \to \rr$
\beq
\lim_{N \to \infty} \frac{1}{N} \tr  \varphi (H)  = \lim_{N \to \infty} \frac{1}{N} \sum_{i=1}^N \varphi ( \lambda_i ) = \int_{\rr} \varphi ( x) \rhosc (x) \d x,
\eeq
almost surely,
where
\beq
\rhosc (x) := \frac{1}{ 2 \pi } \sqrt{ (4 - x^2 )_+}.
\eeq
We will be interested in the asymptotic fluctuations of these linear spectral statistics. 
 The centered linear spectral statistic associated to $\varphi$ and $H$ is,
\beq
\N_N ( \varphi ) :=  \tr \varphi (H) - \ee[ \tr \varphi (H) ] .
\eeq
Generically, one expects that $\N_N ( \varphi )$ converges as $N \to \infty$ in distribution to a Gaussian random variable with variance
\begin{align} \label{eqn:intro-var}
V_H ( \varphi ) &:= \frac{1}{ 2 \pi^2} \int_{-2}^2 \int_{-2}^2 \frac{ ( \varphi (x) - \varphi (y) )^2}{ (x-y)^2} \frac{ 4 - xy}{ \sqrt{ 4-x^2} \sqrt{ 4 - y^2 } } \d x \d y \notag\\
&+ \frac{s_4}{ 2 \pi^2} \left( \int_{-2}^2 \varphi (x) \frac{ 2 - x^2}{ \sqrt{ 4 - x^2}} \d x \right)^2 ,
\end{align}
where $s_4$ is the fourth cumulant of the matrix entries, which is $0$ in the case of the Gaussian Orthogonal Ensemble (GOE).  Early works on the fluctuations of linear spectral statistics include the work of Arkharov \cite{arkharov1971limit} who announced the convergence of normalized traces of sample covariance matrices to Gaussian random variables, the work of Jonsson \cite{jonsson1982some} who gave a proof of this result,  and Girko \cite{girko2012theory} who studied the fluctuations of the traces of resolvents of Wigner and sample covariance matrices.  Since these works there has emerged a large literature proving central limit theorems for various classes of test functions and classes of random matrices. We refer the interested reader to, e.g.,  the works \cite{lytova2009central,he2017mesoscopic,guionnet2002large,cipolloni2020functional,li2021central,az,bai2008clt,chatterjee2009fluctuations,costin1995gaussian} and the references therein.

A common feature  of all these results is that the limiting variance contains a term whose most singular part is a double integral of
\beq \label{eqn:int-singular}
\frac{(\varphi(x)-\varphi(y))^2}{(x-y)^2} .
\eeq
This extends to models even beyond Wigner or covariance matrices, such as beta ensembles \cite{johansson1998fluctuations, bekerman2018mesoscopic}, discrete beta ensembles \cite{borodin2017gaussian}, as well as various determinantal point processes whose kernels are given by or approximate a spectral projection \cite{soshnikov2000gaussian, soshnikov2002gaussian, deleporte2021universality}. 

Of particular interest is the extent to which these results depend on the regularity of $\varphi$.  Johansson \cite{johansson1998fluctuations} obtained the central limit theorem for a class of unitarily invariant Hermitian matrices and general $\varphi \in H^{2+\eps} ( \rr)$. Here, we recall that the homogeneous and inhomogeneous Sobolev spaces, $\dot{H}^s (\rr) = \dot{H}^s  $  and  $H^s = H^s ( \rr)$ are the spaces of all functions such that the norms
 \beq
 \| \varphi \|^2_{\dot{H}^s} := \int | \hat{ \varphi} (\xi )|^2 | \xi|^{2s}  \d \xi  , \qquad \| \varphi \|^2_{H^s} := \int | \hat{ \varphi} (\xi )|^2 ( 1 + | \xi|^{2s} ) \d \xi 
 \eeq
 are finite, 
 where $\hat{\varphi}$ denotes the Fourier transform.  The significance of the expression \eqref{eqn:int-singular} in the context of these Sobolev spaces is that,
 \beq \label{eqn:int-2}
\| \varphi \|_{\dot{H}^{1/2} ( \rr ) } = c \int_{\rr^2} \frac{ ( \varphi (x) -\varphi (y) )^2}{ ( x -y )^2} \d x \d y ,
\eeq
for a numerical constant $c>0$.

Johansson further conjectured that the sole condition under which the CLT holds should be that the limiting expression for the variance is finite. Given the similarity of \eqref{eqn:int-2} with the most singular component of \eqref{eqn:intro-var}, it is natural to investigate this conjecture via test functions in the Sobolev spaces $H^s$, the goal being to get as close to $s=1/2$ as possible. The H{\"o}lder spaces $C^s = C^s (\rr)$ are also natural function spaces suited to this purpose.

 Shcherbina \cite{shcherbina2011central} obtained the CLT for the functions in the space $H^{3/2+\eps}$, and Lytova and Pastur obtained a CLT for the GOE and $\varphi \in C^1$ \cite{lytova2009central}.  The second author with Wong \cite{sosoewong} obtained various results for complex Hermitian Wigner matrices, including a CLT for the class $C^{1/2+\eps}$ under moment conditions, as well as an unconditional result for the space $H^{1+\eps}$, the latter being the sharpest available for general Wigner matrices. The sharpest result in terms of regularity is the result of Kopel who showed that in the case of the GUE, the CLT holds if and only if the expression for the limiting variance \eqref{eqn:intro-var} is finite \cite{kopel2015regularity}. 
 
On the other hand, consider $\varphi (x) = \1_{ \{ x \leq E \} }$ for $E \in (-2, 2)$. For the GUE, Gustavsson showed that the linear statistic $\mathcal{N}_N[\varphi]$ quantity has Gaussian fluctuations with variance of order $\log(N)$ \cite{gustavsson2005gaussian}. This result was extended to the GOE by O'Rourke as well as to Wigner matrices under the additional condition that the first four moments match those of the Gaussian ensemble \cite{orourke}. This was also studied by Dallaporta and Vu \cite{dallaporta2011note}. These results were extended to general Wigner matrices by the authors \cite{meso} and independently and in parallel by Bourgade and Mody \cite{bourgade2019gaussian}. Note that for the indicator function of a compact interval, the expression \eqref{eqn:intro-var} is infinite, but in terms of Fourier modes, the divergence is logarithmic:
\[\int_{|\xi|\le N} |\xi||\widehat{1}_I(\xi)|^2\,\mathrm{d}\xi \sim \log N.\]

For general Wigner matrices, a result as sharp as that of Kopel \cite{kopel2015regularity} seems beyond the reach of current methods. However, the aforementioned results \cite{meso,bourgade2019gaussian} on the indicator function give some evidence of the necessity of the limiting expression $V_H ( \varphi)$ being finite. In the present work we will show that, at least in the spectral bulk, finiteness of $V_H ( \varphi)$ is almost sufficient, for a general class of Wigner matrices with smooth distribution.

The question of the dependence of the fluctuations on the regularity of the function $\varphi$ has an essentially complete answer for linear statistics of the \emph{Circular} Unitary Ensemble (CUE) on $\mathbb{T}=[0,2\pi)$, whose joint density is proportional to
\begin{equation}\label{eqn: cue-density}
    \prod_{j<k}|e^{i\theta_j}-e^{i\theta_k}|^2\,\mathrm{d}\mathbf{\theta}.
\end{equation}
The normalized density \eqref{eqn: cue-density} is that of the eigenvalues of a unitary random matrix chosen from the Haar measure.

For the purposes of discussing the CUE, we let, for $k \in \mathbb{Z}$, 
\beq
\widehat{f}(k)=\frac{1}{2\pi}\int_{-\pi}^\pi \e^{-ik\theta}f(\e^{i\theta})\,\mathrm{d}\theta
\eeq
denote the $k$th Fourier coefficient of the function $f$ on the torus (we will not consider such functions outside of the introduction). 
Diaconis and Evans \cite{diaconis2001linear}, building on previous work by Diaconis and Shashahani \cite{diaconis1994eigenvalues}, showed that the  quantity
\beq
\mathcal{N}_{f,\mathrm{CUE}}:=\frac{\sum_{j=1}^N f(\e^{i\theta_j})-\frac{N}{2\pi}\int f(\e^{i\theta})\,\mathrm{d}\theta}{\sqrt{\sum_{-N}^N|k||\widehat{f}(k)|^2}}
\eeq
converges in distribution to normal random variable with mean zero and variance 1, provided the sequence $c_n=\sum_{-n}^n k^2|\widehat{f}(k)|^2$ is regularly varying. In particular, if $f$ belongs to the Sobolev $\dot{H}^{\frac{1}{2}}(\mathbb{T})$ of functions for which 
\beq
\|f\|_{\dot{H}^{1/2}(\mathbb{T})}^2:=\sum_{-\infty}^\infty |k| |\widehat{f}(k)|^2<\infty,
\eeq
 then $\sum_{j=1}^N f(\e^{\i\theta_j})-\frac{N}{2\pi}\int f(\e^{i\theta})\,\mathrm{d}\theta$ converges to a normally distributed random variable whose variance coincides with the $\dot{H}^{1/2}(\mathbb{T})$ norm of $f$. 
We note that the $\dot{H}^{1/2}(\mathbb{T})$ also has a physical space expression analogous to the first term in \eqref{eqn:intro-var} 
\beq
\|f\|^2_{\dot{H}^{1/2}}=c\int_{-\pi}^\pi\int_{-\pi}^\pi \frac{(f(\e^{i\theta})-f(\e^{i\nu}))^2}{\sin\left(\frac{\theta-\nu}{2}\right)^2}\,\mathrm{d}\theta\mathrm{d}\nu.
\eeq
Diaconis and Evans' result also subsumes both the CLT for smooth functions $f$, and the analog for the CUE of the result concerning indicator functions of intervals. When $f$ is such an indicator of $I\subset \mathbb{T}$, the sum $\sum_{-N}^N k^2|\widehat{f}(k)|^2$ diverges logarithmically.

As stated above, the optimal condition for convergence obtained by Kopel for the GUE is more complicated than being in $\dot{H}^{1/2}$, but if the support of $\varphi$ is compactly contained in $(-2, 2)$, then the main (universal) term in \eqref{eqn:intro-var} is bounded above and below by the expression,
\beq \label{eqn:int-1}
\int_{-2}^2 \int_{-2}^2 \frac{ ( \varphi (x) - \varphi (y) )^2}{ (x-y)^2} \d x \d y.
\eeq
Under the above support condition, the integral form of the $\dot{H}^{1/2}$ norm \eqref{eqn:int-2} is a norm-equivalent to the expression \eqref{eqn:int-1}, as well as the inhomogeneous norm $\| \cdot \|_{H^{1/2}}$.  We therefore turn to whether the CLT holds for functions in the space $H^{1/2+\eps} ( \rr)$.

The significance of the $H^{1/2}$ threshold can also be understood by looking at the special case where
\beq
\varphi(x)=\varphi_{E,\eta}(x):=\frac{1}{x-E-\i \eta},
\eeq
with $N$-dependent $\eta$. 
In this case, $\mathcal{N}_N(\varphi)$ is the trace of the resolvent matrix, and the \emph{local semicircle law} (see \eqref{eqn:ll} below), implies that, in a very strong sense, the fluctuations of $\mathcal{N}_N(\varphi_{E,\eta})$ are of order $O(1/\eta)$, up to an $N^\varepsilon$ error:
\begin{equation} \label{eqn: sc-fluct}
    \mathcal{N}_N(\varphi_{E,\eta})-\mathbb{E}[\mathcal{N}_N(\varphi_{E,\eta})] = \O(N^\varepsilon/\eta)
\end{equation}
with very high probability. This is the order of fluctuations predicted by the $H^{1/2}$ norm:
\beq
\|\varphi_{E,\eta}\|^2_{H^{1/2}}\approx \|\partial_x\varphi_{E,\eta}\|_{L^2}\|\varphi_{E,\eta}\|_{L^2}\approx \frac{1}{\eta^2}.
\eeq
In contrast, general purpose concentration tools such as log-Sobolev inequalities do not capture the correct order of fluctuations in the local regime for Wigner matrices because the bounds they provide are expressed in terms of the Lipschitz norms of the functions. 

Our main result, Theorem \ref{thm:main-1} below, is that for any $\eps >0$ there is a constant $C_\eps >0$ so that,
\beq
\Var ( \N_N ( \varphi ) ) \leq C_\eps \| \varphi \|_{H^{1/2+\eps}}^2 
\eeq
for $\varphi$ with support contained in $[-2+\eps, 2-\eps]$ and $H$ a Wigner matrix whose entries have a smooth distribution (the precise definition appears in the next section). The CLT for such $\varphi$ follows by  approximation by smooth functions and the existing literature on linear spectral statistics. 

Up to the $\eps >0$ factor, this is expected to be optimal as the CLT as stated does not hold for functions in $H^{1/2-\eps}$ as shown by the example of the indicator function.$^1${\let\thefootnote\relax\footnotetext{$1$. Given the results of Diaconis and Evans \cite{diaconis2001linear}, it is natural to conjecture is that, given any function $\varphi$, a CLT holds after renormalizing $\varphi$ by the $H^{1/2}$ norm of its projection onto its first $N$ Fourier modes. While it may be possible to prove such a result for the Gaussian cases, this seems out of reach of current techniques for general ensembles.} These are the strongest results available for general Wigner matrices. In fact, before our work, the CLT was not known for even the GOE case for functions $\varphi$ not admitting a classical derivative.

\subsection{Statement of main results}

We first define the matrix ensembles we investigate in the current work.

\bed[Real symmetric Wigner matrices] \label{def:rs-wig} A real symmetric \emph{Wigner matrix} $H$ is an $N \times N$ symmetric matrix such that the upper triangular part $\{ H_{ij} \}_{i \leq j }$ are independent centered random variables such that,
\beq \label{eqn:rs-var}
\ee[ H_{ij}^2] = \frac{ 1 + \delta_{ij}}{N}.
\eeq
We furthermore assume that the random variables $\{ H_{ij} \}_{1\le i < j\le N}$ are iid and that the random variables $\{ H_{ii} \}_{1\le i\le N}$ are iid. 
We furthermore assume that cumulants obey,
\beq \label{eqn:rs-cumu}
\kappa_k (H_{11} ) = 2^{1-k} \kappa_k (H_{12} )
\eeq
for $0 \leq k \leq 4$. Additionally, for any $k>0$ we assume that there is a $C_k >0$ so that,
\beq
\ee[ |H_{ij} |^k ] \leq \frac{C_k}{N^{k/2}}.
\eeq
\eed

The assumption \eqref{eqn:rs-cumu} is mainly for notational convenience in some resolvent expansions, and could be removed with further effort. Note that if $X$ is a (non-symmetric) matrix of iid centered random variables of variance $1/N$ then \eqref{eqn:rs-var} and \eqref{eqn:rs-cumu} is satisfied by $H = \frac{1}{ \sqrt{2}} (X + X^T)$.

Due to use of Wegner estimates as well as the reverse heat flow, we will require the densities of the matrix elements to have some smoothness. We formulate this in the following.

\bed[Smooth densities] \label{def:smooth-density}  We say a probability density $h(x)$ is Wigner-smooth if there is a $\delta >0$ so that
\beq \label{eqn:smooth-tail}
\int_{\rr} \exp\left[ \delta |x|^2 \right] h(x) \d x < \infty
\eeq
and for every $k$ there is a $C>0$ so that,
\beq \label{eqn:smooth-deriv}
 \left| \frac{\d^k}{ \d^k x} \log h(x) \right| \leq C ( 1+ |x| )^C.
\eeq
\eed

\bed[Smooth Wigner matrices] \label{def:rs-smooth-wig}
We say that a real symmetric Wigner matrix is  smooth if there are two Wigner-smooth probability densities $h_o(x)$ and $h_d(x)$ such that the off-diagonal entries are distributed according to $h_o$ and the diagonal entries are distributed according to $h_d$. 
\eed

The smoothness conditions above could be relaxed somewhat. For example, the sub-Gaussian tails \eqref{eqn:smooth-tail} assumption could be significantly relaxed. We include it so that we can directly  cite the Wegner estimates of \cite{wegner}. However, as noted in \cite{fixed-wig}, the relevant estimates hold under weaker conditions. For simplicity we just assume \eqref{eqn:smooth-tail}. The smoothness conditions \eqref{eqn:smooth-deriv} could also be relaxed, as in \cite{bourgade2018extreme}, in that the constants on the RHS could be assumed to diverge like $N^{c_1 k}$ for some $c_1 >0$.   The $c_1 >0$ would then taken to be small depending on the parameter $\eps >0$ in the assumption that $\varphi \in H^{1/2+\eps}$. Again, as our arguments are quite long we do not pursue this generalization. 

Results for discrete ensembles and $\varphi$ in the Sobolev spaces $H^{s}$ seem out of reach of current methods due to the use of Wegner estimates to cut-off very high-frequency modes. However, some results for discrete ensembles and $\varphi \in C^s$ for various $1/2 < s < 1$ can be proven. This requires further arguments and will be discussed in future work.



The main result of the paper is the following variance bound.

\bet \label{thm:main-1}
Let $H$ be a smooth real symmetric Wigner matrix as in Definition \ref{def:rs-smooth-wig}. For any $\eps >0$ and $\kappa >0$, there is a $C>0$ so that the following holds. For any $\varphi : \rr \to \rr$ supported in $(-2+\kappa, 2- \kappa)$ we have,
\beq
\Var \left( \tr ( \varphi (H) ) \right) \leq C_\eps \| \varphi \|^2_{H^{1/2+\eps}}.
\eeq
\eet

By approximating $\varphi$ by smooth functions one easily obtains the following.
\bec
Let $H$ be a smooth real symmetric Wigner matrix and $\varphi$ as in Theorem \ref{thm:main-1}. Then the centered linear spectral statistic $\N_N ( \varphi )$ converges to a Gaussian random variable with variance given by \eqref{eqn:intro-var}.
\eec

Note that under our support conditions, it follows from the integral representation of the $\dot{H}^s$ norm that 
\beq
\| \varphi \|_{H^{1/2+\eps}} \leq C_\delta \| \varphi \|_{C^{1/2+\eps+\delta}}
\eeq
and so the above results hold also for $ \varphi \in C^{1/2+\eps}$ for any $\eps >0$.

Our methods also allow one to calculate the usual corrections to the expectation (i.e., the next order terms after $\int \varphi (x) \rhosc (x) \d x$) for $\varphi \in H^{1/2+\eps}$. We summarize this in the following remark.
\begin{rmk} Let $H$ and $\varphi$ be as above. The methods of our paper can be used to derive the expansion,
\begin{align}
\ee[ \tr \varphi (H) ] - N \int \varphi (x) \rhosc (x) \d x =& - \frac{1}{ 2 \pi } \int_{-2}^2 \varphi (x) \frac{1}{ \sqrt{4-x^2}} \d x \notag\\
 +& \frac{s_4}{ 2 \pi } \int_{-2}^2 \varphi (x) \frac{ x^4 - 4 x^2 +2}{ \sqrt{ 4 - x}} \d x +  o (1)
\end{align}
where $s_4$ is the fourth cumulant of the off-diagonal entries.

This is technically easier than the estimate on the variance and so we do not provide the details. Note that usually there is an additional term on the RHS involving $\varphi ( \pm 2)$. However, by the compact support condition of $\varphi$ this term is absent.  
\end{rmk}

Finally, we also remark that the condition that $\varphi$ vanishes near the spectral edges is not strictly necessary under the condition that $\varphi$ is smooth near the edges. For example, one could obtain a CLT if $\varphi = \varphi_1 + \varphi_2$ where $\varphi_1$ is as in the above results and $\varphi_2$ is sufficiently smooth, but is allowed to be non-zero near the edges $\pm 2$.


\subsection{Methodology} \label{sec:methods}

As our work is rather involved, we now summarize our methods and give a detailed proof sketch.

The starting point of our methods is a modification of the Littlewood-Paley decomposition introduced in the work of the second author with Wong in \cite{sosoewong}. Roughly, this decomposes $\varphi$ into the sum,
\beq \label{eqn:int-lp}
\varphi = \sum_{k \geq -1} \varphi_k = \sum_{k} (P_{2^{-k} } ) \star g_k
\eeq
where $P_\eta$ denotes the Poisson kernel
\beq
P_\eta (x) := \frac{1}{ \pi} \frac{\eta}{ x^2 + \eta^2}
\eeq
and $\varphi_k$ and $g_k$ whose Fourier support is contained in frequencies only of order $2^k$. The quantities,
\beq
\sum_k 2^{2ks} \| \varphi_k \|_2^2, \qquad \sum_k 2^{2ks} \| g_k \|_2^2
\eeq
are then comparable to $H^{s}$ norm of $\varphi$. The main advantage of the last representation on the right side of \eqref{eqn:int-lp} is that when taking the trace against $H$, this becomes a convolution against
\beq
\tr \, P_\eta (H-E) = \frac{1}{ \pi } \Im \, \tr \frac{1}{ H - z}, \qquad z = E + \i \eta ,
\eeq
i.e., a function of the resolvent $G(z) = (H-z)^{-1}$. The resolvent is highly studied in random matrix theory  (see for example \cite{benaych2019lectures} and \cite[Chapter 6]{erdHos2017dynamical}) and there are many techniques available to estimate it. 

The decomposition \eqref{eqn:int-lp} further allows us to deal with different ranges of frequencies with different methods. In particular, the behavior of $\varphi$ near the frequency $2^k$ involves only the behaviour of the eigenvalues on the scale $\eta := 2^{-k}$.  We now discuss how the various frequency ranges are handled.

Due to the embedding of $H^{1/2+\mfa}$ into $L^\infty$ for any $\mfa >0$ and the trivial estimate $\Var ( \tr f(H) ) \leq N^2 \| f\|_\infty^2$, the frequencies larger than $2^k \geq N^C$ for some large $C>0$ depending on $\mfa >0$ are easily dispensed with.  Moreover, the projection of $\varphi$ onto frequencies $2^k \leq N^c$ for small $c>0$ yields a function whose $C^5$ norm grows less than $N^{c'}$ for small $c'>0$ depending on $c>0$. Such a function is sufficiently smooth so that it can be handled by the existing resolvent techniques for proving CLTs in random matrix theory.

This leaves estimating the variance of the logarithmically many $\varphi_k$ for $N^{c} \leq 2^k \leq N^C$. By Cauchy-Schwarz, it suffices to prove the estimate
\beq
\Var ( \tr ( \varphi_k (H) ) ) \leq C N^{\eps} 2^k\left(  \| \varphi_k \|_2^2 + \| g_k \|_2^2 \right)
\eeq
for any $\eps >0$ and 
for $k$ in this range.  These estimates are the main technical contribution of the current work.

 The first simplification is that, by the reverse heat flow \cite{local-relax} and our smoothness assumptions, we can assume that $H$ is a Gaussian divisible ensemble of the form,
\beq
H_t := \e^{ -t/2} W + \sqrt{ 1 - \e^{-t} } G
\eeq
for any $t = o (1)$, $W$ a smooth Wigner matrix and $G$ a matrix from the GOE.  This allows us to use Dyson Brownian motion techniques to analyze the eigenvalue behavior.  In particular, we will use the approach to DBM using the stochastic advection equation developed in the work of Bourgade \cite{bourgade2018extreme}. These results play a fundamental role in our approach as the work \cite{bourgade2018extreme} provides essentially optimal estimates which turn out to be crucial; weaker estimates would not suffice.

Here, we split these estimates into three regimes, using different techniques in each range. The three regimes are:
\begin{enumerate}[label=(\arabic*)]
\item \label{it:meso} The mesoscopic and microscopic regimes $N^c \leq 2^k \leq C_1 N$,
\item \label{it:micro} The submicroscopic regime $C_1 N \leq 2^k \leq N^{2-\eps_1}$,
\item \label{it:ultra} The very high frequency regime $N^{2-\eps_1} \leq 2^k \leq N^C$.
\end{enumerate} 
Above, $\eps_1 >0$ is chosen sufficiently small depending on the $\mfa>0$ such that $\varphi \in H^{1/2+\mfa}$. 
The tightest threshold is for frequences $2^k \approx C_1 N$: the techniques we use in the regime \ref{it:meso} work up to $k$ satisfying $2^k \leq C_1 N$ for any $C_1 >0$ but will not work for $2^k \approx N^{1+\eps}$. On the other hand, the techniques we use in the regime \ref{it:micro} work only for $2^k \geq C_1 N$ for some sufficiently large $C_1 >0$. 

The regime \ref{it:ultra} is easiest. The desired variance bound follows by Cauchy-Schwarz and a Wegner estimate \cite{wegner},
\beq
\sup_{E \in \rr} \ee[ | \lambda_i - E - \i \eta |^{-1} ] \leq C' N^{1+\eps}
\eeq
for any $\eta > N^{-C}$ and $\eps >0$. In fact, this argument does not depend on the Gaussian divisibility of $H_t$. 

The treatment of the remaining regimes \ref{it:meso} and \ref{it:micro} relies heavily on the Gaussian divisibility and Dyson Brownian Motion (DBM). One main input is the work of Bourgade \cite{bourgade2018extreme} that provides an almost-optimal approach to DBM. Roughly, this work exhibits a coupling between the eigenvalues of $H_t$ (denoted by $\lambda_i$) and the eigenvalues of the GOE (denoted by $\mu_i$) such that,
\beq \label{eqn:int-dbm}
\lambda_i = \mu_i + N^{-1}(\Phi_t (W) - \Phi_t (G')) + N^{\eps} \O \left( \frac{1}{ N^2 t} \right),
\eeq
where $\Phi_t (A) := \tr f_t (A)$ and $f_t$ is an approximate step function on the scale $t$ centered at $\gamma_i$, $i$th $N$-quantile of $\rhosc$, $W$ is the smooth Wigner matrix in the definition of $H_t$ and $G'$ is an auxilliary GOE matrix. In particular, $W$ is independent of $\mu_i$ and $G'$.

If $t$ is close to $1$, the estimate \eqref{eqn:int-dbm} allows one to access the behavior of $\lambda_i$ down to the scale $N^{-2}$. Given that the reductions outlined above reduce our proof to the regime $2^k \leq N^{2 - \eps}$ it makes sense that Bourgade's result gives some hope of proving our desired estimates. We now discuss the implementation of these ideas. It turns out to be rather different in the regimes \ref{it:meso} and \ref{it:micro}, but we first make some general remarks applicable to both cases.

The estimate \eqref{eqn:int-dbm} says that $\lambda_i$ is given by a universal random variable, $\mu_i - \Phi_t ( G')$ up to the quantity $\Phi_t (W)$. If one can show that the behavior of $\Phi_t (W)$ is independent of the entry distribution, or otherwise control it, then one can relate the behavior of $\lambda_i$ directly to the universal behavior. However, this approach cannot be implemented in such a simple minded fashion. There are many highly non-trivial obstacles:
\begin{enumerate}[label=(\roman*)]
\item The behavior of $\Phi_t(W)$ is in fact not universal. Asymptotically, it is Gaussian but the variance is non-universal as it depends on $s_4$ through the formula \eqref{eqn:intro-var}
\item Since we need to control $\lambda_i$ down to the scale $N^{-2+\eps_1}$, we would need to calculate, for example, 
\beq
\ee[ \e^{ \i \lambda \Phi_t (W) } ]
\eeq
for $\lambda $ up to $N^{1-\eps_1}$. Resolvent techniques are available to study this quantity, but we expect that such fine estimates would require significant technical effort and numerous iterated expansions.  Moreover, one also requires estimates for large $\lambda$ with extremely small errors as well.
\item Even if one can show that $\Var ( \tr \varphi_k (H_t ) ) \approx \Var ( \tr \varphi_k (G ) )$ for $G$ from the GOE, this still leaves one the task of calculating the quantity on the RHS. To our knowledge, the required estimates for the GOE are lacking from the literature (estimates for the GUE seem to be available).
\item It is difficult to even get started with the use of \eqref{eqn:int-dbm} to replace $\lambda_i$ by the $\mu_i$ in $\Var ( \tr \varphi_k (H_t ) )$ with an acceptable error. This in fact occupies a significant component of our work.
\end{enumerate}
Let us discuss how these obstacles are overcome. They appear in different guises in each of the regimes \ref{it:meso} and \ref{it:micro}. First, recalling $g_k$ above, we have:
\begin{align} \label{eqn:int-var-rep}
\Var ( \tr ( \varphi_k (H_t ) ) &= \frac{1}{ \pi^2} \int_{\rr^2} g_k (x) g_k (y) \Cov \left( \Im[ \tr G(x + \i \eta ) ] , \Im [ \tr G ( y + \i \eta ) ] \right) \d x \d y \notag\\
&= \sum_{i, j}  \frac{1}{ \pi^2} \int_{\rr^2} g_k (x) g_k (y) \Cov \left( \Im [ ( \lambda_i - x - \i \eta )^{-1} , \Im [ ( \lambda_j - y - \i \eta )^{-1} \right) \d x \d y.
\end{align}
Note that there is a well-known estimate,
\beq
\Var (  \Im[ \tr G(x + \i \eta )] ) \leq N^{\eps } \eta^{-2}
\eeq
available for any $\eps >0$ and large ranges of $\eta >0$. Using this in \eqref{eqn:int-var-rep} would yield $\Var ( \tr ( \varphi_k (H) ) ) \leq N^{\eps} \eta^{-2} \| \varphi_k \|_2^2$. This is too lossy, as one needs an estimate with a factor $\eta^{-1}$ instead of $\eta^{-2}$.

\paragraph{Regime \ref{it:micro}.}  We first discuss the submicroscopic regime $2^k \gg N$. 
In the introduction, we mentioned the conjecture that the CLT should hold for any $\varphi$ after renormalizing by the projection $\sum_{k\le N}\varphi_k$ of $\varphi$ onto its first $N$ Fourier modes. We would then hope that the frequences with $2^k \gg N$ do not contribute to the asymptotic fluctuations. This is supported by the following heuristic. 

We will focus on the contribution of $\ee[ ( \tr \varphi_k (H) )^2]$ to the variance, as this will explain our approach. Consider the representation on the last line of \eqref{eqn:int-var-rep} and apply the estimate \eqref{eqn:int-dbm} twice, once each for $\lambda_i$ and $\lambda_j$, yielding two mesoscopic linear statistics of $W$ which we denote by $\Xi_i$ and $\Xi_j$. Assuming that the $\lambda_i$ can be replaced by the expression on the right side at an acceptable error, one would find by Fourier duality,
\begin{align} \label{eqn:int-micro-a1}
 & \int_{\rr^2} g_k (x ) g_k (y)  \ee\left[ \Im [ ( \lambda_i - x - \i \eta )^{-1}]  \Im [ ( \lambda_j - y - \i \eta )^{-1}] \right] \d x \d y \notag\\ 
\approx & \int_{\rr^2} g_k (x ) g_k (y)  \ee\left[ \Im [ ( \hat{\mu}_i - x - N^{-1}\Xi_i - \i \eta )^{-1}]  \Im [ ( \hat{\mu}_j - y - N^{-1}\Xi_j- \i \eta )^{-1}] \right] \d x \d y  \notag\\
= & \int_{\rr^2} \hat{g}_k ( \xi_1 ) \hat{g}_k ( \xi_2 ) \hat{F} ( \xi_1, \xi_2 ) \ee[ \e^{ \i ( \xi_1 N^{-1} \Xi_1 + \xi_2 N^{-1} \Xi_2 )] } ] \d \xi_1 \d \xi_2
\end{align}
where 
\beq
F(x, y) = \ee \left[ \Im [ ( \hat{\mu}_i - x  - \i \eta )^{-1}]  \Im [ ( \hat{\mu}_j - y - \i \eta )^{-1}] \right]
\eeq
 and $\hat{\mu}_i$ involves only $\mu_i$ and $G'$ and is independent of $W$.  It turns out that the $\Xi_i$ are asymptotically Gaussian with variance growing as $\log(N)$ and so one expects that $\ee[ \e^{ \i \xi_i N^{-1} \Xi_i } ] \leq \e^{ - c \xi^2 \log(N) N^{-2}}$. For $|\xi_i | \geq C_1 N$ for $C_1$ sufficiently large, this decay is sufficient to completely remove any contribution of the $g_k$'s to the overall variance, because the $g_k$'s are supported only in frequences of order $2^k \geq C_1 N$ in the regime \ref{it:micro}. 
 
 However, it is difficult to implement this directly. First, in fact the $\Xi_i$ are highly correlated if $t$ is of the same order as $\gamma_i - \gamma_j$ (here, we denote the $N$-quantiles or classical locations of the semicircle distribution by $\gamma_i$, which are the typical locations of the eigenvalues $\lambda_i$). This will destroy the decay of the joint Fourier transform of $(\Xi_i, \Xi_j)$. Therefore, one sees that one must apply \eqref{eqn:int-dbm} using a time $t_{ij}$ adapted to the distance $\gamma_i - \gamma_j$, i.e., $t_{ij} := N^{-c_1} | \gamma_i - \gamma_j|$, for some $c_1>0$.  Note that then if $i$ and $j$ are quite close, $t_{ij}$ is quite short. That the error incurred by taking this $t_{ij}$ in the first approximation in \eqref{eqn:int-micro-a1} is acceptable, of order $\O ( \eta^{-1} \| \varphi_k \|_2^2$), seems to be a minor miracle in itself (here, it is required to take $c_1 >0$ adapted to the $\mfa >0$ in the assumption $\varphi \in H^{1/2+\mfa}$). Establishing this requires significant effort and delicate estimates. This is the content of Proposition \ref{prop:sub-2} below.
 
 Still, this approach then requires handling the Fourier transform of $\Xi_i$ for large frequencies, as well as obtaining very precise estimates with a small error term. Since $\Xi_i$ is a mesoscopic spectral statistic, this could in principle be handled by resolvent expansions. Such an approach may be possible, but it would require significant effort and resolvent expansions to arbitrarily high order. Even for global statistics where $t_{ij} \approx 1$, obtaining optimal estimates would require significant effort. See, e.g. \cite{bao2021quantitative}, where strong estimates are obtained, but would still not be sufficient for our purposes.
 
 Instead, we observe that we require only \emph{decay} of the joint Fourier transform of the $(\Xi_i, \Xi_j )$ and not a precise asymptotic estimate. We therefore apply Dyson Brownian Motion again, and realize $W$ itself as a Gaussian divisible ensemble. The form of $\Xi_i$ is that of an $\arctan(\cdot)$ function, or $\Im[ \log(\cdot)]$. It turns out that the SDE obeyed by $\tr \log (H-z)$ under DBM is relatively simple and amenable to study using the method of characteristics as in \cite{huang2019rigidity}.  This is carried out in Section \ref{sec:char}. 
 
 This allows for a decomposition of $\Xi_i \approx Z_i + A_i$ where $Z_i$ is an explicit Gaussian random variable and $A_i$ depends on the initial data and is independent of $Z_i$. The error in this approximation turns out to be the same order as that incurred in \eqref{eqn:int-dbm} and so we make this approximation before applying the Fourier duality. As $(Z_i, Z_j)$ are jointly Gaussian and weakly correlated, we find the desired decay. This part of the argument appears in Proposition \ref{prop:sub-3} below.
 
\paragraph{Regime \ref{it:meso}} We now discuss the regime \ref{it:meso}.  In this case, we use the representation on the right side of the first line of \eqref{eqn:int-var-rep}. Let us first discuss what is expected. Based on the form of the covariance between two resolvents centered at $x$ and $y$ for large $\eta$, we expect that the estimate,
\beq \label{eqn:int-cov}
\left|  \Cov \left( \Im[ \tr G(x + \i \eta ) ] , \Im [ \tr G ( y + \i \eta ) ] \right) \right| \leq \frac{C}{\eta^2 + (x-y)^2} ,
\eeq
holds. 
Note that this bound (or a weaker one losing some small polynomial in $N$ factors) would imply the desired result as by using a simple Schwarz inequality $|g_k (x) g_k (y) | \leq \frac{|g_k (x) |^2 + |g_k (y) |^2}{2}$ one easily obtains the desired estimate of $C\eta^{-1} \| g_k \|_2^2$ for $\Var ( \tr ( \varphi_k (H) ) )$.

Let us discuss possible approaches to obtaining the estimate \eqref{eqn:int-cov}. First, in the regime $2^k \leq N^{1-\eps}$, or $\eta \geq N^{\eps-1}$, resolvent expansions are stable and in principle the estimate could be obtained in this fashion. Indeed, the paper of Knowles and He \cite{he2020mesoscopic} obtains an expansion for the covariance kernel for $x$ and $y$ in the bulk. However, this still leaves aside the regime $2^k \approx N$. The method we use in this regime turns out to easily apply to the entire regime $N^{\eps} \leq 2^k \leq N$ and so we proceed in this fashion.

Our method is to use a modification of the estimate \eqref{eqn:int-dbm} in which the family of mesoscopic linear spectral statistics centered at $\gamma_i$ for each $i$ is replaced by two statistics, one centered at $x$ and one centered at $y$. This estimate is applied to the second line of the representation \eqref{eqn:int-micro-a1}. This modified estimate has an additional error term on the right side of \eqref{eqn:int-dbm} adapted to the distance $|x-\gamma_i|$. Since the function $\Im[ \tr ( G (x + \i \eta ) )]$ is sharply peaked at $x$, this error turns out to be acceptable.  Here, we choose one $t$ for all $x$ and $y$ which is taken to be large, almost equal to $1$.  This part of the argument appears as Proposition \ref{prop:meso-homog-1} below.

In the regime under consideration, decay from the Fourier transform of the $\Phi_t$ cannot help us, as the fluctuations of $\Phi_t$ are essentially on the scale $\eta$, and so they instead must be computed directly. As specified above, the hope is that, if we can show that the fluctuations of $\Phi_t$ are universal, then the covariance kernel can be compared to that of the Gaussian ensembles. There are numerous obstacles to this, as mentioned above.

The first issue is calculating the characteristic function of the (almost-global) mesoscopic linear spectral statistic. Here, we must determine it up  to an error of $\O (N^{-1})$ in order for the error to be acceptable for $2^k \approx N$. To our knowledge, such an estimate does not exist in the literature. In fact, the naive guess that,
\beq \label{eqn:int-false}
\ee[ \e^{ \i \xi ( \Phi_t (W) - \ee[ \Phi_t (W) ] ) } ] = \e^{ - \xi^2 V_W ( f_t ) } + \O (N^{-1} )
\eeq
where $V_W$ is the functional \eqref{eqn:intro-var} turns out to be false.  It turns out that there is an additional correction to $V_W$ of order $N^{-1/2}$ depending on $s_3$ and moreover a non-Gaussian (multiplicative) correction to the characteristic function of the form $\e^{\i \xi^3 N^{-1/2} s_3 a}$ for some coefficient $a$ depending on the test function. 

That an estimate as strong as \eqref{eqn:int-false} cannot hold was observed by Bao and He \cite{bao2021quantitative}. In this work they considered global linear spectral statistics and subtracted off random correction terms which allowed for an estimate of $\O (N^{-1})$. This work does not appear applicable to our setting, as we need to calculate the characteristic function for $\Phi_t$ itself, without any correction terms. Moreover, they handle only global statistics; ours is almost-global but it is not stated how the estimates of \cite{bao2021quantitative} depend on the scale of the test function. We expect that it would enter in a mild fashion, but nonetheless this work does not appear directly applicable.

In the work \cite{meso} we provided an approach to calculating characteristic functions of linear spectral statistics based on cumulant expansions and Stein's method. In the present work we return to this method and simply calculate a few further terms that appear in the cumulant/resolvent expansions. We find the desired estimate of order $\O (N^{-1})$ after including more non-universal correction terms. This is the content of Section \ref{sec:fine-lss} (which uses methods of Section \ref{sec:global}).

At this point, we must figure out how to handle the non-universal corrections. In fact, one term, the $s_4$ term in the formula \eqref{eqn:intro-var}, is in fact of size $\O (N^{-1})$, the same scale as one extreme of the range considered in the regime \ref{it:meso}. 

Here, our main observation is that the non-universal terms can, in some sense, be viewed as an additional random shift to the locations $x$ and $y$ in the covariance function \eqref{eqn:int-cov}.  We implement this observation by de-convolving the function $g_k$ with a Gaussian of variance $A_1 N^{-1}$ for some large $A_1 >0$, and consequently convolving the centering $x$ of $\Im [ \tr (G (x + \i \eta ) ) ]$ with this Gaussian. Due to the fact that $2^k \leq C_1 N$, this deconvolution action is bounded on $g_k$ in the $L^2$ sense. Here, we see the necessity of the restriction of the upper end of the range \ref{it:meso}. It is possible that this could be extended to something like $2^k \leq c_1 \sqrt{ \log(N) } N$ for some $c_1 >0$ adapted to the regularity scale $\mfa >0$ but this is not necessary for us.

The convolution by a Gaussian of the centering $x$  in the covariance function allows us to absorb the non-universal corrections. Essentially, this shows that the covariance of two resolvents for a Gaussian divisible ensemble can be controlled by the covariance of two resolvents of the GOE/GUE after a random shift is of order $\O (N^{-1} )$.  This is essentially the content of Propositions \ref{prop:meso-a1} and \ref{prop:meso-a2} below. It follows that if we can show that the estimate \eqref{eqn:int-cov} holds for the Gaussian ensembles, we will obtain it also  for our matrix $H$, and complete the proof; this is due to the fact that the random shift is of smaller or equal order than $\eta$, so introducing the shift on the RHS of \eqref{eqn:int-cov} does not affect the resulting estimate.

To our knowledge the estimate \eqref{eqn:int-cov} does not appear in the literature for the Gaussian ensembles down to the scales that we require, even for the GUE. One contribution of our work is to establish such an estimate in both the GUE and GOE cases. We establish estimates for the covariance function using the well-known representations for the correlation functions of the Gaussian ensembles in terms of the Hermite polynomials, see, e.g., \cite{mehta2004random}. Here, a significant simplification is afforded by the fact that we are restricted to the spectral bulk in which the asymptotics for the Hermite polynomials are straightforward. In fact, the GUE case is dispensed with quite quickly by relying on these asymptotics and the Christoffel-Darbox formula.

The GOE case is more involved and must be split into the $N$ odd and $N$ even cases.  Here, the cluster functions are essentially ``perturbations'' of the GUE kernel. There does not seem to be a unified approach to these additional terms and they must instead be checked on a case-by-case basis. For the GOE case, we in fact do not obtain an estimate as strong as \eqref{eqn:int-cov} and obtain on the right side a function of the form $\O ( \eta^{-1} (|x-y| + \eta)^{-1})$ which turns out to be sufficient for our purposes. We expect that this is significantly suboptimal and that \eqref{eqn:int-cov} should hold for the GOE, but we do not need this explicit form and do not attempt to prove it. In our estimates, it turns out that the $N$-odd case is more challenging; here we observe a non-obvious (to us, at least) cancellation between two terms that is only evident after some asymptotic analysis of some integrals of Hermite polynomials. We suspect that this has been observed before in the literature but did not find it after a cursory search.  The estimates for the Gaussian ensembles appear in Section \ref{sec:gaussians}.

 \subsubsection{Further discussion}

Some of the above methods are tangentially related to the recent work of Cipolloni, Erd{\H{o}}s and Schr{\"o}der \cite{cipolloni2021spectral} on the spectral form factor defined as,
\beq
S(t) := \frac{1}{N^2} \sum_{i, j=1}^N \e^{ \i t ( \lambda_i - \lambda_j ) }.
\eeq
In the work \cite{cipolloni2021spectral}, they derive asymptotics for the expectation and variance of this quantity in the regimes $ t \ll N^{5/11}$ and $t \ll N^{5/17}$, respectively, using extensive resolvent expansions. Via Fourier duality, one could relate $\Var ( \tr ( \varphi_k (H) ))$ to $S(t)$. However, our work is primarily focused on larger frequency regimes corresponding to $t \geq N$  and consequently there is no overlap in methods or results with those of \cite{cipolloni2021spectral}.

\subsection{Organization of paper}

In Section \ref{sec:background} we collect the results from the random matrix literature that we will use in our work. This includes the local and isotropic semicircle laws, Wegner estimates as in \cite{wegner}, the homogenization work of Bourgade \cite{bourgade2018extreme}, the Littlewood-Paley theory applied to random matrix theory as in \cite{sosoewong}, as well as some standard results on the Gaussian ensembles. In Section \ref{sec:technical} we collect the main technical estimates proven in this work. In Section \ref{sec:sub-micro}, we establish the desired estimates on $\Var ( \tr \varphi_k (H))$ for the submicroscopic frequences $2^k \geq C_1 N$ (i.e., the regime \ref{it:micro} as discussed in Section \ref{sec:methods}). In Section \ref{sec:meso} we deal with the frequences in the range $N^{\eps} \leq 2^k \leq C_1 N$, i.e., the regime \ref{it:meso} as in Section \ref{sec:methods}.  In Section \ref{sec:global} we establish estimates for the covariance of resolvents for large frequencies including near the edge. In Section \ref{sec:char}  we apply the method of characteristics to study the mesoscopic linear spectral statistic that arises in the application of Bourgade's homogenization result. In Section \ref{sec:gaussians} we apply the exact formulas for the cluster functions of the Gaussian ensembles to derive estimates on the covariance of two resolvents with small imaginary part for the GOE and GUE. Section \ref{sec:main} contains the proof of Theorem \ref{thm:main-1}. In Section \ref{sec:fine-lss} we establish fine estimates on the characteristic function of almost-global linear spectral statistics via resolvent expansions. The appendices collect various auxilliary results, including a few calculations involving the semicircle distribution, as well as some standard Hermite polynomial asymptotics and an asymptotic calculation of some integrals of Hermite polynomials that arise in the $N$-odd GOE case of Section \ref{sec:gaussians}. 

\subsection{On the symmetry class of the ensemble}

Throughout this paper we will restrict ourselves to the real symmetric Wigner matrices. We will require a few results that are based on the exact calculations for the Gaussian ensembles using the formulas for the correlation functions or cluster functions, which are  different for the GOE and GUE. Where we cite results from the literature, we state them for both the GOE/GUE. Our primary use of exact formulas is in our estimates for the covariance of two resolvents. There, we give the proofs for both ensembles.

We state and prove results on the characteristic functions of almost-global linear statistics only in the real symmetric class. The complex Hermitian case can be treated by the exact same method and requires only additional notation, as it based only on resolvent expansions.

For the rest of the paper, the symmetry class of the matrix ensemble under consideration does not make a significant difference to the methods.

\subsection{Acknowledgements}

The work of Philippe Sosoe is partially supported by NSF grant DMS-1811093.

\section{Background} \label{sec:background}


In this section we collect the various results of random matrix theory that we will use in our proofs. This section also serves to collect the various notation that we use.

We will use the following notion of high-probability events.
\bed
We say that a family of events $\A_i$ indexed by $i \in \I$ some index set hold with overwhelming probability if for all $D>0$ there is a $C>0$ so that
\beq
\sup_{i \in \I} \pp [ \A^c_i ] \leq C N^{-D},
\eeq
for all $N$.
\eed
If $a_i$ and $b_i$ are two positive quantities indexed by some set $\I$ (for example, $N$-dependent sequences, or positive functions on the upper half-plane) then the notation
\beq
a_i \asymp b_i
\eeq
means that there is a $C>0$ so that
\beq
\frac{1}{C} a_i \leq b_i \leq C a_i ,
\eeq
for all $i$.

\subsection{Local semicircle laws and rigidity}

Let $H$ be a Wigner matrix with eigenvalues $\lambda_1 \leq \lambda_2 \leq \dots \leq \lambda_N$. We introduce the following notation,
\beq
m_N (z) := \frac{1}{N} \sum_{i=1}^N \frac{1}{ \lambda_i - z}
\eeq
and
\beq
\msc (z) := \int \frac{1}{ x -z } \rhosc (x) \d x = \frac{ -z + \sqrt{ z^2 -4 }}{2}
\eeq
where Wigner's semicircle density is given by
\beq
\rhosc (x) := \frac{ \sqrt{4- x^2}}{ 2 \pi } \1_{\{ |x| \leq 2 \} }.
\eeq
Note that,
\beq
\msc(z)^2 + z \msc (z) + 1 = 0.
\eeq
We will denote the resolvent by,
\beq
G_{ij} (z) = \left( \frac{1}{H -z } \right)_{ij}
\eeq
For any $\tau >0$ define,
\beq \label{eqn:Dtau-def}
\D_{\tau, 1} := \{ z = E + \i \eta : |E| \leq \tau, N^{\tau-1} <  \eta \leq \tau^{-1} \}.
\eeq
and
\beq
\D_{\tau, 2} := \{ z = E + \i \eta : |E| \leq \tau,0 <  \eta \leq \tau^{-1} \}.
\eeq

\bet
Let $H$ be a Wigner matrix. Let $\tau >0$ and $\eps >0$. Then, with overwhelming probabiltity we have for all $z \in \D_{\tau, 2}$ that,
\beq \label{eqn:ll}
|m_N (z) - \msc (z) | \leq \frac{ N^{\eps}}{N \Im[z] }
\eeq
as well as,
\beq \label{eqn:entry-ll}
| G_{ij} (z) - \delta_{ij} (z) | \leq N^{\eps} \left( \sqrt{ \frac{ \Im[\msc (z) ] }{ N \Im[z] } } + \frac{1}{ N \Im[z] } \right)
\eeq
for all $i,j$.
\eet
\bet
Let $H$ be a Wigner matrix with resolvent $G(z)$. Let $\tau >0$, $\eps >0$ and $D>0$. There is a $C>0$ so that,
\beq \label{eqn:iso}
\pp\left[ \exists z \in \D_{\tau, 1} : \left| \bx^T G (z) \by - \bx^T \by \msc(z) \right|  >  N^{\eps} \left( \sqrt{ \frac{ \Im[\msc (z) ] }{ N \Im[z] } } + \frac{1}{ N \Im[z] } \right)\right] \leq C N^{-D}
\eeq
for all unit vectors $\bx, \by \in \cc^N$.
\eet
From the above results, we see that it will be useful to introduce the control parameter,
\beq \label{eqn:Psi-def}
\Psi (z) = \sqrt{ \frac{  \Im [ \msc (z) ] }{ N \Im [z]} } + \frac{1}{N | \Im[z] | }.
\eeq

We will denote the quantiles of the semicircle distribution by $\gamma_i$. They are defined by,
\beq
\frac{i}{N} = \int_{-2}^{\gamma_i } \rhosc (x) \d x ,
\eeq
where $\gamma_N = 2$. 

\bet Let $H$ be a Wigner matrix with eigenvalues $\lambda_1 \leq \lambda_2 \leq \dots \leq \lambda_N$. 
Let $\eps >0$. With overwhelming probability,
\beq
| \lambda_i - \gamma_i | \leq \frac{N^{\eps}}{N^{2/3} \min\{ i^{1/3}, (N+1-i)^{1/3} \} },
\eeq
for every $i$.
\eet

For any $z \in \cc$ we define the distance to the spectral edge as,
\beq
\kappa (z) := ||\Re[z]| - 2 | | .
\eeq
For the Stieltjes transform of the semicircle distribution we have the following, see for example \cite{benaych2019lectures}, \cite{erdHos2017dynamical}.
\bel
Let $ z= E + \i \eta$, for $\eta >0$. If $|E| < 2$ then,
\beq
\Im [ \msc (z) ] \asymp \sqrt{ \kappa(z) + \eta }
\eeq
and if $|E| > 2$ then,
\beq
\Im [ \msc(z) ] \asymp \frac{ \eta}{ \sqrt{ \kappa(z) + \eta } }.
\eeq
\eel

\subsection{Homogenization of Dyson Brownian motion} \label{sec:homog}

A Gaussian divisible ensemble is a Wigner matrix $H$ such that $H$ can be written as independent sum of
\beq
H  = \e^{-t/2} W + \sqrt{1 - \e^{-t} } G
\eeq
where $t >0$, $W$ is another real symmetric Wigner matrix and $G$ is from the GOE. Note that $H$ satisfies \eqref{eqn:rs-cumu} iff $W$ does. 

Dyson Brownian motion is the solution $\{x_i (t) \}_{1 \leq i \leq N, t \geq 0}$ to the following system of stochastic differential equations,
\beq \label{eqn:dbm-def}
\d x_i (t) = \sqrt{\frac{2}{N  } } \d B_i (t) + \frac{1}{N} \sum_{ j \neq i } \frac{1}{ x_i (t) - x_j (t) } \d t - \frac{x_i (t) }{2} \d t,
\eeq
where $\{ B_i (t) \}_i$ is an iid family of standard Brownian motions. 
 If the initial data $\{ x_i (0) \}_{i=1}^N$ are distributed as the eigenvalues of a symmetric matrix $H_0$ with real entries then $x_i(t)$ has the same distribution as the eigenvalues of the Gaussian divisible ensemble
\beq
H_t := \e^{-t/2} H_0 + \sqrt{ 1- \e^{-t}} G
\eeq
where $G$ is a GOE matrix independent from $H_0$.

Associated with DBM are the characteristics:
\begin{align} \label{eqn:char-def}
z_t &:= \frac{ \e^{t/2} ( z + \sqrt{ z^2 -4 }  ) + \e^{-t/2} (z - \sqrt{ z^2-4} )}{2} \notag\\
&= \cosh(t/2) z + \sinh(t/2) \sqrt{z^2-4}.
\end{align}
A straightforward calculation shows that
\beq \label{eqn:delz-bwd}
\del_t z_t = \msc (z_t) + \frac{z_t}{2} .
\eeq
The characteristic is well defined for $\Im[z] >0$. For $E \in \rr$ we define,
\beq \label{eqn:char-real-line-def}
E_t = \lim_{ \eta \dto 0 } (z+ \i \eta )_t.
\eeq

Consider now two DBM process $\{ \mu_i (t) \}_{i,t}$ and $\{ \lambda_i (t) \}_{i, t}$ evolving according to \eqref{eqn:dbm-def}, with the same Brownian motion terms. For the initial data, we assume that $\lambda_i (0)$ are the eigenvalues of a Wigner matrix $W$ and the $\mu_i (0)$ are the eigenvalues of a Gaussian ensemble of the same symmetry class as $W$. 
 Define,
\beq \label{eqn:XW-def}
X^W(z, t) := \frac{ \Im \sum_{k} \log ( \lambda_k (0) - z_t ) }{ \Im \msc (z_t ) } - N \frac{ \int \Im \log (x -z_t ) \rhosc (x) \d x }{ \Im \msc (z_t ) }
\eeq
and
\beq \label{eqn:XG-def}
X^G(z, t) := \frac{ \Im \sum_{k} \log ( \mu_k (0) - z_t ) }{ \Im \msc (z_t ) } - N \frac{ \int \Im \log (x -z_t ) \rhosc (x) \d x }{ \Im \msc (z_t ) }
\eeq
 The following Theorem is a reformulation of Theorem 3.1 of \cite{bourgade2018extreme}.  The exact statement first appeared in \cite{blz} (but see also \cite{bourgade2016fixed} for a similar formulation where a solution to a discrete parabolic equation coming from homogenization for DBM was rewritten in terms of a mesoscopic linear spectral statistic).  It follows in a straightforward manner from Theorem 3.1 of \cite{bourgade2018extreme} as well as the proof of Lemma 3.4 that appears there.  This reformulation is useful as the functions $X^W$ and $X^G$ are mesoscopic linear spectral statistics which may be studied via resolvent and dynamical methods.
\bet \label{thm:sae-homog}
Let $\{ \mu_i (t) \}_i$ and $\{ \lambda_i (t) \}_i$ be as above. Let $\eps >0$, $\delta >0$ and $\kappa >0$.  For any $D>0$ there exists a constant $C>0$ so that for any $k \in \ldbrack 1,N \rdbrack$ and $E$ such that  $E, \gamma_k \in (-2+\kappa, 2 - \kappa)$ we have for any $ N^{\delta-1} < t < 1$ that
\beq
\pp \left[ \left| \lambda_k (t) - \mu_k (t) - \frac{X^W (E, t) - X^G (E, t) }{N} \right| > N^{\eps} \frac{ N|E- \gamma_k | + 1 }{N^2 t}  \right] \leq CN^{-D}.
\eeq
\eet
\proof Theorem 3.1 of \cite{bourgade2018extreme} gives an estimate for $\lambda_k (t) - \mu_k (t)$ in terms of a function $\bar{u}_k (t)$. We then use (3.3) of \cite{bourgade2018extreme} to replace it by $(N \Im [ \msc ( z_t ) ] )^{-1} \int_0^1 \Im[f_0 (z_t) ] \d \nu$ where $f_0 (z_t)$ is defined in (1.11) of \cite{bourgade2018extreme} and contains a parameter $\nu$ omitted from the notation. The statement of Lemma 3.4 of \cite{bourgade2018extreme} contains a restriction $\eta \gg N^{-1}$ in $z_t = (E+ \i \eta )_t$.  However one can see from the proof that in our case $s=t$ the estimate nonetheless holds with $| \eta + t - s | = 0$. Now, by the definition of $f_0$ one can see that the integral over $\nu$ equals our $X^W - X^G$. \qed

\subsection{Littlewood-Paley theory} \label{sec:lpt}

In this section we review elements of the Littlewood-Paley theory. This was first applied to random matrix theory in the work \cite{sosoewong}, and our treatment is basically identical to the one given there. We will denote the Fourier transform of $f : \rr \to \cc$ by
\beq
\hat{f} ( \xi ) := \int_{\rr} \e^{  \i \xi x } f(x) \d x.
\eeq
 Given $\varphi \in \S$ where $\S$ is the space of Schwarz distributions, the inhomogeneous Littlewood-Paley decomposition is the sum,
\beq
\varphi = \sum_{k \geq -1 } \varphi_k
\eeq
where,
\beq
\varphi_{-1} = h \star \varphi
\eeq
and 
\beq
\varphi_k = 2^k \omega ( 2^k \cdot ) \star \varphi , \qquad k \geq 0
\eeq
where $h$ and $\omega$ are smooth functions satisfying the following. The Fourier transforms $\hat{h}$ and $\hat{\omega}$ are even, smooth and of compact support. The function $\hat{h}$ is supported in the ball $\{ | \xi | \leq 4/3 \}$ and the function $\hat{\omega}$ is supported in the annulus $\C := \{ 3/4 \leq | \xi | \leq 8/3 \}$. Moreover, $ 0 \leq \hat{h}, \hat{\omega} \leq 1$ and,
\beq
\hat{h} ( \xi ) + \sum_{k \geq 0 } \hat{\omega} ( 2^{-k} \xi ) = 1 .
\eeq
We have the following, see, e.g., \cite{bahouri2011fourier}.
\bet
For any $s >0$ there is a $C>0$ so that,
\beq
\frac{1}{C} \| \varphi \|_{H^s}^2 \leq \sum_k 2^{2 k s } \| \varphi_k \|_{L^2}^2 \leq C \| \varphi \|^2_{H^s} 
\eeq
and for any $0 < \alpha < 1$ there is a $C>0$ so that,
\beq
\frac{1}{C} \| \varphi \|_{C^\alpha } \leq \sup_{k} 2^{ \alpha k } \| \varphi_k \|_{L^\infty} \leq C \| \varphi \|_{C^{\alpha}}
\eeq
\eet
Following \cite{sosoewong}, the Littlewood-Paley decomposition allows for a convenient decomposition of $\varphi$ into convolutions with Poisson kernels at different scales. Define,
\beq
P_\eta (x) := \frac{1}{ \pi} \frac{ \eta}{ \eta^2 + x^2}.
\eeq
Recall,
\beq
\hat{P}_\eta (x) = \e^{ - \eta | \xi | }.
\eeq
Define now,
\beq
\hat{g}_k (\xi ) = \e^{ 2^{-k} | \xi | } \hat{\varphi}_k ( \xi).
\eeq
Then, the following is Theorem 5 of \cite{sosoewong}.
\bet \label{thm:poisson-lp}
 Let $s \geq 0$ and $0 < \beta < \alpha < 1$. If $\varphi \in H^s$ or $\varphi \in C^\alpha$, then $\varphi$ has a representation,
 \beq
 \varphi (x) = \sum_{k=-1}^\infty ( P_{2^{-k}} \star g_k ) (x)
 \eeq
 where the sum on the right is convergent in $H^s$ if $\varphi \in H^s$ and in $C^\beta$ if $\varphi \in C^\alpha$.  There are constants $C_s$ and $C_{\alpha, \beta}$ such that,
 \beq
 \sum_{k=-1}^\infty 2^{2 k s } \| g_k \|_{L^2}^2 \leq C_s \| \varphi \|_{H^s}^2
 \eeq
 and
 \beq
 \sum_{k=-1}^\infty 2^{\beta k} \| g_k \|_{L^\infty} \leq C_{\alpha, \beta} \| \varphi \|_{C^\alpha}.
 \eeq
\eet

In general we will assume that $\varphi$ is of compact support. However, the various $g_k$ and $\varphi_k$ will not be of compact support. The following controls their decay. 
\bep \label{prop:lp-support}
Let $\varphi$ be in $H^s$ or $C^\alpha$ and be supported in $(-2+\kappa, 2 - \kappa)$. Let $|x| > 2 - \kappa$.  Then, for every $M>0$ there is a $C_M >0$ so that for any $k \geq 0$,
\beq
|\varphi_k (x) | \leq C_M \frac{2^k}{1 +2^{kM} ||x|-(2-\kappa ) |^M} \| \varphi \|_{L^2}.
\eeq
The same estimate holds for $g_k(x)$.
\eep
\proof We prove the estimate for $g_k$, the estimate for $\varphi_k$ being simpler.  Let $\psi(x)$ be defined by,
\beq
\hat{\psi} (\xi ) := \hat{\omega} ( \xi) \e^{ | \xi | }.
\eeq
Since $\hat{\omega}$ vanishes near the origin, this is a smooth function of compact support. Therefore, for any $M >0$ there is a $C_M >0$ such that,
\beq
| \psi (x) | \leq \frac{C_M}{ (1 + |x| )^M}.
\eeq
On the other hand,
\beq
g_k (x) = 2^k \int \psi (2^k (x - y)) \varphi (y) \d y.
\eeq
For $y \in (-2 + \kappa, 2 - \kappa)$ and $|x| > 2 - \kappa$ we have,
\beq
|\psi (2^k (x - y)) | \leq \frac{C_M}{ 1 + 2^{Mk} | |x|- (2 - \kappa ) |^M }
\eeq
and the claim follows using $\| \varphi \|_{L^1} \leq 2 \| \varphi \|_{L^2}$ since $\varphi$ is of compact support. \qed

%

Let $\chi$ be a smooth even bump function such that $\chi (x) =1$ for $|x| \leq 1$ and $\chi (x) = 0$ for $|x| > 2$. We have the following.
\bep
Let $\varphi \in H^s$ or $C^\alpha$. Let $\lambda \geq 1$ and define $f$ by
\beq
\hat{f} ( \xi ) = \chi ( \lambda^{-1} \xi ) \hat{\varphi} ( \xi).
\eeq
Then for any $k \geq 1$ there is a $C_k >0$ such that,
\beq
\| f^{(k)} \|_{L^\infty} \leq C_k  \lambda^{k+1/2} \| \varphi \|_{L^2}.
\eeq
Moreover,
\beq
\| f \|_{H^s} \leq C \| \varphi \|_{H^s}, \qquad \| f \|_{C^\alpha } \leq C \| \varphi \|_{C^{\alpha }}
\eeq
\eep
\proof From Lemma 2.1 of \cite{bahouri2011fourier} we see,
\beq
\| f^{(k)} \|_{L^\infty} \leq C_k \lambda^{k+1/2} \| f \|_{L^2}
\eeq
Note that if $\psi(x)$ is the inverse Fourier transform of $\chi$ then
\beq
f = ( \lambda \psi (\lambda \cdot ) ) \star \varphi ,
\eeq
and so the first estimate follows from Young's inequality for convolutions. The last inequality of the Proposition follows also from Young's inequality, as well as the fact that if $f_k$ is the $k$th Littlewood-Paley projection for $f$ then,
\beq
f_k = ( \lambda \psi (\lambda \cdot ) ) \star \varphi_k .
\eeq
This completes the proof. \qed

\subsection{Wegner-type estimates}

We recall the following from \cite{erdHos2010wegner}. 
\bet \label{thm:weg}
Let $H$ be a smooth Wigner matrix. Let $\kappa >0$. Let $1/N > \eta > 0$.  There exists a constant $C>0$ so that for all $N \geq 10$ we have,
\beq
\pp \left[ \exists i : | \lambda_i - E | < \eta \right] \leq C N \eta
\eeq
for all $E \in (-2+\kappa, 2- \kappa)$.
\eet
\remark This is stated as Theorem 3.4 of \cite{erdHos2010wegner}. There, the estimate is given for complex Hermitian Wigner matrices, but as stated in the introduction of the paper, the method and estimate above applies equally well to the real symmetric case. 

 Note that we do not require such a sharp estimate, i.e., we could allow a weaker estimate with an $N^{\eps}$ factor on the RHS. The primary role in the sub-Gaussian assumption \eqref{eqn:smooth-tail} is so that \cite{erdHos2010wegner} can be directly cited; in that paper, the primary role of this role is so that the above estimate is obtained with no $N^{\eps}$ or $\log(N)^C$ factors. We do not need such a tight estimate for our methods. \qed

\bep \label{prop:weg-1}
Let $H$ be a smooth Wigner matrix and let $A>0, \kappa >0$ and $ \eps >0$. There exists a $C >0$ so that
\beq
\ee \left[ \frac{1}{N | \lambda_i - E - \i \eta | } \right] \leq C N^{\eps}.
\eeq
holds for all $E \in \rr$ and $\eta \geq N^{-A}$ and $ \kappa N \leq i \leq (1- \kappa )N$ .
\eep
\proof If $|E - \gamma_i | > N^{\eps-1}$ then the result holds by rigidity.  If $\eta \geq N^{-1}$ then the result is trivial. So we may assume that $E \in (-2 + \kappa', 2 - \kappa')$ for some $\kappa' >0$ and $ \eta \leq N^{-1}$. Then, for
\beq
K := \lceil (2 \eps -1) \log_2 (N) - \log_2 ( \eta ) \rceil
\eeq
we have,
\begin{align}
\ee \left[ \frac{1}{N | \lambda_i - E - \i \eta | } \right] &= \ee \left[  \frac{1}{N | \lambda_i - E - \i \eta | }  \1_{ \{ | \lambda_i - E| < \eta \} } \right] \notag \\
&+ \sum_{k=1}^K  \ee \left[  \frac{1}{N | \lambda_i - E - \i \eta | }  \1_{ \{ 2^{k-1} \eta <| \lambda_i - E| < 2^{k}\eta \} } \right]  \notag \\
&+ \ee \left[  \frac{1}{N | \lambda_i - E - \i \eta | }  \1_{ \{ | \lambda_i - E| > 2^{K} \eta \} } \right] .
\end{align}
By the choice of $K$ and rigidity we have that the last line is $\O (1)$ as,
\beq
\pp \left[ | \lambda_i - E| > 2^{K} \eta  \right] \leq N^{-D}
\eeq
for any $D>0$. For the first two terms on the RHS we have,
\begin{align}
  & \ee \left[  \frac{1}{N | \lambda_i - E - \i \eta | }  \1_{ \{ | \lambda_i - E| < \eta \} } \right]  + \sum_{k=1}^K  \ee \left[  \frac{1}{N | \lambda_i - E - \i \eta | }  \1_{ \{ 2^{k-1} \eta <| \lambda_i - E| < 2^{k}\eta \} } \right]  \notag \\
 \leq & \frac{C}{N} \sum_{k=0}^K \frac{1}{ \eta 2^k } \pp \left[ \exists i : | \lambda_i - E | < 2^k \eta \right] \leq C K ,
\end{align}
where the last line follows from Theorem \ref{thm:weg}. This yields the claim. \qed

\bep \label{prop:weg-2}
Let $H$ be a smooth Wigner matrix. Let $A>0$, $\kappa >0$ and $\eps >0$.  There exists a $C>0$ so that,
\beq
\sup_{E \in \rr} \sum_{ \kappa N \leq i \leq (1- \kappa)N } \ee\left[ \frac{1}{ N | \lambda_i - E - \i \eta | } \right] \leq C N^{\eps}. 
\eeq
\eep
\proof For any $E$, there are at most $C N^{\eps}$ indices such that $| \gamma_i -E | \leq N^{\eps-1}$. The contribution from these indices is estimated via Proposition \ref{prop:weg-1}. The contribution from the remaining indices is estimated via rigidity and is bounded above by,
\beq
\sum_{ \kappa N \leq i \leq (1- \kappa)N } \1_{ \{ |\gamma_i - E | > N^{\eps-1}  \} } \ee\left[ \frac{1}{ N | \lambda_i - E - \i \eta | } \right] \leq C \sum_{i=1}^N \frac{1}{ i} \leq C \log(N)
\eeq
and the claim follows. \qed

\subsection{Density of states of Gaussian ensembles}

The following results on the asymptotic density of states in the bulk for the Gaussian ensembles can be found in  \cite{garoni2005asymptotic, forrester2006asymptotic}. 
\bet \label{thm:gaussian-bulk}
Let $\kappa >0$.  Uniformly for $E \in (-2+\kappa, 2- \kappa)$, the density of states of the GUE satisfies,
\beq
\rho (E) = \rhosc (E) + \frac{ (-1)^N}{ 4 \pi^3 N \rhosc(E)^2} \cos \left[ N \left( \frac{E}{2} \sqrt{ 4 - E^2} + 2 \arcsin \frac{E}{2} \right) \right] + \O ( N^{-2} )
\eeq
and the density of states of the GOE satisfies,
\beq
\rho (E) = \rhosc (E) - \frac{1}{ 4 \pi^2 N \rhosc (E) } + \O (N^{-2} ).
\eeq
\eet
In Section \ref{sec:gauss-expect} we derive the following which is useful at the edge.
\bep \label{prop:gaussian-edge}
Let $\varphi : \rr \to \rr$ and $\eps >0$. There exists $C >0$ so that, for the GOE and GUE we have,
\beq
\left| \ee[ \tr ( \varphi (H) ) ] - N \int \varphi (x) \rhosc (x) - \left( 1- \frac{2}{ \beta} \right) \int_{-2}^2 \frac{ \varphi(x)}{ \sqrt{4-x^2} } \d x \right|  \leq C N^{\eps} \| \varphi \|_{C^5}.
\eeq
where $\beta=1$ in the GOE case and $\beta=2$ for GUE.
\eep

\section{Main technical estimates} \label{sec:technical}

In this section, we collect the various technical estimates that we prove throughout the paper. 

We have the following asymptotics of the covariance of the imaginary part of the empirical Stieltjes transform which is useful on scales $\gg N^{-1/2}$. 
\bet \label{thm:global-cov}
Let $H$ be a Wigner matrix and $\eps >0$. For any $ z = x \pm \i y$ and $w = u + \pm \i v$ with $u, v \geq N^{-1/2}$ we have,
\begin{align}
& \left| N^2 \Cov \left( \Im [ m_N (z) ], \Im [ m_N (w) ] \right) \right| \leq \frac{C}{ \sqrt{ \kappa(x) + y} \sqrt{ \kappa(u) + v } } \frac{1}{ (x-u)^2 + y^2 + v^2} \notag\\
+ & \frac{ N^{\eps}}{ \sqrt{ \kappa(x) + y } + \sqrt{ \kappa(u) + v } } \left( \frac{1}{N y^2 v} + \frac{1}{ N v^2 y} \right).
\end{align}
\eet
This theorem is proven in Section \ref{sec:cov-ests}. 
Introduce the norm,
\beq \label{eqn:weighted-norm}
\| \varphi \|_{1, w} := \int_{\rr} | \varphi(x)| |4 -x^2|^{-1/2} \d x.
\eeq
For a Wigner matrix $H$ with third and fourth cumulants $s_3$ and $s_4$, we introduce a number of functionals which will be used to describe the fluctuations of linear statistics.  Let
\begin{align} \label{eqn:main-V-def}
V( \varphi ) &= \frac{1}{ 2 \pi^2} \int_{-2}^2 \int_{-2}^2 \frac{ (\varphi(x) - \varphi(u) )^2}{ (x-u)^2} \frac{ 4- xu}{ \sqrt{ 4 - x^2} \sqrt{ 4 - u^2} } \d x \d u \notag\\
&+ \frac{s_4}{2 \pi^2} \left( \int_{-2}^2 \varphi (x) \frac{ 2 -x^2}{ \sqrt{4-x^2} } \d x \right)^2 \notag\\
&- \frac{2 s_3}{  \pi^2 N^{1/2}} \left( \int_{-2}^2 \varphi (x) \frac{ 2 -x^2}{ \sqrt{4-x^2} } \d x \right) \left( \int_{-2}^2 \varphi (x) \frac{x}{ \sqrt{4-x^2} } \d x \right)
\end{align}
and
\begin{align}
B(f) =- \frac{1}{ 12 \pi^3} \left( \int_{-2}^2  f(x) \frac{x}{\sqrt{4-x^2}} \d x \right)^3,
\end{align}
and
\begin{align} \label{eqn:main-expect-correct}
e(f) &=  e_1(f) + e_2 (f)
\end{align}
where
\begin{align}
e_1 (f) &=- \frac{1}{ 2 \pi }\int f(x) \frac{1}{\sqrt{4-x^2} } \d x + \frac{ f(2) + f(-2) }{4} 
\end{align}
and
\begin{align}
e_2 (f) &= \frac{s_4}{2 \pi} \int_{-2}^2 f(x) \frac{ x^4 - 4 x^2 + 2}{\sqrt{4-x^2}} \d x + \frac{2 s_3}{N^{1/2}} \frac{1}{ \pi} \int_{-2}^2 f(x) \frac{ 3x -  x^3}{ \sqrt{4 -x^2}} 
\end{align}
The same resolvent expansions used to prove the above result will yield the following. 
\bet \label{thm:wig-smooth-lss}
Let $H$ be a Wigner matrix and let $\varphi : \rr \to \rr$ be given, and assume $\| \varphi'' \|_{1, w} \leq N$ and that the support of $\varphi$ is in $[-5, 5]$. For all $\eps >0$ there is a $C>0$ so that,
\beq
\left| \Var ( \tr \varphi (H) ) - V(\varphi )  \right| \leq CN^{-1+\eps}  \| \varphi''\|_{1, w},
\eeq
and
\beq
\left| \ee[ \tr \varphi (H) ] - N \int \varphi (x) \rhosc (x) \d x - e(\varphi) \right| \leq CN^{-1+\eps}  \| \varphi''\|_{1, w},
\eeq
and for all $\xi \in \rr$ we have,
\begin{align}
&\left| \ee\left[ \exp \left( \i \xi ( \tr \varphi (H) - N \int \varphi (x) \rhosc (x) - e ( \varphi ) ) \right) \right] - \exp \left( - \xi^2 V(f)/2 + \i \xi^3 s_3 N^{-1/2} B(f) \right) \right|  \notag\\
\leq & CN^{\eps} (1 + |\xi|^6+ \| \varphi' \|_{1, w}^7 ) N^{-1} \| \varphi''\|_{1, w}
\end{align}
\eet
\remark If $\varphi$ is not of compact support, then the above estimates still hold with additional error terms.  Let $A_1 := \| \varphi \vert_{[-4, 4]^c} \|_\infty$.  For any $D>0$, the variance estimate holds with an additional error of order $(1+A_1)^2 N^{-D}$, the expectation estimate with an additional error of order $A_1 N^{-D}$ and the estimate for the characteristic function with an additional error of order $(1+ | \xi | )(1+A_1)N^{-D}$. \qed

Theorem \ref{thm:wig-smooth-lss} is proven in Section \ref{sec:wig-smooth-lss-proof}. 
Consider now $\varphi$ and its Littlewood-Paley decomposition,
\beq
\varphi = \sum_k \varphi_k.
\eeq
We will derive estimates for $k$ in various ranges.

\bet \label{thm:meso-var}
Let $H$ be a Wigner matrix of the form,
\beq
H = \e^{-T/2}W + \sqrt{ 1 - \e^{-T} } G
\eeq
where $G$ is an independent GOE matrix.  Assume that,
\beq
T = N^{- \omega}, \qquad 0 < \omega < \frac{1}{100}.
\eeq
Fix a small constant $\mfb$ satisfying,
\beq
0 < \mfb < \frac{1}{100}.
\eeq
Let $k$ satisfy,
\beq
\mfb \log_2 (N) \leq k \leq \log_2 (N \mfb^{-1} ).
\eeq
Then for all $\eps >0$ there is a $C>0$ (depending also on $\mfb$ and $\omega >0$) so that,
\beq
\Var ( \tr \varphi_k (H) ) \leq C N^{\eps} T^{-1} \| \varphi_k\|_2^2 2^{k} + C N^{-100} \| \varphi \|_2^2
\eeq

\eet
This theorem is proven at the end of Section \ref{sec:meso}.

\bet \label{thm:sub}
Let $T = N^{-\omega}$, for $0 < \omega < \frac{1}{10}$. Let $H$ be a Wigner matrix of the form
\beq
H = \e^{-T/2} W + \sqrt{ 1 - \e^{-T} } G
\eeq
where $G$ is an independent Gaussian matrix and $W$ is a smooth Wigner matrix. Let $\eps >0$. Let $\delta > \omega$.  Let $A>0$ and $M>0$. There is a $B_1 >0$ depending on $\delta - \omega$ and a $C>0$ depending on all of these parameters so that for $k$ satisfying,
\beq
\log_2 (B_1 N ) \leq k \leq 2 \log_2 (N)
\eeq
we have,
\beq
\Var ( \tr ( \varphi_k (H) ) ) \leq C \left(  N^{\delta} N^{\eps} 2^{k} \| \varphi_k \|_2^2 + 2^{-Mk } \| \varphi \|_2^2 \right)
\eeq
Moreover, for $k$ satisfying
\beq
2 \log_2(N) \leq k \leq A \log_2 (N)
\eeq
we have for any smooth Wigner matrix,
\beq
\Var ( \tr ( \varphi_k (H) ) ) \leq C  \left( N^{\eps} 2^{k} \| \varphi_k \|_2^2 + 2^{-Mk } \| \varphi \|_2^2 \right)
\eeq
for any $\eps >0$.
\eet
The above theorem is proven in Section \ref{sec:sub-proof}. 

\bet \label{thm:main-uv}
Let $\varphi \in H^{1/2+s}$ for some $s>0$. Let $s> \eps >0$. There is a $C>0$ so that for any $M \geq 10^4$ we have,
\beq
\Var \left( \sum_{k\geq M} \tr \varphi_k (H) \right) \leq C  N^2 2^{M(\eps-2s)} \| \varphi \|_{H^{1/2+s}}^2
\eeq
\eet

\proof Let $0 < \rho < s$. If $\varphi$ is a function whose Fourier transform is supported only on frequencies larger than $100 \times  2^{L}$ for $L \geq 10$, then for $\rho > s$ we have,
\beq
\| \varphi\|_{H^{1/2+\rho}}^2 \leq C \sum_{k=L}^\infty \| \varphi_k \|_2^2 2^{2k(1/2+ \rho) } \leq C 2^{2L(\rho -s ) } \| \varphi\|_{H^{1/2+s}}^2.
\eeq
By the Sobolev embedding theorem (e.g., Theorem 1.66 of \cite{bahouri2011fourier}) we have,
\beq
\| f \|_\infty \leq C_\rho \| f \|_{H^{1/2+\rho} },
\eeq
for any $\rho >0$. 
We have the trivial inequality,
\beq
\Var ( \tr f(H) ) \leq N^2 \| f\|_{\infty}^2 .
\eeq
The claim follows by taking $L= M-100$.  \qed

\bet \label{thm:main-gaussian}
Let $G$ be a GOE or GUE matrix and let $\kappa >0$. Let $\eps >0$ be small and $\mfc >0$ be large. There is a constant $C>0$ so that for all $|x|, |y| < 2 - \kappa$, and $ 1 > \eta > N^{-\mfc}$
\beq
N^2 \Cov \left( \Im [ m_N (x + \i \eta ) ] , \Im [ m_N (y + \i \eta ) ] \right) \leq C N^{\eps} \left( \frac{1}{ \eta ( |x-y| + \eta ) } \right).
\eeq
\eet
\proof This follows from Propositions \ref{prop:gue-cov-calc}, \ref{prop:goe-even} and \ref{prop:goe-odd}.  \qed

\section{Sub-microscopic frequencies} \label{sec:sub-micro}

Let $H$  be a Wigner matrix. 
Let $\mfs >0$. We consider $\varphi \in H^{1/2+ \mfs}$ or $C^{1/2+\mfs}$. Recall the definition of $\varphi_k$ and $g_k$ as in Section \ref{sec:lpt}.  In this section we choose a large $B_1 >0$ and $A>0$ and consider $k$ in the range,
\beq
\log_2 (B_1 N) \leq k \leq A \log_2 (N).
\eeq
We will assume $B_1 \geq 1$ for definiteness. Constants in this section will be allowed to depend on $A$ and $B_1$, but the dependence on $B_1$ will be stated explicitly, as its role will become apparent only later in the section.


We consider the representation
\begin{align}
 & \Var ( \tr ( \varphi_k (H) ) ) \notag\\
= & \sum_{i, j} \int g_k (x) g_k (y)\, \Cov ( \Im [ ( \lambda_i - x - \i \eta )^{-1} ], \Im [ ( \lambda_j - x - \i \eta )^{-1} ] ) \d x \d y ,
\end{align}
where, throughout this section,
\beq
\eta := 2^{-k}.
\eeq
Due to Proposition \ref{prop:lp-support} we will see that it will suffice to consider pairs such that both $\gamma_i$ and $\gamma_j$ are in the spectral bulk. To facilitate this, introduce  the notation,
\beq \label{eqn:J-def}
J_\alpha := [\![ \alpha N, (1- \alpha)N ]\!],
\eeq
for an $\alpha >0$. 
Fixing a $\delta_1 >0$ we treat separately the contribution from indices $i, j$ such that $|i-j| \leq N^{\delta_1} / (N \eta)$ or otherwise. The contribution from the indices that are close together is estimated in the following Proposition.
\bep \label{prop:sub-diag-cut-off}
Let $H$ be a smooth Wigner matrix, and let $\delta_1 >0$ and $\eps >0$. Then, for
\beq
\log_2 (N) \leq k \leq A \log_2 (N)
\eeq
we have,
\begin{align}
 & \sum_{( i, j) \in J_\alpha^2} \1_{ \{ |i-j| < \frac{ N^{\delta_1}}{N \eta} \} } \int_{\rr^2} |g_k(x) g_k (y) | \ee[ \Im[ ( \lambda_i - x - \i \eta  )^{-1}] \Im[ ( \lambda_j - y - \i \eta  )^{-1} ] ] \d x \d y \notag\\
\leq & C \frac{ N^{\delta_1+\eps}}{\eta} \| g_k \|_2^2
\end{align}
and
\begin{align}
 & \sum_{( i, j) \in J_\alpha^2} \1_{ \{ |i-j| < \frac{ N^{\delta_1}}{N \eta} \} } \int_{\rr^2} |g_k(x) g_k (y) | \ee  \big[ \Im[ ( \lambda_i - x - \i \eta  )^{-1}]  \big]\ee\big[ \Im[ ( \lambda_j - y - \i \eta  )^{-1} ] \big] \d x \d y \notag\\
\leq &C \frac{ N^{\delta_1+\eps}}{\eta} \| g_k \|_2^2 
\end{align}
\eep
\proof We prove only the first estimate, the second being similar. By estimating $|g_k (x) g_k (y) | \leq |g_k(x)|^2 + |g_k (y) |^2$ and symmetry it suffices to estimate,
\begin{align}
 & \sum_{( i, j) \in J_\alpha^2} \1_{ \{ |i-j| < \frac{ N^{\delta_1}}{N \eta} \} } \int_{\rr^2} |g_k(x)|^2 \ee\big[ \Im[ ( \lambda_i - x - \i \eta  )^{-1}] \Im[ ( \lambda_j - y - \i \eta  )^{-1} ] \big] \d x \d y \notag\\ 
\leq & C \sum_{( i, j) \in J_\alpha^2} \1_{ \{ |i-j| < \frac{ N^{\delta_1}}{N \eta} \} } \int_{\rr} |g_k(x)|^2 \ee\big[ \Im[ ( \lambda_i - x - \i \eta  )^{-1}]  \big] \d x \notag\\
\leq &  C \frac{ N^{\delta_1}}{N \eta}  \sum_{i \in J_\alpha} \int_{\rr} |g_k(x)|^2 \ee\big[ \Im[ ( \lambda_i - x - \i \eta  )^{-1}]  \big] \d x .
\end{align}
In the first inequality we used the fact that,
\beq
\int_{\rr} \frac{ \eta}{ (s-y)^2 + \eta^2 } \d y \leq C
\eeq
for any $s$, and in the second inequality we did the sum over $j$, using that the quantity no longer depends on $j$.  Now by Proposition \ref{prop:weg-2} we have,
\beq
\frac{1}{N} \sum_{i \in J_\alpha } \ee[ \Im[ ( \lambda_i - x - \i \eta  )^{-1} ] ] \leq CN^{\eps},
\eeq
for all $x$.  The estimate with the norm $\| g_k \|_2^2$ follows after integration over the $x$ variable.\qed


The following deals with the contribution from the edge.
\bep \label{prop:sub-edge} Let $\eps >0$ and $M>0$.
Assume $\varphi$ is supported in $[-2+\kappa, 2 - \kappa]$ and let $k$ satisfy,
\beq
\log_2(N) \leq k \leq A \log_2 (N).
\eeq
Assume that $\alpha >0$ is so small that $\gamma_{\alpha N} < -2 + \kappa/4$. Then,
\begin{align}
& \left| \sum_{ (i, j) \notin J_\alpha^2 } \int_{\rr^2}  g_k (x) g_k (y)   \Cov \left( \Im [ ( \lambda_i - x - \i \eta )^{-1} ] , \Im [ ( \lambda_j - y - \i \eta )^{-1} ] \right) \d x \d y \right| \notag\\
\leq & C N^{\eps} \left( \eta^{-1} \| g_k \|_2^2 + 2^{-Mk } N^2  \| \varphi \|_2^2 \right).
\end{align}
\eep
\proof By Proposition \ref{prop:lp-support} it suffices to consider the integral over $D :=[-2+\kappa/2,2- \kappa/2 ]^2$. Clearly,
\begin{align}
& \left| \sum_{ (i, j) \notin J_\alpha^2 } \int_{D}  g_k (x) g_k (y)   \Cov \left( \Im [ ( \lambda_i - x - \i \eta )^{-1} ] , \Im [ ( \lambda_j - y - \i \eta )^{-1} ] \right) \d x \d y \right| \notag\\ 
\leq & 2  \left|\int_{D}  g_k (x) g_k (y)   \Cov \left( N \Im [ m_N (x + \i \eta ) ] ,  \sum_{j \notin J_\alpha}  \Im [ ( \lambda_j - y - \i \eta )^{-1} ] \right) \d x \d y \right| \notag \\
+ & 4 \int_D |g_k (x) g_k (y) |  \Var \left( \sum_{i \notin J_\alpha} \Im[ ( \lambda_i - x - \i \eta)^{-1} ] \right)  \d x \d y
\end{align}
By the rigidity estimates and the fact that if $(x, y) \in D$, then $x$ is separated from $\gamma_i$ for $i \notin J_\alpha$ by a constant independent of $N$, we easily see that,
\beq
 \Var \left( \sum_{i \notin J_\alpha} \Im[ ( \lambda_i - x - \i \eta)^{-1} ] \right)  \leq N^{\eps}
\eeq
for any $\eps >0$ and $N$ large enough.  By the local law, $\Var ( N \Im [ m_N ( x + \i \eta ) ] ) \leq N^{\eps} / \eta^2$. Therefore by Cauchy-Schwarz,
\beq
\left|  \Cov \left( N \Im [ m_N (x + \i \eta ) ] ,  \sum_{j \notin J_\alpha}  \Im [ ( \lambda_j - y - \i \eta )^{-1} ] \right) \right| \leq C N^{\eps} \eta^{-1}
\eeq
and the claim follows. \qed

In what follows we will focus on pairs of indices $(i, j) \in J_\alpha^2$. 
Fix now $\omega>0$ and assume that $H$ is a Gaussian divisible ensemble with Gaussian component of size
\beq
T = N^{- \omega} ,
\eeq
and fix $\delta_1 > \delta_2 > \delta_3 > \omega > 0$. For each $i, j$ we define
\beq
t_{ij} := N^{-\delta_2} | \gamma_i - \gamma_j | , \qquad T_{ij} := N^{-\delta_3} | \gamma_i - \gamma_j |.
\eeq
For $i$ and $j$ satisfying $|i -j | \geq N^{\delta_1} (N \eta)^{-1}$ note that the error in Theorem \ref{thm:sae-homog} satisfies,
\beq \label{eqn:sub-eta}
\frac{1}{N^2 t_{ij}} \asymp \frac{N^{\delta_2}}{N |i-j| } \leq  N^{\delta_2 - \delta_1} \eta \ll \eta
\eeq
under the assumption that $\delta_1 > \delta_2$.  Note also that,
\beq \label{eqn:sub-tij}
t_{ij} = N^{-\delta_2} | \gamma_i - \gamma_j | \asymp \frac{N^{-\delta_2}}{N} |i-j| \geq \frac{ N^{\delta_1-\delta_2}}{N} \frac{1}{N \eta} \geq \frac{N^{\delta_1-\delta_2}}{ N} \gg \frac{1}{N}.
\eeq
This partially motivates the above choices.

The goal of the remainder of the section is to estimate the quantity,
\beq
\Cov ( \Im [ ( \lambda_i - x - \i \eta )^{-1} ], \Im [ ( \lambda_j - x - \i \eta )^{-1}] )
\eeq
for fixed $i$ and $j$.  In order to do so we take a somewhat non-trivial coupling of the eigenvalues $\lambda_k (t)$ to an equilibrium process which we will now outline. In particular, the coupling depends on the pair of indices chosen, $(i, j)$. 

Recall that our eigenvalues $\{ \lambda_k \}_k$ are the eigenvalues of a Gaussian divisible ensemble with Gaussian component of size $T$, which we denote by $H_T$, with 
\beq
H_t := \e^{- t/2} H_0 + \sqrt{ 1- \e^{-t} } G 
\eeq
where $G$ is a Gaussian matrix independent from $H_0$, some fixed Wigner matrix. 
 These have the same distribution as the solution to \eqref{eqn:dbm-def}, which we denote by $\{ \lambda_k (s) \}_k$  at time $ s= T_{ij} + t_{ij}$ with initial data have the same distribution as an independent copy of $H_{T- T_{ij} - t_{ij} }$.  Denote the Brownian motions driving \eqref{eqn:dbm-def} by $B_k (s)$.
 
 Consider now an auxilliary process $\mu_k (s)$ which solves \eqref{eqn:dbm-def} for times $T_{ij} \leq s \leq T_{ij} + t_{ij}$ with the same Brownian terms $B_k (s)$ (but only for these times $s$), with initial data $\mu_k (T_{ij})$ being the eigenvalues of an independent Gaussian ensemble.  
 
 Clearly, the three families of random variables,
 \beq
 \{ \mu_k (s) \}_{k, s \geq T_{ij} }, \qquad \{ \lambda_k (0) \}_{k} , \qquad \{ B_k (s) \}_{k, 0 \leq s \leq T_{ij} }
 \eeq
 are mutually independent.


By Theorem \ref{thm:sae-homog} we have for any $\eps >0$ (we will choose also $\eps < \delta_2 - \delta_1$) that,
\beq
| \lambda_i (T_{ij} + t_{ij} ) -  \mu_i (T_{ij} + t_{ij} ) - N^{-1} ( X_1 ( \gamma_i, t_{ij} ) - X_2 ( \gamma_i, t_{ij} ) ) | \leq \frac{N^{\eps}}{N^2 t_{ij}},
\eeq
where
\beq
X_1 ( z, t ) = \frac{ \Im \sum_n \log ( \lambda_n (T_{ij} ) - z_t ) }{ \Im [ \msc (z_t ) ] }  - N \frac{ \int \Im [ \log ( x -z_t ) ] \rhosc (x) \d x }{  \Im [ \msc (z_t ) ] }
\eeq
and
\beq
X_2 ( z, t ) = \frac{ \Im \sum_n \log ( \mu_n (T_{ij} ) - z_t ) }{ \Im [ \msc (z_t ) ] }  - N \frac{ \int \Im [ \log ( x -z_t ) ] \rhosc (x) \d x }{  \Im [ \msc (z_t ) ] } ,
\eeq
and an analogous estimate changing $i$ for $j$.  By Proposition \ref{prop:char-est}  we have for $X_1 ( \gamma_i, t_{ij})$ and $X_1 ( \gamma_j, t_{ij} )$ the estimates,
\begin{align}
X_1 ( \gamma_i, t_{ij} ) = Z_i + Y_i + N^{\eps} \O ( (N t_{ij} )^{-1} )
\end{align}
for any $\eps >0$ with overwhelming probability where
\begin{align}
Y_i :=& \frac{ \sum_n \Im [ \log ( \lambda_n (0) - ( \gamma_i )_{T_{ij} + t_{ij} } )]}{ \Im [ \msc (  ( \gamma_i )_{t_{ij} }) ]} - \frac{ N \int \Im \log (x - ( \gamma_i )_{T_{ij} + t_{ij} } ) \rhosc (x) \d x}{ \Im [ \msc (  ( \gamma_i )_{t_{ij} }) ] } \notag\\
+& \frac{1}{2} \left( 1 - \frac{2}{ \beta} \right) \int_{0}^{T_{ij} } \Im [ \msc' (  ( \gamma_i )_{t_{ij} +s } ) ] \d s
\end{align}
and
\begin{align}
Z_i := \frac{1}{N^{1/2}} \sum_n  \sqrt{ \frac{2}{ \beta}} \int_0^{T_{ij}} \Im[ ( \gamma_n - ( \gamma_i )_{T_{ij} + t_{ij} -s } )^{-1} ] \d B_n (s)
\end{align}
and analogous definitions for $Y_j$ and $Z_j$.  Above, $\beta =1$ in the case of real symmetric Wigner matrices that we consider. It is $\beta=2$ for the complex Hermitian case.

Here $(\gamma_i)_{t_{ij}}$ denotes the solution at time $t_{ij}$ of the characteristic flow \eqref{eqn:char-def} with initial point $z=\gamma_i + \i 0$ (see also \eqref{eqn:char-real-line-def}).

Making the additional definition
\beq
\muij_i = \mu_i (T_{ij} + t_{ij} ) - N^{-1} X_2 (\gamma_i, t_{ij} )
\eeq
we conclude from the above discussion as well as 
  Proposition \ref{prop:char-cov} the following.
 \bep \label{prop:sub-1}
 Let $\eps >0$. With overwhelming probability, we have
 \beq \label{eqn:sub-sae-1}
 | \lambda_i (T_{ij} +t_{ij} ) - \muij_i - N^{-1}(Y_i - Z_i) | \leq \frac{N^{\eps}}{N^2 t_{ij}}
 \eeq
 and an analogous estimate for $\lambda_j$. The three families of random variables $(Z_i, Z_j)$, $(Y_i, Y_j)$ and $( \muij_i, \muij_j)$ are mutually independent, and the $Z_i, Z_j$ are jointly Gaussian. Moreover, there is a $c>0$ so that,
\beq
c ( \delta_2 - \delta_3)  \log(N) \leq \Var (Z_i ), \Var (Z_j) \leq c^{-1} ( \delta_2 - \delta_3)  \log(N)
\eeq
and,
\beq
| \Cov ( Z_i, Z_j ) | \leq \frac{1}{ c \log(N) }.
\eeq
\eep
For notational simplicity, define also
\beq
\lamij_i := \muij_i - N^{-1} (Y_i + Z_i ).
\eeq
The following is one of the main technical arguments in the paper, as obtaining an acceptable error estimate is delicate. 
\bep \label{prop:sub-2}   Let $k$ satisfy
\beq
\log_2(N) \leq k \leq A \log_2 (N).
\eeq
For any $\eps >0$ and $M >0$ there is a $C>0$ so that,
\begin{align}
& \sum_{(i, j )\in J^2_\alpha} \1_{ \{ |i-j| > \frac{ N^{\delta_1}}{N \eta} \}} \bigg| \int_{\rr^2} g_k (x) g_k (y) \Cov ( \Im [ ( \lambda_i - x - \i \eta )^{-1} ], \Im [ ( \lambda_j - y - \i \eta )^{-1} ) \d x \d y \notag \\
& - \int_{\rr^2} g_k (x) g_k (y) \Cov \big( \Im [ ( \lamij_i - x - \i \eta )^{-1} ], \Im [ ( \lamij_j - y - \i \eta )^{-1}] \big) \d x \d y \bigg| \notag\\
& \leq C \left(  \frac{N^{\eps+\delta_2}}{\eta} \| \varphi_k \|_2^2 + 2^{-Mk} N^{2} \| \varphi \|_2^2 \right)
\end{align}
\eep
\proof We first remark that by using Proposition \ref{prop:lp-support} we can restrict both of the integrals above to the square 
\beq
R:= [-5, 5]^2
\eeq
 at the cost of $C_M 2^{-M k } N^2 \| \varphi \|_2^2$, for any $M>0$.  In what follows we will only consider the integration over such $x, y$ and ignore the remaining contributions.

From the estimate \eqref{eqn:sub-sae-1} and the fact that our parameters satisfy \eqref{eqn:sub-eta} we see that for any $\eps >0$ we have,
\begin{align}
\bigg| &\Cov ( \Im [ ( \lambda_i - x - \i \eta )^{-1} ], \Im [ ( \lambda_j - y - \i \eta )^{-1} ) \notag\\
 -& \Cov ( \Im [ ( \lamij_i - x - \i \eta )^{-1} ], \Im [ ( \lamij_j - y - \i \eta )^{-1} ) \bigg| \notag\\
\leq & C \frac{N^{\eps}}{N^2 t_{ij}} \ee\left[ \frac{1}{ | \lamij_i - x - \i \eta|^2} \left| \frac{1}{ \lamij_j - y - \i \eta } - \ee \left[  \frac{1}{ \lamij_j - y - \i \eta }\right] \right|  \right] \label{eqn:sub-err-1} \\
+ & C \frac{N^{\eps}}{N^2 t_{ij}} \ee\left[ \frac{1}{ | \lamij_i - x - \i \eta|^2} \right] \ee\left[ \left| \frac{1}{ \lamij_j - y - \i \eta } - \ee \left[  \frac{1}{ \lamij_j - y - \i \eta }\right] \right|  \right]  \label{eqn:sub-err-2} \\
+ &C \frac{N^{\eps}}{N^2 t_{ij}} \ee\left[ \frac{1}{ | \lamij_j - y - \i \eta|^2} \left| \frac{1}{ \lamij_i - x - \i \eta } - \ee \left[  \frac{1}{ \lamij_i - x - \i \eta }\right] \right|  \right] \label{eqn:sub-err-3} \\
+ & C \frac{N^{\eps}}{N^2 t_{ij}} \ee\left[ \frac{1}{ | \lamij_j - y - \i \eta|^2} \right] \ee\left[ \left| \frac{1}{ \lamij_i - x - \i \eta } - \ee \left[  \frac{1}{ \lamij_i - x - \i \eta }\right] \right|  \right] \label{eqn:sub-err-4} \\
+ & C \frac{N^{\eps}}{ N^4 t_{ij}^2} \ee\left[ \frac{1}{ | \lamij_i - x - \i \eta |^2  } \frac{1}{ | \lamij_j - y - \i \eta |^2} \right] \label{eqn:sub-err-5}\\
+ & C \frac{N^{\eps}}{ N^4 t_{ij}^2} \ee\left[ \frac{1}{ | \lamij_i - x - \i \eta |^2  }\right] \ee \left[ \frac{1}{ | \lamij_j - y - \i \eta |^2} \right] \label{eqn:sub-err-6}
\end{align}
Of the terms above, \eqref{eqn:sub-err-1}, \eqref{eqn:sub-err-2}, \eqref{eqn:sub-err-3} and \eqref{eqn:sub-err-4} are handled similarly, so we just estimate \eqref{eqn:sub-err-1}.  Analogously, we will then estimate \eqref{eqn:sub-err-5} as \eqref{eqn:sub-err-6} is similar.

 We begin with \eqref{eqn:sub-err-1}.  Fixing a small $\eps_1 >0$, we break the integral over $x$ and $y$ into four terms, depending on whether $|x-\gamma_i | $ and $|y- \gamma_j|$ are less or greater than $N^{\eps_1-1}$.  That is, we expand the identity,
 \beq
1 = ( \1_{ \{ | x- \gamma_i  | < N^{\eps_1-1} \} } + \1_{ \{ | x- \gamma_i  | \geq N^{\eps_1-1} \} } ) ( \1_{ \{ | y- \gamma_j  | < N^{\eps_1-1} \} } + \1_{ \{ | y- \gamma_j  | \geq N^{\eps_1-1} \} } 
 \eeq
 and estimate the contribution of each of the four resulting terms. 
   We begin with the case that $|x - \gamma_i | , |y - \gamma_j | < N^{\eps_1-1}$. For this term, we use the fact that,
\beq \label{eqn:sub-Zi-second}
\sup_{s, t} \ee\left[ \frac{1}{ | Z_i - s + \i N \eta | |Z_j - t + \i N \eta | } \right] \leq C \log(N)^2
\eeq
 by the estimates of the covariance and variance of $Z_i, Z_j$ of Proposition \ref{prop:sub-1}.  The same estimate holds for,
\beq \label{eqn:sub-Zi-expect}
\sup_{s, t} \ee[ | Z_i - s + \i N \eta |^{-1}] + \ee[ | Z_j - t + \i N \eta |^{-1}] \leq C \log(N) .
\eeq Therefore,
\begin{align}
 & \sum_{(i, j) \in J^2_\alpha} \1_{ \{ |i-j| > \frac{ N^{\delta_1}}{N \eta} \}} \int_{R} |g_k (x) g_k (y) | \1_{ \{ |x - \gamma_i | < N^{\eps_1-1} \} } \1_{ \{ |y - \gamma_j | < N^{\eps_1-1} \} } \notag\\
 \times & \frac{N^{\eps}}{N^2 t_{ij}} \ee\left[ \frac{1}{ | \lamij_i - x - \i \eta|^2} \left| \frac{1}{ \lamij_j - y - \i \eta } - \ee \left[  \frac{1}{ \lamij_j - y - \i \eta }\right] \right|  \right]  \d x \d y \notag\\
 \leq & \sum_{(i, j) \in J^2_\alpha} \1_{ \{ |i-j| > \frac{ N^{\delta_1}}{N \eta} \}} \int_{R} |g_k (x) g_k (y) | \1_{ \{ |x - \gamma_i | < N^{\eps_1-1} \} } \1_{ \{ |y - \gamma_j | < N^{\eps_1-1} \} } \frac{N^{2\eps}}{N^2 t_{ij}} \frac{N^2}{\eta} \d x \d y \notag\\
  \leq & \frac{N^{2\eps+\delta_2}}{\eta} \sum_{(i, j) \in J^2_\alpha} \1_{ \{ |i-j| > \frac{ N^{\delta_1}}{N \eta} \}} \int_{R} |g_k (x) g_k (y) | \1_{ \{ |x - \gamma_i | < N^{\eps_1-1} \} } \1_{ \{ |y - \gamma_j | < N^{\eps_1-1} \} } \frac{1}{ |x-y| + N^{-1} } \d x \d y \notag\\
   \leq & \frac{N^{2\eps+\delta_2+ 2 \eps_1}}{\eta} \int_{R} |g_k (x) g_k (y) | \frac{1}{ |x-y| + N^{-1} } \d x \d y
 \end{align}
 The first inequality bounds $| \lamij_i - x - \i \eta|^{-2} \leq \eta^{-1} | \lamij_i - x - \i \eta|^{-1}$ and then uses \eqref{eqn:sub-Zi-second} and \eqref{eqn:sub-Zi-expect}. 
 In the second inequality we used the fact that
 \beq
 t_{ij} \geq N^{-\delta_2} |x-y| + N^{-1}
 \eeq
 which holds for all $|x-\gamma_i| < N^{\eps_1-1}$ and $|y - \gamma_j | < N^{\eps_1-1}$ as long as $\eps_1 < \delta_1- \delta_2$, due the inequalities \eqref{eqn:sub-tij}. The final inequality is from the fact that for fixed $x$ and $y$, we have
 \beq
\sum_{( i, j) \in J^2_\alpha} \1_{ \{ | x- \gamma_i | < N^{\eps_1-1} \} } \1_{ \{ |y- \gamma_j | < N^{\eps_1-1} \} } \leq CN^{2 \eps_1}.
 \eeq
For the integral, we have that 
\begin{align}
\int_{R} |g_k (x) g_k (y) | \frac{1}{ |x-y| + N^{-1} } \d x \d y  \leq \int_{R} |g_k (x)|^2 \frac{1}{ |x-y| + N^{-1} } \d x \d y \leq \log(N) \| \varphi_k \|_2^2 .
\end{align}
We now consider the contribution when $|x-\gamma_i | > N^{\eps_1-1}$ and $|y - \gamma_j | < N^{\eps_1-1}$.  For this term, we estimate,
\begin{align}
& \sum_{( i, j) \in J^2_\alpha} \1_{ \{ |i-j| > \frac{ N^{\delta_1}}{N \eta} \}} \int_{R} |g_k (x) g_k (y) | \1_{ \{ |x - \gamma_i | > N^{\eps_1-1} \} } \1_{ \{ |y - \gamma_j | < N^{\eps_1-1} \} } \notag\\
 \times & \frac{N^{\eps}}{N^2 t_{ij}} \ee\left[ \frac{1}{ | \lamij_i - x - \i \eta|^2} \left| \frac{1}{ \lamij_j - y - \i \eta } - \ee \left[  \frac{1}{ \lamij_j - y - \i \eta }\right] \right|  \right]  \d x \d y \notag\\
 & \leq \sum_{( i, j) \in J^2_\alpha} \1_{ \{ |i-j| > \frac{ N^{\delta_1}}{N \eta} \}} \int_{R} |g_k (x) g_k (y) | \1_{ \{ |x - \gamma_i | > N^{\eps_1-1} \} } \1_{ \{ |y - \gamma_j | < N^{\eps_1-1} \} } \notag\\
 \times & \frac{N^{\eps}}{N^2 t_{ij}} \ee\left[ \frac{1}{ | \lamij_i - x - \i N^{-1}|^2}\left(   \frac{1}{ \left|  \lamij_j - y - \i \eta \right| } + \ee \left[  \frac{1}{ \left| \lamij_j - y - \i \eta  \right| }\right]  \right) \right]  \d x \d y \label{eqn:sub-err-1a}.
\end{align}
We will not use any cancellation between $\lamij_j- y- \i \eta$ and its expectation, so we just estimate the contribution of $\ee [| \lamij_i - x - \i N^{-1} |^{-2} | \lamij_j - y - \i \eta |^{-1}]$ to the summation and integral, the term where the expectation hits each term separately being similar.  For this term, we use 
\beq \label{eqn:sub-err-1b}
|g_k (x) g_k (y) | \leq |g_k(x) |^2 +|g_k (y) |^2,
\eeq resulting in two integrals. For the second, we estimate,
\begin{align}
& \sum_{( i, j) \in J^2_\alpha} \1_{ \{ |i-j| > \frac{ N^{\delta_1}}{N \eta} \}} \int_{R} |g_k (y)|^2 \frac{1}{N^2 t_{ij}} \1_{ \{ |y - \gamma_j | < N^{\eps_1-1}\} }\ee\left[ \frac{1}{ | \lamij_i - x - \i N^{-1}|^2} \left| \frac{1}{ \lamij_j - y - \i \eta }\right| \right] \d x \d y \notag\\
& \leq  C \sum_{( i, j) \in J^2_\alpha}\1_{ \{ |i-j| > \frac{ N^{\delta_1}}{N \eta} \}} \int_{|y|<5} | g_k (y) |^2 \frac{1}{N t_{ij}} \1_{ \{ |y - \gamma_j | < N^{\eps_1-1}\} } \ee[ | \lamij_j - y - \i \eta |^{-1} ] \d y \notag\\
& \leq C N^{\eps}\sum_{( i, j) \in J^2_\alpha}\1_{ \{ |i-j| > \frac{ N^{\delta_1}}{N \eta} \}} \int_{|y|<5} | g_k (y) |^2 \frac{1}{t_{ij}} \1_{ \{ |y - \gamma_j | < N^{\eps_1-1}\} }  \d y \notag\\ \notag\\
& \leq CN^{1+2\eps+\delta_2} \sum_{j \in J_\alpha}  \int_{|y|<5} |g_k (y) |^2 \1_{ \{ |y - \gamma_j | < N^{\eps_1-1}\} } \d y \leq C N^{1+ 2\eps+\eps_1+\delta_2} \| \varphi_k \|_2^2
\end{align}
In the first inequality we did the integral over $x$ which contributes a factor of $N$. The second inequality follows from estimating the expectation using \eqref{eqn:sub-Zi-expect}. The third inequality uses the fact that for all $j \in J_\alpha$,
\beq \label{eqn:sub-t-sum}
\sum_{i\in J_\alpha} \1_{ \{ |i-j| > \frac{ N^{\delta_1}}{N \eta} \}}  \frac{1}{ t_{ij} } \leq N^{1+ \eps+\delta_2} .
\eeq
The last inequality uses that, for all $y$,
\beq \label{eqn:sub-y-sum}
\sum_{j \in J_\alpha} \1_{ \{ |y - \gamma_j | < N^{\eps_1-1}\} }  \leq C N^{\eps_1}.
\eeq
The other contribution to \eqref{eqn:sub-err-1a} after the Schwarz inequality \eqref{eqn:sub-err-1b} is bounded above by
\begin{align}
& \sum_{( i, j) \in J^2_\alpha} \1_{ \{ |i-j| > \frac{ N^{\delta_1}}{N \eta} \}} \int_{ R} |g_k (x)|^2 \frac{1}{N^2 t_{ij}} \1_{ \{ |y - \gamma_j | < N^{\eps_1-1}\} }\ee\left[ \frac{1}{ | \lamij_i - x - \i N^{-1}|^2} \left| \frac{1}{ \lamij_j - y - \i \eta }\right| \right] \d x \d y \notag\\
& \leq \sum_{( i, j) \in J^2_\alpha} \1_{ \{ |i-j| > \frac{ N^{\delta_1}}{N \eta} \}} \int_{|x|<5} |g_k (x)|^2 \frac{N^{\eps}}{N^2 t_{ij}} \frac{1}{ ( \gamma_i  - x)^2 + N^{-2} } \d x \notag\\
&\leq \int_{|x|<5} |g_k(x)|^2 \sum_{i \in J_\alpha }  \frac{N^{2 \eps+\delta_2}}{N} \frac{1}{ ( \gamma_i  - x)^2 + N^{-2} }  \d x \leq C N^{2 \eps+\delta_2} N \| \varphi_k \|_2^2.
\end{align}
In the first inequality we did the integral over $y$ which contributes a factor of $\log(N)$ at worst, and then used rigidity to estimate $| \lamij_i - x - \i N^{-1} | \geq N^{-\eps} | \gamma_i - x - \i N^{-1} |$ with overwhelming probability. In the second inequality we used again \eqref{eqn:sub-t-sum} but for sums over $j$. The final inequality follows from the fact that,
\beq \label{eqn:sub-imgamma}
\sum_{i \in \J_\alpha }  \frac{1}{N} \frac{1}{ ( \gamma_i  - x)^2 + N^{-2} } \leq CN,
\eeq
for all $ x\in \rr$. 

Returning to \eqref{eqn:sub-err-1} we now turn to the case $|x-\gamma_i | < N^{\eps_1-1}$ and $|y - \gamma_j | > N^{\eps_1-1}$.  For such $y$, we see by rigidity that, with overwhelming probability
\beq
\left| \frac{1}{ \lamij_j - y - \i \eta } - \ee\left[ \frac{1}{ \lamij_j - y - \i \eta } \right] \right| \leq \frac{N^{\eps}}{N} \frac{1}{ ( \gamma_j - y)^2 + N^{-2} } .
\eeq 
So, we are led to estimate,
\begin{align} \label{eqn:sub-hh1}
  & \sum_{( i, j) \in J^2_\alpha} \1_{ \{ |i-j| > \frac{ N^{\delta_1}}{N \eta} \}} \int_{R} |g_k (x) g_k (y) | \1_{ \{ |x - \gamma_i | < N^{\eps_1-1} \} } \frac{N^{\eps}}{N^3 t_{ij} } \frac{1}{( \gamma_j - y)^2+N^{-2} } \ee[ | \lamij_i - x - \i \eta |^{-2} ] \d x \d y \notag\\
 & \leq  \sum_{( i, j) \in J^2_\alpha}  \1_{ \{ |i-j| > \frac{ N^{\delta_1}}{N \eta} \}} \int_{R} |g_k (x) g_k (y) | \1_{ \{ |x - \gamma_i | < N^{\eps_1-1} \} } \frac{N^{2 \eps}}{N^2 \eta t_{ij} } \frac{1}{( \gamma_j - y)^2+N^{-2} } \d x \d y \notag\\
\end{align}
where we used again \eqref{eqn:sub-Zi-expect}.  To bound the last line of \eqref{eqn:sub-hh1},  we apply again the Schwarz inequality \eqref{eqn:sub-err-1b} and treat the two integrals separately. We first have,
\begin{align}
 & \sum_{( i, j) \in J^2_\alpha} \1_{ \{ |i-j| > \frac{ N^{\delta_1}}{N \eta} \}} \int_{R} |g_k (x) |^2 \1_{ \{ |x - \gamma_i | < N^{\eps_1-1} \} } \frac{N^{2 \eps}}{N^2 \eta t_{ij} } \frac{1}{( \gamma_j - y)^2+N^{-2} } \d x \d y \notag\\
\leq & \sum_{( i, j) \in J^2_\alpha} \1_{ \{ |i-j| > \frac{ N^{\delta_1}}{N \eta} \}} \int_{|x|<5} |g_k (x) |^2 \1_{ \{ |x - \gamma_i | < N^{\eps_1-1} \} } \frac{N^{2 \eps}}{N \eta t_{ij} }  \d x \notag\\
\leq & \frac{N^{3 \eps+\delta_2}}{\eta} \sum_{i \in J_\alpha} \int_{|x|<5} |g_k(x) |^2 \1_{ \{ |x - \gamma_i | < N^{\eps_1-1} \} }  \d x  \leq \frac{ N^{3\eps+\eps_1+\delta_2}}{\eta} \| \varphi_k \|_2^2 .
\end{align}
In the first inequality we did the $y$-integration, which contributes a factor of $N$. In the second inequality we did the sum over $j$ and used \eqref{eqn:sub-t-sum}, and used finally \eqref{eqn:sub-y-sum} in the third inequality. The other contribution to the last line of \eqref{eqn:sub-hh1} is, 
\begin{align}
 & \sum_{( i, j) \in J^2_\alpha}\1_{ \{ |i-j| > \frac{ N^{\delta_1}}{N \eta} \}} \int_{ R} |g_k (y) |^2 \1_{ \{ |x - \gamma_i | < N^{\eps_1-1} \} } \frac{N^{2 \eps}}{N^2 \eta t_{ij} } \frac{1}{( \gamma_j - y)^2+N^{-2} } \d x \d y \notag\\ 
 \leq & \sum_{( i, j) \in J^2_\alpha} \1_{ \{ |i-j| > \frac{ N^{\delta_1}}{N \eta} \}} \int_{|y|<5} |g_k (y) |^2 \frac{N^{2 \eps+\eps_1}}{N^3 \eta t_{ij} } \frac{1}{( \gamma_j - y)^2+N^{-2} }  \d y \notag\\
 \leq & \sum_{ j \in J_\alpha} \int_{|y|<5} |g_k (y) |^2 \frac{N^{3 \eps+\eps_1+\delta_2}}{N^2 \eta  } \frac{1}{( \gamma_j - y)^2+N^{-2} }  \d y \leq C \frac{N^{3 \eps + \eps_1+\delta_2}}{ \eta} \| \varphi_k \|_2^2.
\end{align}
In the first inequality we did the integration over $x$ which is restricted to an interval of size $\O( N^{\eps_1-1} )$. In the second inequality we did the summation over $i$ using \eqref{eqn:sub-t-sum}. The final inequality uses \eqref{eqn:sub-imgamma}.  

Finally, if $|x - \gamma_i | > N^{\eps_1-1}$ and $|y - \gamma_j | > N^{\eps_1-1}$ then,
\begin{align}
 & \1_{ \{ |x - \gamma_i |  > N^{\eps_1-1} \} } \1_{ \{ | \gamma_j - y | > N^{\eps_1-1} \} } \ee\left[ \frac{1}{ | \lamij_i - x - \i \eta|^2} \left| \frac{1}{ \lamij_j - y - \i \eta } - \ee \left[  \frac{1}{ \lamij_j - y - \i \eta }\right] \right|  \right] \notag\\
\leq & \frac{N^{\eps}}{N} \frac{1}{ ( \gamma_i - x)^2 + N^{-2} } \frac{1}{ ( \gamma_j - y )^2 + N^{-2} } ,
\end{align}
by rigidity. 
Now,
\begin{align}
 & \sum_{( i, j) \in J^2_\alpha} \1_{ \{ |i-j| > \frac{ N^{\delta_1}}{N \eta} \}}  \int_{R} \frac{1}{N^3 t_{ij} } |g_k(x) g_k (y) |  \frac{1}{ ( \gamma_i - x)^2 + N^{-2} } \frac{1}{ ( \gamma_j - y )^2 + N^{-2} }\d x \d y \notag\\
\leq & \sum_{( i, j) \in J^2_\alpha}   \1_{ \{ |i-j| > \frac{ N^{\delta_1}}{N \eta} \}}  \int_{R} \frac{1}{N^3 t_{ij} } |g_k(x)|^2   \frac{1}{ ( \gamma_i - x)^2 + N^{-2} } \frac{1}{ ( \gamma_j - y )^2 + N^{-2} }\d x \d y  \notag\\
\leq & \sum_{( i, j) \in J^2_\alpha}  \1_{ \{ |i-j| > \frac{ N^{\delta_1}}{N \eta} \}}  \int_{|x|<5} \frac{1}{N^2 t_{ij} } |g_k(x)|^2   \frac{1}{ ( \gamma_i - x)^2 + N^{-2} } \d x  \notag\\
\leq & N^{\eps+\delta_2} \sum_{i \in J_\alpha} \int_{|x|<5} \frac{1}{N} |g_k (x) |^2 \frac{1}{ ( \gamma_i -x )^2 + N^{-2} } \d x \leq N^{1+\eps+\delta_2} \| \varphi_k \|_2^2. 
\end{align}
In the first inequality we applied \eqref{eqn:sub-err-1b}. In the second we did the integration over $y$. In the third we used \eqref{eqn:sub-t-sum} to do the summation over $j$. In the last inequality we then used \eqref{eqn:sub-imgamma}. This completes the estimation of \eqref{eqn:sub-err-1}. It remains to consider \eqref{eqn:sub-err-5}.  We have,
\begin{align}
 & \sum_{i, j \in J_\alpha}   \1_{ \{ |i-j| > \frac{ N^{\delta_1}}{N \eta} \}}  \int_{R} \frac{1}{N^4 t^2_{ij} } |g_k(x) g_k (y) | \ee[ | \lamij_i - x - \i \eta |^{-2} | \lamij_j - y - \i \eta|^{-2} ] \d x \d  y \notag\\
 \leq  &  \sum_{i, j \in J_\alpha}   \1_{ \{ |i-j| > \frac{ N^{\delta_1}}{N \eta} \}}  \int_{R} \frac{1}{N^4 t^2_{ij} } |g_k(x)|^2  \ee[ | \lamij_i - x - \i \eta |^{-2} | \lamij_j - y - \i \eta|^{-2} ] \d x \d  y \notag\\
  \leq  &  \sum_{i, j \in J_\alpha}   \1_{ \{ |i-j| > \frac{ N^{\delta_1}}{N \eta} \}}  \int_{|x|<5} \frac{1}{N^4 t^2_{ij} \eta } |g_k(x)|^2  \ee[ | \lamij_i - x - \i \eta |^{-2}  ] \d x \notag\\
    \leq  & N^{2 \delta_2 - \delta_1}  \sum_{i \in J_\alpha}  \int_{|x|<5} \frac{1}{N} |g_k(x)|^2  \ee[ | \lamij_i - x - \i \eta |^{-2}  ] \d x \notag\\
     \leq & N^{2 \delta_2 - \delta_1}    \sum_{i \in J_\alpha}  \int_{|x|<5} \frac{1}{N \eta} |g_k(x)|^2  \ee[ | \lamij_i - x - \i \eta |^{-1}  ] \d x \leq  \frac{C N^{\eps+ 2 \delta_2 - \delta_1}}{ \eta}  \| \varphi_k \|_2^2
\end{align}
In the first inequality, we used \eqref{eqn:sub-err-1b}. In the second inequality, we did the integration over $y$. In the third inequality, we used
\beq
\frac{1}{N^3} \sum_{j}   \1_{ \{ |i-j| > \frac{ N^{\delta_1}}{N \eta} \}}  \frac{1}{ t_{ij}^2 } \leq \frac{C}{N} \sum_{j}   \1_{ \{ |i-j| > \frac{ N^{\delta_1}}{N \eta} \}} \frac{N^{ 2 \delta_2}}{|i-j|^2} \leq C N^{2 \delta_2 - \delta_1} \eta .
\eeq
In the last inequality we used \eqref{eqn:sub-Zi-expect}. This completes the proof. \qed

\bep \label{prop:sub-3}
Let $\delta_1, \delta_2, \delta_3$ be as above.  Then, for $B_1 >0$ large enough, depending on $\delta_1, \delta_2$ and $\delta_3$ we have for
\beq
10 \log_2 (N) \geq k \geq \log_2 (B_1 N)
\eeq
that,
\beq
\left| \int_{\rr^2} g_k (x) g_k (y) \Cov \big( \Im [ ( \lamij_i - x - \i \eta )^{-1} ], \Im [ ( \lamij_j - y - \i \eta )^{-1} )]\big) \d x \d y \right| \leq C N^{-100} \| \varphi_k \|_2^2
\eeq

\eep
\proof By Fourier duality and the independence of $Z_i, Z_j$ from $\muij_i, \muij_j, Y_i, Y_j$ (see Proposition \ref{prop:sub-1}) we have,
\begin{align}
& \int_{\rr^2} g_k (x) g_k (y) \ee\big[ \Im [ ( \lamij_i - x - \i \eta )^{-1} ]  \Im [ ( \lamij_j - y - \i \eta )^{-1} ]\big] \d x \d y \notag\\
= & \int_{\rr^2} \hat{ \varphi}_k (\xi_1 ) \hat{ \varphi}_k ( \xi_2 ) \ee\left[ \exp\left( \i \xi_1 (N^{-1} Y_i - \muij_i ) + \i \xi_2 (N^{-1} Y_j - \muij_j ) \right) \right] \notag\\
& \times \ee \left[ \exp \left( \i \xi_1 N^{-1} Z_i + \i \xi_2 N^{-1} Z_j \right) \right] \d \xi_1 \d \xi_2.
\end{align}
By Proposition \ref{prop:sub-1} we have,
\beq
\left| \ee \left[ \exp \left( \i \xi_1 N^{-1} Z_i + \i \xi_2 N^{-1} Z_j \right) \right] \right| \leq \exp \left[ - c ( \delta_3 - \delta_2 ) \log(N) N^{-2} ( \xi_1^2 + \xi_2^2 ) \right]
\eeq
By definition, $\hat{\varphi}_k ( \xi)$ is non-zero only if $| \xi |\geq B_1 N/10$.  By taking $B_1$ large enough we see that,
\beq
\left| \ee \left[ \exp \left( \i \xi_1 N^{-1} Z_i + \i \xi_2 N^{-1} Z_j \right) \right] \right| \leq N^{-10^7}
\eeq
for all such $\xi_1, \xi_2$. Using this and the fact that 
\beq
\| \hat{ \varphi}_k  \|_1 \leq C N^{10} \| \hat{ \varphi}_k \|_2 = C N^{10} \| \varphi_k \|_2
\eeq
by the fact that $\hat{\varphi}_k$ is supported in a set of size at most $\O (N^{10})$, we conclude that
\beq
\left| \ee[ \varphi_k ( \lamij_i ) \varphi_k ( \lamij_j ) ] \right| \leq N^{-100} \| \varphi_k \|_2^2 .
\eeq
The estimate for
\beq
\ee[ \varphi_k ( \lamij_i )] \ee[ \varphi_k ( \lamij_j ) ]
\eeq
follows from similar reasoning. \qed

\subsection{Proof of Theorem \ref{thm:sub}} \label{sec:sub-proof}

Let $\omega >0$ and $\delta >0$ be as in the statement of the theorem. Let $\eps>0$ and assume $100 \eps < \delta - \omega$. Let $\delta_1 = \delta$.  Let $\kappa >0$ so that $\varphi$ is supported in $[-2+\kappa, 2 - \kappa]$ and choose $\alpha>0$ such that $\gamma_{ \alpha N } \leq -2 + \kappa/10$. For $k$ that satisfy,
\beq
 (2- \delta_1+\eps) \log_2 (N) \leq k \leq A \log_2 (N)
\eeq
we see that $\eta = 2^{-k}$ satisfies,
\beq
\frac{N^{\delta_1}}{N \eta} \geq N^{1+\eps} \gg 1.
\eeq
Therefore, by Propositions \ref{prop:sub-diag-cut-off} and \ref{prop:sub-edge}, we have for 
\beq
(2- \delta_1+\eps) \log_2 (N) \leq k \leq 2 \log_2 (N)
\eeq
that,
\beq
\Var ( \tr ( \varphi_k (H) ) ) \leq C N^{\delta+\eps} \| \varphi_k \|_2^2 2^{k}  + C 2^{-Mk } \| \varphi \|_2^2
\eeq
For
\beq
2 \log_2 (N) \leq k \leq A \log_2 (N),
\eeq
we can apply Propositions \ref{prop:sub-diag-cut-off} and \ref{prop:sub-edge} with $\delta_1 = \eps$ to conclude that,
\beq
\Var ( \tr ( \varphi_k (H) ) ) \leq C N^{2\eps} \| \varphi_k \|_2^2 2^k + C 2^{-Mk } \| \varphi \|_2^2.
\eeq
Let $\delta_2$ and $\delta_3$ satisfy $\omega < \delta_3 < \delta_2 < \delta_1$ and be equally spaced between $\omega$ and $\delta_1$. By Propositions \ref{prop:sub-diag-cut-off} and \ref{prop:sub-2} we see that for $k$ satisfying,
\beq
\log_2 ( N) \leq k \leq (2- \delta_1 + \eps ) \log_2 (N)
\eeq
we have,
\begin{align}
& \Var ( \tr \varphi_k (H) ) ) \notag\\
= & \sum_{i, j} \1_{ \{ |i-j| > \frac{ N^{\delta_1}}{N \eta} \} } \int_{\rr^2} g_k (x) g_k (y) \Cov ( \Im [ ( \lamij_i - x - \i \eta )^{-1} ], \Im [ ( \lamij_j - y - \i \eta )^{-1}]) \d x \d y \notag\\
+ & \O ( N^{\delta+\eps} \| \varphi_k \|_2^2 2^k  + N^2 2^{-M k } \| \varphi \|_2^2).
\end{align}
Due to the lower bound $k \geq \log_2 (N)$ we may absorb the $N^2$ into the $2^{-Mk}$. 
Now, the first term on the right side is estimated by Proposition \ref{prop:sub-3}, after taking $B_1$ large enough, which completes the proof. \qed

\section{Micro and mesoscopic frequencies} \label{sec:meso}

Let $H$  be a Wigner matrix. 
Let $\mfs >0$. We consider $\varphi \in H^{1/2+ \mfs}$. Recall the definition of $\varphi_k$ and $g_k$ as in Section \ref{sec:lpt}.  In this section, we fix a small $\frac{1}{100} > \mfb >0$ and consider $k$ in the range 
\beq
\mfb \log_2 (N) \leq k \leq \log_2 ( N \mfb^{-1} ).
\eeq
Fixing a large $A_1 >0$, we will also introduce the following,
\beq
\hat{f}_k ( \xi ) := \e^{ A_1 N^{-2}  \xi^2/2} \hat{g}_k ( \xi).
\eeq
We will see below that $A_1$ is chosen large depending on the fourth cumulant $s_4$.  We first show the following.
\bep  \label{prop:meso-a4} Let $1\leq p\leq \infty$. Let $0 \leq k \leq \log_2 (N \mfb^{-1} )$.  There is a $C>0$ (depending on $A_1$ and $\mfb$) so that,
\beq
\| f_k \|_p \leq C \| \varphi_k \|_p
\eeq
for all such $k$.
\eep
\proof This is similar to the proof of Theorem \ref{thm:poisson-lp}. Let $\hat{\psi}$ be a bump function adapted to a ball of radius $R$, with $R$ chosen large enough depending on $\mfb$ so that
\beq
\hat{\psi} ( \xi/N  ) \hat{g}_k ( \xi ) = \hat{g}_k (\xi)
\eeq
for all $\xi$ and $k \leq \log_2 ( N \mfb^{-1})$. Then, if $h$ is the inverse Fourier transform of the function,
\beq
\hat{h} (\xi ) = \hat{\psi} ( \xi/N) \e^{ A_1 \xi^2 / (2 N^2) },
\eeq
the norm $\| h \|_1$ is a constant not depending $N$. Since
\beq
f_k = h \star g_k
\eeq
we see that $\| h_k \|_p \leq \|  h\|_1 \| g_k \|_p$ by Young's inequality. On the other hand, for $k \geq 0$ it was shown in the proof of Theorem \ref{thm:poisson-lp} (i.e., \cite[Theorem 5]{sosoewong} ) that,
\beq
\| g_k \|_p \leq C \| \varphi_k \|_p
\eeq
for all $k \geq 0$. The claim follows. \qed

In a similar way to Proposition \ref{prop:lp-support} we have the following.
\bep \label{prop:meso-lp-supp}
Let $\varphi \in H^s$ or $C^\alpha$, and let $k$ satisfy
\beq
\mfb \log_2 (N) \leq k \leq \log_2 (N \mfb^{-1})
\eeq
Let $\kappa>0$ and assume that $\varphi$ is supported in $(-2+\kappa, 2- \kappa)$. Let $D>0$. There is a constant $C>0$ so that for any $x \notin (-2 + \kappa/2, 2-\kappa/2)$ and $k$ as above we have
\beq
|f_k (x) | \leq \frac{C N^{-D}}{1+|x|^D} \| \varphi \|_2.
\eeq
\eep
Now, let $Z_1$ and $Z_2$ be independent centered Gaussians with mean $0$ and variance $A_1 / N^2$. We use the notation $\tilZ = (Z_1, Z_2)^T$. Our starting point is the following representation for the variance associated to $\varphi_k$. Obviously,
\beq
\Var ( \tr  \varphi_k (H) ) \leq \ee \left[ \left( \tr \varphi_k (H) - N \int \varphi_k (x) \rhosc (x) \d x \right)^2 \right],
\eeq
and 
\begin{align} \label{eqn:meso-var-rep}
 &\ee \left[ \left( \tr \varphi_k (H) - N \int \varphi_k (x) \rhosc (x) \d x \right)^2 \right]\notag\\
= & \int_{\rr^2} f_k (x) f_k (y) N^2 \ee^{(\tilZ ) } \ee \big[ \Im [ (m_N - \msc) ( x +Z_1+ \i \eta ) ] \Im [ (m_N - \msc) ( y+Z_2 + \i \eta )] \big]  \d x \d y,
\end{align}
where the outer expectation is with respect to $\tilZ$ and the inner expectation is conditional on $\tilZ$.  The role of $\tilZ$ will only become apparent later. For now, we note that due to the Gaussian decay of $\tilZ$ and Proposition \ref{prop:meso-lp-supp}, we can restrict our attention to the case that $x+Z_1$ and $y+Z_2$ are well-separated from the edges of the semicircle distribution.

We now turn to the application of Theorem \ref{thm:sae-homog}. We now fix $\omega >0$ and let
\beq
T = N^{-\omega}.
\eeq
The results of this section will be applied with $\omega >0$  small, depending on the regularity of $\varphi$ that is quantified in terms of the parameter $\mfs >0$ above. We consider Gaussian divisible ensembles, 
\beq \label{eqn:meso-hh1}
H_T = \e^{-T/2} W + \sqrt{ 1- \e^{-T}} G.
\eeq
where $W$ is a Wigner matrix. 
We consider the eigenvalues of $H_T$ to be those of a DBM  process as in \eqref{eqn:dbm-def} which we denote by $\{ \lambda_i (t) \}_{i=1}^N$, $0\leq t\leq T$. We couple this to an equilibrium process $\{ \mu_i (t) \}_{i=1}^N$ as described in Section \ref{sec:homog}.

Define $J_\alpha$ as in \eqref{eqn:J-def}, and for $\kappa >0$ we consider $\alpha >0$ so that $\gamma_{ \alpha N } < -2 + \kappa / 10$. 
Fixing $x \in (-2+\kappa, 2- \kappa)$, we consider the estimate
\beq \label{eqn:meso-sae}
| \lambda_i (T) - \hatmu_i (T,x) - N^{-1} Y (x) | \leq N^{\eps} \frac{ N |x- \gamma_i | + 1}{N^2 T}
\eeq
which, thanks to Theorem \ref{thm:sae-homog}, holds with overwhelming probability for $i \in J_\alpha$, where we defined
\beq
Y(x) := X^W (x, T)
\eeq
and
\beq
\hatmu_i (T,x) := \mu_i (T) - N^{-1} X^G(x, T),
\eeq
with the $X$'s defined in \eqref{eqn:XW-def} and \eqref{eqn:XG-def}. Note that the random variable $X^W (x, T)$ is independent of the $\hatmu_i$'s, and that by rigidity,
\beq
|N^{-1} Y(x) | + | N^{-1} X^G(x, T) | \leq N^{\eps-1}
\eeq
with overwhelming probability, for any $\eps >0$. 
 We also have the analogous estimates to all of those above with $x$ replaced by $y$.  
  We now use the estimate \eqref{eqn:meso-sae} to prove the following.
\bep \label{prop:meso-homog-1}
Let $\eps_1 >0$ and $\eps >0$. Fix $x, y \in (-2 + \kappa, 2 - \kappa)$.  Let $\hatx$ and $\haty$ satisfy
\beq
| x - \hatx | + | y - \haty | \leq \frac{N^{\eps_1/10}}{N}.
\eeq
 Then,
 \begin{align}
 &N^2 \ee[  \Im [ (m_N - \msc ) ( \hatx + \i \eta ) ] \Im [(m_N - \msc ) ( \haty + \i \eta )]] \notag\\
 = & \ee\bigg\{ \sum_i \Im [ ( \hatmu_i(x) - N^{-1}Y(x) - \hatx - \i \eta )^{-1}  - \msc ( \hatx + N^{-1} Y(x) + \i \eta ) ] \notag\\
 &\times   \sum_i \Im [ ( \hatmu_i(y) - N^{-1} Y(y) - \haty - \i \eta )^{-1}  - \msc ( \haty + N^{-1} Y(y) + \i \eta ) ] \bigg\} \notag\\
 + & \eta^{-1} N^{\eps_1+\eps} \O \left(  1+ T^{-1} \right).
 \end{align}
\eep
\proof The estimate \eqref{eqn:meso-sae} implies,
\beq
\left| \frac{1}{ \lambda_i(T) - \hatx - \i \eta } - \frac{1}{\hatmu_i (T, x) - \hatx - N^{-1} Y(x) - \i \eta  } \right| \leq N^{\eps_1+\eps} \frac{N | x - \gamma_i | + 1}{N^2 T} \frac{1}{ | \gamma_i - x - \i \eta|^2},
\eeq
with overwhelming probability for $i \in J_\alpha$.  For $i \notin J_\alpha$ we have that
\beq
|\lambda_i (T) - \gamma_i | +| \gamma_i- ( \mu_i (T, x) - N^{-1} Y (x)) | \leq \frac{N^{\eps}}{N^{2/3} \min \{ i^{1/3}, (N+1-i)^{1/3} \} } 
\eeq
with overwhelming probability, which implies, since $\gamma_i$ is well-separated from $\hatx$ that
\beq
\left| \frac{1}{ \lambda_i (T) - \hatx - \i \eta } - \frac{1}{ \hatmu_i (T, x) - \hatx - N^{-1} Y(x) - \i \eta } \right| \leq \frac{C N^{\eps}}{N^{2/3} \min \{ i^{1/3}, (N+1-i)^{1/3} \}  }
\eeq
  Summing over $i$ we see that with overwhelming probability that,
\begin{align}
\left| N \Im [ m_N ( \hatx + \i \eta ) ] - \sum_i \Im ( \hatmu_i (T, x) - \hatx - N^{-1} Y(x) - \i \eta )^{-1} \right| \leq N^{\eps_1+2\eps} \frac{1}{T } + N^{\eps}
\end{align}
By the smoothness of $\msc$ we see that,
\beq
| \msc ( \hatx + \i \eta ) -  \msc ( \hatx + N^{-1} Y(x) + \i \eta )  | \leq N^{\eps-1},
\eeq
with overwhelming probability. 
Therefore, with overwhelming probability,
\begin{align}
 & N \Im [ m_N ( \hatx + \i \eta ) ] - N \Im [\msc ( \hatx + \i \eta ) ] \notag\\
= & \sum_i \Im ( \hatmu_i (T, x) - \hatx - N^{-1} Y(x) - \i \eta )^{-1} -N \Im\msc ( \hatx + Y(x) + \i \eta ) \notag\\
+ & N^{\eps_1+2\eps} \O \left(  1+ \frac{1}{T } \right). \label{eqn:meso-1}
\end{align}
By rigidity, the second last line of \eqref{eqn:meso-1} is $N^{\eps} \O ( \eta^{-1} )$ with overwhelming probability. We also arrive at similar estimates changing $x$ and $\hatx$ for $y$ and $\haty$. The claim then follows by expanding out the expectation in terms of the quantities on the second last line of \eqref{eqn:meso-1}. \qed

Eventually we are going to use the estimates of Theorem \ref{thm:wig-smooth-lss} to replace the appearance of $Y(x)$ and $Y(y)$ in the estimates above by explicit distributions. In order to do this, we first require some estimates that will be useful when we apply Fourier duality. For this, let $\chi$ be a smooth bump function that is $1$ on $[-1, 1]$ and $0$ outside of $[-2, 2]$. Fix now an
\beq
\eps_2 >0
\eeq
and let
\beq
\rho_2 (x) = \chi \big( x / (N^{\eps_2} + (N \eta ) )\big) 
\eeq
and consider the functions, defined for $E \in \rr$ by,
\beq
F_{1, \hatx,x } (E) := \sum_i \Im \left[ (N( \hatmu_i (T, x) - \hatx) - E - \i N\eta )^{-1}\right] -\Im\left[ \msc ( \hatx + N^{-1} E + \i \eta ) \right]
\eeq
and 
\beq
F_{2, \hatx,x} (E) :=  \rho_2 (E) F_{1, \hatx, x} (E).
\eeq
We have the following estimate for the Fourier transform of the (random) function $F_2$,
\beq
\hat{F}_{2, x, \hatx} ( \xi ) = \int \e^{ \i x \xi } F_{2, x, \hatx} ( z ) \d z.
\eeq
\bep \label{prop:meso-obs-1} Let $\kappa>0$ and let $|x|, | \hat{x} | \leq 2 - \kappa$. 
For any $k >0$ and $\eps >0$ there is a constant $C$ such that, with overwhelming probability,
\beq
| \hat{F}_{2, \hatx, x} (\xi ) | \leq C N^{\eps_2+\eps} \frac{1}{ 1 + | \xi N \eta |^k }.
\eeq
as well as,
\beq
| \hat{F}_{2, \hatx, x} (\xi ) | \leq CN^{\eps+\eps_2}.
\eeq
\eep
\proof For the proof, denote
\beq
\hat{m} (w) := \frac{1}{N} \sum_i \frac{1}{ \hatmu_i (T, x) - w },
\eeq
so that
\beq
F_{1, \hatx, x} (E) = \Im \left[ \hat{m} ( \hatx + \i \eta + E N^{-1}) - \msc ( \hatx + \i \eta + EN^{-1} ) \right].
\eeq
We have,
\beq
| \hatmu_i (T, x) - \gamma_i | \leq \frac{N^{\eps}}{N^{2/3} \min \{ i^{1/3}, (N+1-i )^{1/3} \} }
\eeq
with overwhelming probability. Using these estimates, it is not hard to establish that
\beq \label{eqn:meso-hh2}
| \hat{m} (E + \i \eta ) - \msc (E + \i \eta ) | \leq \frac{N^{\eps}}{N \eta }
\eeq
for any $E \in \rr$ with overwhelming probability. The second estimate of the Proposition now follows because the function $F_{2, \hatx, x}$ is supported in a set of size $\O ( N^{\eps_2} +  N \eta )$.  

By the Cauchy integral formula and the estimates \eqref{eqn:meso-hh2},
\beq
| \del^k_z ( \hat{m} (E + \i \eta ) - \msc (E + \i \eta )  )| \leq C_k \frac{ N^{\eps}}{ N \eta^{k+1}}
\eeq
with overwhelming probability. For holomorphic $f$ we have by the Cauchy-Riemann equations,
\beq
\del_E \Im[f] = \Im [ f'] ,
\eeq
and so  the estimate,
\beq
| \del^k_E F_{1, \hatx, x} (E) | \leq \frac{C_k}{(N \eta )^{k+1} }
\eeq
holds. The first estimate now follows. \qed

For $z_1, z_2 \in \rr$, define $H(z_1, z_2)$ by,
\beq
H (z_1, z_2 ) := \ee^{( \mu ) } \left[ F_{2, x, x} (z_1) F_{2, y, y} (z_2 ) \right].
\eeq
where the notation $\ee^{( \mu)}$ means the expectation over the $\hatmu_i$ that appear in the definition of the functions $F_2$ above.  The previous proposition implies the following Corollary for the behavior of the Fourier transform of $H$.
\bec \label{cor:H-fourier} 
Let $\eps >0$ and $D>0$. Then for $|\xi_1| + | \xi_2 | \geq N^{\eps} ( N \eta )^{-1} $ we have,
\beq
| \hat{H} (\xi_1, \xi_2 ) | \leq C N^{-D} \frac{1}{ 1 + | \xi_1|^4 + | \xi_2 |^4 }.
\eeq
Moreover, for any $\xi_1, \xi_2$ we have for any $\eps >0$,
\beq
| \hat{H} ( \xi_1, \xi_2 ) | \leq C N^{\eps+2 \eps_2}.
\eeq
\eec
\proof On the event of $\O (N^{-D} )$ probability that the estimates of the previous Proposition do not hold, we have the deterministic estimate,
\beq
| \del_E^k  F_{2, x, x} ( E) | \leq C_k \eta^{-1} (N \eta)^{-k}
\eeq
for any $k \geq 0$ and so we conclude the proof. \qed


Given $x$, $y$ and $T>0$ we now define a few matrices. Let $s_3$ and $s_4$ denote the third and fourth cumulants of the off-diagonal entries $\sqrt{N} W_{12}$ of the Wigner matrix $W$ in the definition \eqref{eqn:meso-hh1} of $H_T$.  For functions $\phi_i$, we define the quadratic forms
\beq
V_1 ( \phi_1, \phi_2 ) := \frac{1}{ 2 \pi^2}  \int_{-2}^2 \int_{-2}^2 \frac{ \left( \phi_1 (x) - \phi_1 (y) \right) \left( \phi_2 (x) - \phi_2 (y) \right) }{(x-y)^2} \frac{ 4- x y}{ \sqrt{4-x^2}\sqrt{ 4 - y^2}} \d x \d y
\eeq
and
\begin{align}
V_2 ( \phi_1, \phi_2 ) &= \frac{s_4}{2 \pi^2} \left( \int_{-2}^2 \phi_1 (x) \frac{ 2 -x^2}{ \sqrt{4-x^2} } \d x \right) \left( \int_{-2}^2 \phi_2 (x) \frac{ 2 -x^2}{ \sqrt{4-x^2} } \d x \right) \notag\\
&- \frac{s_3}{  \pi^2 N^{1/2}} \left( \int_{-2}^2 \phi_1 (x) \frac{ 2 -x^2}{ \sqrt{4-x^2} } \d x \right) \left( \int_{-2}^2 \phi_2 (x) \frac{x}{ \sqrt{4-x^2} } \d x \right) \notag\\
& - \frac{s_3}{  \pi^2 N^{1/2}} \left( \int_{-2}^2 \phi_2 (x) \frac{ 2 -x^2}{ \sqrt{4-x^2} } \d x \right) \left( \int_{-2}^2 \phi_1 (x) \frac{x}{ \sqrt{4-x^2} } \d x \right)
\end{align}
and the third-order homogeneous polynomial
\begin{align}
P_1 (\xi_1, \xi_2 , \phi_1, \phi_2 ) = - \frac{s_3}{12 \pi^3} \left( \xi_1 \int_{-2}^2 \phi_1 (x) \frac{ x}{ \sqrt{4-x^2}} \d x + \xi_2 \int_{-2}^2 \phi_2 (x) \frac{ x}{ \sqrt{4-x^2}} \d x \right)^3
\end{align}
For $x \in (-2 + \kappa, 2 - \kappa)$ define,
\beq
\phi_x (E) := \frac{  \Im [ \log ( E - x_T ) ]}{ \Im [ \msc (x_T) ] } .
\eeq
We form now the $2 \times 2$ matrices and $M_1$ and $M_2$ by,
\beq
M_i = \left( \begin{matrix} V_i ( \phi_x, \phi_x) & V_i ( \phi_x, \phi_y ) \\ V_i ( \phi_x , \phi_y ) & V_i ( \phi_y, \phi_y ) \end{matrix}\right)
\eeq
and the third-order homogeneous polynomial,
\beq
P ( \xi_1, \xi_2 ) = P_1 ( \xi_1, \xi_2, \phi_x, \phi_y ).
\eeq
We define also,
\beq
e_x = e ( \phi_x), \qquad e_y = e ( \phi_y), \qquad e_{i, x} = e_i ( \phi_x) , \qquad e_{i, y} = e_i ( \phi_y),
\eeq
where the functionals $e (f)$ and $e_i (f)$ are  as in \eqref{eqn:main-expect-correct}, for $i=1, 2$. Note that it is clear that,
\beq
| e_{i, x}| + | e_{i, y} | \leq C
\eeq
for $i=1, 2$. 
 We also denote,
\beq
\xi = ( \xi_1, \xi_2)^T \in \rr^2 .
\eeq 
The following is an immediate consequence of Theorem \ref{thm:wig-smooth-lss}. 
\bep \label{prop:meso-lss}
Let $\eps >0$. Then for $|\xi_1 | + | \xi_2 | \leq N^{1/10}$ we have
\begin{align}
 & \left| \ee \left[ \exp \left( \i \xi_1 ( Y (x) - e_x ) + \i \xi_2 ( Y(y) - e_y ) \right) \right] - \exp\left( - \xi^T(M_1 + M_2) \xi + \i N^{-1/2} P( \xi ) \right) \right| \notag\\
\leq & \,C \frac{N^{\eps}}{NT} (1 + | \xi |^7 ).
\end{align}
\eep
This motivates the following definitions. Fix an $\eps_3 >0$ satisfying
\beq \label{eqn:meso-eps3}
\frac{1}{ 100} > \eps_3 >0
\eeq
and let,
\beq
\psi_1 ( \xi ) = \exp \left( - \xi^T M_1 \xi \right)
\eeq
and
\beq
\psi_2 ( \xi ) = \exp \left( - \xi^T M_2 \xi - \xi_1^2 A_1 - \xi_2^2 A_1 + N^{-1/2} \i P ( \xi ) \right) \chi ( \xi_1 / N^{\eps_3} ) \chi ( \xi_2 / N^{\eps_3 } ).
\eeq
Note that $\psi_1$ is the ``universal component'' of the asymptotic joint distribution of $(Y(x), Y(y))$, and $\psi_2$ is the non-universal component, which also includes our Gaussian convolution. Note also that the matrix $M_1$ is positive semi-definite since,
\beq
\xi^T M_1 \xi = V_1 ( \xi_1 \phi_x + \xi_2 \phi_y , \xi_1 \phi_x + \xi_2 \phi_y  ) \geq 0.
\eeq
Whether or not $\psi_1$ decays will not play any role due to the estimates we have derived for $\hat{H}$ above. 

We now choose $A_1$ sufficiently large that,
\beq \label{eqn:meso-A1}
| \xi^T M_2 \xi^T + \xi_1^2 A_1 + \xi_2^2 A_1 - N^{-1/2} \i P ( \xi ) | \geq 10( \xi_1^2 + \xi_2^2)
\eeq
for all $\xi$ satisfying $|\xi_1| + | \xi_2 | \leq N^{1/20}$.  Note that this is possible since the entries of $M_2$ are bounded independently of $T$ and $|P (\xi ) | \leq C N^{-1/4}$ for such $\xi$. 
Now, let $q(z_1, z_2)$ be the inverse Fourier transform of $\psi_2$, that is
\beq
\psi_2 ( \xi_1, \xi_2 ) = \int_{\rr^2} \e^{ \i ( \xi_1 z_1 + \xi_2 z_2 ) } q ( z_1, z_2) \d z_1 \d z_2.
\eeq
 The function $q(z_1, z_2)$ may not be a probability distribution, but it has sufficient properties suitable for our purposes, summarized in the following.
\bep \label{prop:q-est}
The function $q(z_1, z_2)$ is real-valued and satisfies,
\beq
\int_{\rr^2} q(z_1, z_2) \d z_1 \d z_2 = 1 .
\eeq
For any $k >0$, there is a $C>0$ so that,
\beq
| q(z_1, z_2) | \leq C \frac{1}{1 + | z_1|^k + |z_2|^k }.
\eeq
\eep
\proof The function is real-valued since its Fourier transform satisfies $\bar{ \psi}_2 (\xi ) = \psi_2 ( - \xi )$.  The integral of $q$ is $1$ since $\psi_2 ( 0,0) = 1$.   
The decay estimates follow from the fact that the choices \eqref{eqn:meso-eps3} and \eqref{eqn:meso-A1} implies that,
\beq
\| ( \del^k_{ \xi_1}  \del^l_{\xi_2} ) \psi_2 \|_p \leq C_{l,k,p}
\eeq
for any $p, k, l$. \qed

Introduce also $p(z_3, z_4)$, a Gaussian distribution with covariance $M_1$, and define the function $F(z_1, z_2)$ by
\beq
F (z_1, z_2) := \int p(z_3, z_4 ) \ee[ F_{1, x, x} (z_1 + z_3 + e_{1,x}) F_{1, y, y} (z_2 + z_4 +e_{1,y}) ] \d z_3 \d z_4.
\eeq
Eventually we will see that the function $F(z_1, z_2)$ well-approximates the covariance of two resolvents of the GOE. First, we require the following elementary estimate on the distribution $p(z_3, z_4)$. 
\bel
For any $D>0$ and $\eps >0$ there exists $C>0$ so that,
\beq
\int_{ |z_3 | > N^{\eps} } p (z_3, z_4) \d z_3 \d z_4 \leq C N^{-D}
\eeq
and an analogous estimate for $z_4$.
\eel
\proof It suffices to find an upper bound on the diagonal entries of $M_1$.  We note that for $|x| \leq 2 - \kappa$ that
\beq
V_1 ( \phi_x, \phi_x) \leq C \left( 1 + \int_{ |u|, |v| \leq 2- \kappa/2} \left( \frac{ \phi_x (u ) - \phi_x (v) }{u-v} \right)^2 \d u \d v \right).
\eeq 
The integral on the RHS can be estimated via the same method as in the proof of Lemma 5.19 of \cite{landon2021single}. This results in an upper bound of $\O ( \log(N) )$, which is sufficient for our purposes. \qed

The following result replaces the kernel appearing in the represenation \eqref{eqn:meso-var-rep} for the variance associated to $\varphi_k$ with the function $F(z_1, z_2)$ defined above. Later, the function $F$ will be related to the covariance kernel of the GOE. 
\bep \label{prop:meso-a1}
Let $ \eps >0$ be given, and let $\eps_3 >0$ and $\eps_2 >0$ be as above.  We have,
\begin{align}
 & \ee^{(\tilZ)}[ N^2 \ee^{(H)}[  \Im [ (m_N - \msc ) ( x+Z_1 + \i \eta ) ] \Im [(m_N - \msc ) ( y+Z_2 + \i \eta )]] ]  \notag \\
 = & \int q(z_1, z_2) F(z_1 + e_{2,x}, z_2+e_{2,y}) \d z_1 \d z_2 \notag\\
 + & N^{\eps+2\eps_2} \O (\eta^{-1}  T^{-1})
\end{align}
\eep
\proof By applying Proposition \ref{prop:meso-homog-1} with $\hatx = x + Z_1$ and $\haty = y + Z_2$, we first see that 
\begin{align}
 & \ee^{(\tilZ)}[ N^2 \ee^{(H)}[  \Im [ (m_N - \msc ) ( x+Z_1 + \i \eta ) ] \Im [(m_N - \msc ) ( y+Z_2 + \i \eta )]] ]  \notag \\
 = & \ee^{(\tilZ, Y)} N^2 \ee^{(\mu)} [ F_{1, x, x} ( N Z_1+ Y(x) ) F_{1, y, y} ( NZ_2 + Y(y) ) ] \notag\\
 + & N^{\eps} \O ( \eta^{-1}  T^{-1}).
 \end{align}
 It is straightforward to conclude from the rigidity estimates that $|Y(x)| + |Y(y) | \leq N^{\eps}$ with overwhelming probability, for any $\eps >0$. Therefore,  we have 
\begin{align}
& \ee^{(\tilZ, Y)} N^2 \ee^{(\mu)} [ F_{1, x, x} ( NZ_1+ Y(x) ) F_{1, y, y} (N Z_2 +  Y(y) ) ] \notag\\ 
= &\ee^{(\tilZ, Y)} N^2 \ee^{(\mu)} [ F_{2, x, x} (NZ_1+ Y(x) ) F_{2, y, y} ( NZ_2 + Y(y) ) ] + \O (N^{-100}).
\end{align}
Recall that $H(z_1, z_2)$ is defined by
\beq
H(z_1, z_2) = \ee^{(\mu)} [ F_{2, x, x} (z_1) F_{2, y, y} ( z_2 ) ] .
\eeq
By Fourier duality,
\begin{align}
 & \ee^{(\tilZ, Y)} N^2 \ee^{(\mu)} [ F_{2, x, x} (NZ_1+ Y(x) ) F_{2, y, y} ( NZ_2 + Y(x) ) ] \notag\\
= & N^2 \int_{\rr^2} \hat{H} ( \xi_1, \xi_2 )\e^{ - A_1 ( \xi_1^2 + \xi_2^2 )} \psi ( \xi_1, \xi_2 ) \e^{ \i \xi_1 e_x + \i \xi_2 e_y } \d \xi_1 \d \xi_2,
\end{align}
where we denoted
\beq
\psi ( \xi_1, \xi_2) = \ee \left[ \exp \left( \i \xi_1 ( Y(x) - e_x ) + \i \xi_2 ( Y(y) - e_y ) \right) \right].
\eeq
By Proposition \ref{prop:meso-lss} and the estimates for $\hat{H}$ from Corollary \ref{cor:H-fourier} we have,
\begin{align}
& N^2 \int_{\rr^2} \hat{H} ( \xi_1, \xi_2 )\e^{ - A_1 ( \xi_1^2 + \xi_2^2 )} \psi ( \xi_1, \xi_2 ) \e^{ \i \xi_1 e_x + \i \xi_2 e_y } \d \xi_1 \d \xi_2 \notag\\
= &N^2 \int_{\rr^2} \hat{H} ( \xi_1, \xi_2 ) \psi_1 ( \xi_1, \xi_2 ) \psi_2 ( \xi_1, \xi_2 ) \e^{ \i \xi_1 e_x + \i \xi_2 e_y } \d \xi_1 \d \xi_2 +  N^{10 \eps+2\eps_2} \O ( \eta^{-2} (N T)^{-1} )  \notag\\
= & N^2 \int_{\rr^4} q (z_1, z_2 ) p (z_3, z_4) H (z_1+z_3 + e_x, z_2+z_4 + e_y) ] \d z_1 \d z_2 \d z_3 \d z_4 + N^{10 \eps+2 \eps_2} \O ( \eta^{-2} (N T)^{-1} ) ,
\end{align}
where we used Fourier duality again in the last line. By the decay of $p$ and $q$ we see that,
\begin{align}
& N^2 \int_{\rr^4} q (z_1, z_2 ) p (z_3, z_4) H (z_1+z_3 + e_x, z_2+z_4+ e_y ) ] \d z_1 \d z_2 \d z_3 \d z_4  \notag\\
= & N^2 \int_{\rr^4} q (z_1, z_2 ) p (z_3, z_4)  \ee^{(\mu)} [ F_{1, x, x} (z_1+z_3 + e_x) F_{1, y, y} ( z_2+z_4 + e_y ) ]\d z_1 \d z_2 \d z_3 \d z_4 \notag\\
+ & \O (N^{-100} ).
\end{align}
But, by definition,
\begin{align}
& N^2 \int_{\rr^4} q (z_1, z_2 ) p (z_3, z_4)  \ee^{(\mu)} [ F_{1, x, x} (z_1+z_3 + e_x) F_{1, y, y} ( z_2+z_4 + e_y ) ]\d z_1 \d z_2 \d z_3 \d z_4 \notag\\
= & N^2 \int_{\rr^2} q (z_1, z_2 ) F(z_1 + e_{2,x} , z_2 + e_{2, y}  ) \d z_1 \d z_2.
\end{align}
This yields the claim. \qed

Let $Y^{G'}(x) := X^{G'} (x, T)$ where $G'$ is an auxilliary Gaussian matrix, coupled to the process $\mu_i (t)$ in the same manner as $\lambda_i (t)$. 

By Theorem \ref{thm:wig-smooth-lss} we have,
\bep
Let $G'$ be a GOE matrix. Then,
\beq
\left| \ee \left[ \exp \left( \i \xi_1 ( Y(x) - e_{1,x} ) + \i \xi_2 ( Y(y) - e_{1,y} ) \right) \right] - \exp \left( - \xi^T M_1 \xi \right) \right| \leq C \frac{N^{\eps}}{NT} (1 + |\xi|^7 ).
\eeq
\eep

Consider now repeating the above arguments with the matrix,
\beq
G_t := \e^{-t/2} G' + \sqrt{ 1  - \e^{-t} } G
\eeq
where $G'$ is an independent Gaussian matrix, in the place of $H_t$ as defined in \eqref{eqn:meso-hh1} (i.e., replace the matrix $W$ there with a GOE matrix). In this case, we do not introduce the auxilliary variables $Z_1$ and $Z_2$, and instead of introducing the function $\psi_2$, it is replaced by the constant function $\psi_2 ( \xi ) = 1$, and $q$ is just the $\delta$ function at $0$ in $\rr^2$.  Then, we see that the analog of Proposition \ref{prop:meso-a1} will be the following, obtained by almost exactly the same arguments.
\bep \label{prop:meso-a2}
Let $G'$ be a Gaussian matrix. Let $\eps >0$ and $\eps_1 >0$. Assume that,
\beq
|z_1| + |z_2| \leq N^{\eps_1/100}
\eeq
Then,
\begin{align}
 & N^2 \ee^{(G')}[  \Im [ (m_N - \msc ) ( x+z_1/N + \i \eta ) ] \Im [(m_N - \msc ) ( y+z_2/N + \i \eta )]] ]  \notag\\
 = & N^2 F(z_1, z_2) + N^{\eps+\eps_1} \O ( \eta^{-1} T^{-1} ).
\end{align}
\eep
The above will allow us to begin estimating the covariance. We have the following for the Gaussian ensembles.
\bep \label{prop:meso-a3}
Let $G$ be a matrix from the Gaussian ensemble. Then for any $1 \geq \eta \geq N^{-2}$ and $\eps >0$ and $x, y \in (-2+\kappa, 2- \kappa )$ we have,
\beq
N^2 \left|  \Cov_G ( \Im [ m_N (x + \i \eta )], \Im [ m_N (y + \i \eta ) ] ) \right| \leq C N^{\eps} \frac{1}{ \eta ( |x-y| + \eta ) }
\eeq
Moreover, 
\beq
N \left| \ee^G[ \Im [( m_N - \msc) ( x + \i \eta )  ] \right| \leq C
\eeq
\eep
\proof The first estimate follows from Theorem \ref{thm:main-gaussian}. The second follows from writing the Poisson kernel as 
\beq
P_\eta (x) = \chi_1 (x) P_\eta (x) + (1- \chi_1(x) ) P_\eta (x) = f_1(x) + f_2(x).
\eeq
where $\chi_1$ is a function that is $1$ on $(-2+\kappa/2, 2- \kappa/2)$ and $0$ outside of $(-2-\kappa/4 , 2- \kappa/4)$. From Theorem \ref{thm:gaussian-bulk} we have
\beq
\left| \ee[ \tr f_1 (G)] - N \int f_1(x) \rhosc (x) \right| \leq C.
\eeq
From Proposition \ref{prop:gaussian-edge} we have,
\beq
\left| \ee[ \tr f_2 (G)]- N \int f_2 (x) \rhosc (x) \right| \leq C.
\eeq
The claim follows. \qed

We can now conclude the main technical estimate of this section. 
\bet \label{thm:meso-main}
For any $\eps >0$ the following holds. Let $\eta = 2^{-k}$ and $k$ satisfy,
\beq
\mfb \log_2 (N) \leq k \leq \log_2 (N \mfb^{-1} ).
\eeq
We have
\begin{align}
& \left|  \ee^{(\tilZ)}[ N^2 \ee^{(H)}\big[  \Im [ (m_N - \msc ) ( x+Z_1 + \i \eta ) ] \Im [(m_N - \msc ) ( y+Z_2 + \i \eta )] \big] \right| \notag\\
\leq & C N^{\eps} \left( \frac{1}{\eta ( |x-y| + \eta )} + \frac{1}{ \eta T} \right). 
\end{align}
\eet
\proof Thanks to Proposition \ref{prop:meso-a1}, it suffices to estimate the function $F(z_1, z_2)$. For $|z_1| + |z_2| \leq N^{\eps/1000}$ we have, by Proposition \ref{prop:meso-a2},
\begin{align}
\left| N^2 F(z_1, z_2) \right| & \leq N^2 \left|  \Cov_G ( \Im [ m_N (x+N^{-1} z_1 + \i \eta )], \Im [ m_N (y +N^{-1}z_2 + \i \eta ) ] ) \right| \notag\\
&+ N^2\bigg|  \ee^G \big[ \Im[(m_N - \msc) ( x + N^{-1} z_1 + \i \eta )]\big]  \notag\\
&\times \ee^G \big[ \Im[(m_N - \msc) ( y + N^{-1} z_2 + \i \eta )]\big]  \bigg| \notag\\
&+ N^{\eps} \frac{1}{ \eta T}.
\end{align}
The claim now follows from Proposition \ref{prop:meso-a3}, as well as the estimates in Proposition \ref{prop:q-est}.  \qed

\noindent{\bf Proof of Theorem \ref{thm:meso-var}}. By \eqref{eqn:meso-var-rep} and Proposition \ref{prop:meso-lp-supp} we have,
\begin{align}
 & \Var ( \tr \varphi_k (H) ) \notag\\
 \leq & \int_{I_{\kappa/2}^2} \left| f_k (x) f_k (y) N^2 \ee^{(\tilZ ) } \big[ \Im [ (m_N - \msc) ( x +Z_1+ \i \eta ) ]  \Im [ (m_N - \msc) ( y+Z_2 + \i \eta ) ]\big] \right|  \d x \d y \notag\\
 + & C N^{-100} \| \varphi \|_2^2
\end{align}
where 
\beq
I_{a} = (-2+a, 2-a)
\eeq
and the support of $\varphi$ is in $(-2+\kappa, 2- \kappa)$.  For $x, y \in I_{\kappa/2}$ we see from Theorem \ref{thm:meso-main} that
\begin{align}
& \int_{I_{\kappa/2}^2} \left| f_k (x) f_k (y) N^2 \ee^{(\tilZ ) } [ \Im [ (m_N - \msc) ( x +Z_1+ \i \eta ) ] ( \Im [ (m_N - \msc) ( y+Z_2 + \i \eta ) ] \right|  \d x \d y  \notag\\
\leq & N^{\eps} \int_{I_{\kappa/2}^2} |f(x) f(y) |  \left( \frac{1}{ \eta ( |x-y| + \eta )} + \frac{1}{ \eta T} \right) \d x \d y
\end{align}
Clearly,
\begin{align}
\frac{1}{ \eta T} \int_{I_{\kappa/2}^2} |f_k(x) f_k(y) | \d x \d y \leq & \frac{1}{ \eta T} \int_{I_{\kappa/2}^2} |f_k(x)|^2 + |f_k(y)|^2 \d x \d y  \notag\\
 \leq & \frac{C}{ \eta T} \| f_k \|_2^2 \leq \frac{C}{ \eta T} \| \varphi_k \|_2^2
\end{align}
where in the final inequality we used Proposition \ref{prop:meso-a4}.   Since,
\beq
\int_{|x|<5} \frac{1}{ \eta ( |x-y| + \eta ) } \d x \leq  C \frac{| \log ( \eta ) |}{ \eta}
\eeq
for any $y$, we see that 
\begin{align}
\int_{I_{\kappa/2}^2} \frac{|f_k (x) f_k (y) |}{ (x-y)^2 + \eta^2 } \d x \d y \leq & \int_{I_{\kappa/2}^2} \frac{|f_k (x)|^2 }{ (x-y)^2 + \eta^2 } \d x \d y  \notag \\
& \leq C N^{\eps} \eta^{-1} \| f_k \|_2^2 \leq C N^{\eps} \eta^{-1} \| \varphi_k \|_2^2.
\end{align}
The result follows. \qed

\section{Large frequencies} \label{sec:global}

In this section we prove Theorems \ref{thm:global-cov} and the variance and expectation estimates of Theorem \ref{thm:wig-smooth-lss} via resolvent expansions.  Fix $\tau >0$ and let us define,
 \beq
 \D = \D_{\tau, 1} \cup \bar{\D}_{\tau, 1} ,
 \eeq
 where $\D_{\tau, 1}$ is defined in \eqref{eqn:Dtau-def}. Recall the definition of the control parameter $\Psi(z)$ in \eqref{eqn:Psi-def}.

 In the rest of this section we will use the notation,
\beq
z = x \pm \i y , \qquad w = u \pm \i v
\eeq
where $y, v >0$.
 
 The main tool is the following self-consistent equation for expectations of products of Stieltjes transforms.
\bep \label{prop:cov-self}
For $H$ a real symmetric Wigner matrix and $z, w \in \D$ we have for any $\eps >0$,
\begin{align}
& (z+ 2 \msc (z) )\ee[ m_N (z) ( m_N (w) - \ee[ m_N (w) ] ) ] \notag\\
= & - \frac{2}{N^2}\del_w \frac{ \msc (z) - \msc (w) }{ z -w } + \frac{4 s_3 ( \msc(z)^2 \msc'(w) + \msc(z) \msc(w) \msc'(w)) }{N^{5/2}} \notag  \\
- & \frac{2 s_4 \msc'(w) \msc (w) \msc(z)^2}{N^2} \notag\\
+ & N^{\eps} \O ( N^{-3} ( y^{-2} v^{-1} + (y+v)^{-1} v^{-2} )) \notag \\
+ & N^{\eps} \O ( N^{-3/2} v^{-1} ( \Psi (z)^2 \Psi (w) + \Psi(z) \Psi(w)^2 ) +N^{-2} \Psi(z) \Psi (w) v^{-1} )
\end{align}
\eep

The proof of this proposition is via cumulant expansions, and is given over the course of Section \ref{sec:global-self}.  

\subsection{Self-consistent equation for covariance} \label{sec:global-self}

By the cumulant expansion we have for any $\eps >0$ that,
\begin{align} \label{eqn:cov-cum-exp}
 & \frac{1}{N^2} \sum_{i, j} z \ee[ G_{ii} (z) ( G_{jj} (w) - \ee[ G_{jj} (w) ] ) ] \notag\\
 =& \frac{1}{N^2} \sum_{i, j, a} \ee[ H_{ia} G_{ai} (z) ( G_{jj} (w) - \ee[ G_{jj} (w) ] ) ] \notag\\
 = & \frac{1}{N^2} \sum_{i, j, a} \sum_{k=2}^7 \frac{ (1+ \delta_{ia} )^{k-1} s_k }{(k-1)!N^{k/2} } \ee[ \del_{ia}^{k-1} (G_{ai} (z) ( G_{jj} (w) - \ee[ G_{jj} (w) ] ) )] \notag \\
 +& \O (N^{\eps-3} ).
\end{align}
Here, the notation $\del_{ia}$ means the derivative with respect to the matrix entry $(i, a)$.  Since we have symmetric matrices, note that this involves the derivative both with respect to the entry $(i,a)$ and $(a, i)$. 

For a treatment of the error term, see the discussion in \cite{lee2018local} (the expansion with error term appears there as Lemma 3.2); there is no substantial difference from the treatment appearing there and here.  We will refer to the $k$th term in the sum \eqref{eqn:cov-cum-exp} as the $k$th-order term. Each of these terms is expanded and dealt with in each of the next few subsubsections. We will also require expressions for derivatives of resolvent entries with respect to matrix entries. These are collected in Appendix \ref{a:derivs}.

\subsubsection{Second order term}

\bel \label{lem:cov-2o}
We have for any $\eps >0$,
\begin{align}
& \frac{1}{N^3} \sum_{ija} (1+\delta_{ia} ) \ee[ \del_{ia} (G_{ai} (z) ( G_{jj} (w) - \ee[ G_{jj} (w) ] ) )] \notag \\
= &-2 \msc (z) \ee[ ( m_N (z) - \ee[ m_N (z) ] ) ( m_N (w) - \ee[ m_N (w) ] ) ]\\
- & \frac{2}{N^2} \del_w \frac{ \msc (z) - \msc (w) }{ z -w } + \O ( N^{\eps-3} (y^{-2} v^{-1}+ (y+v)^{-1} v^{-2} )
\end{align}
\eel
\proof We have, using \eqref{eqn:delGia} and \eqref{eqn:delGjj},
\begin{align}
& \frac{1}{N^3} \sum_{ija} (1+\delta_{ia} ) \ee[ \del_{ia} (G_{ai} (z) ( G_{jj} (w) - \ee[ G_{jj} (w) ] ) )] \notag \\
= & - \frac{1}{N^3} \sum_{ija} \ee[ G_{ii} (z) G_{aa} (z) ( G_{jj} (w) - \ee[ G_{jj} (w) ] ) ] \label{eqn:cov-1o1} \\
 -& \frac{1}{N^3} \sum_{ija} \ee[ G_{ia} (z)^2 ( G_{jj} (w) - \ee[ G_{jj} (w) ] ) ] \label{eqn:cov-1o2} \\
 - & \frac{2 }{N^3} \sum_{ija} \ee[ G_{ia} (z) G_{ja} (w) G_{ji} (w) ] \label{eqn:cov-1o3}
\end{align}
The first term \eqref{eqn:cov-1o1} equals,
\begin{align}
 &- \frac{1}{N^3} \sum_{ija} \ee[ G_{ii} (z) G_{aa} (z) ( G_{jj} (w) - \ee[ G_{jj} (w) ] ) ] \notag \\
 =& - \ee[ m_N (z)^2 (m_N (w) - \ee[ m_N (w) )] \notag \\
 =&  -2 \ee[ m_N (z) ] \ee[ ( m_N (z) - \ee[ m_N (z) ] ) ( m_N (w) - \ee[ m_N (w) ] ) ] + \O ( N^{\eps-3} y^{-2} v^{-1} ) \notag \\
 = & -2 \msc (z) \ee[ ( m_N (z) - \ee[ m_N (z) ] ) ( m_N (w) - \ee[ m_N (w) ] ) ] + \O ( N^{\eps-3} y^{-2} v^{-1} ).
\end{align}
In the second to last line we wrote $m_N (z) = ( m_N (z) - \ee[ m_N (z) ] ) + \ee[ m_N (z) ]$, expanded the square and used the local law \eqref{eqn:ll}. In the last line we used again the local law \eqref{eqn:ll}. 

The term \eqref{eqn:cov-1o2} equals,
\begin{align}
&-\frac{1}{N^3} \sum_{ija} \ee[ G_{ia} (z)^2 ( G_{jj} (w) - \ee[ G_{jj} (w) ] ) ]  = - \frac{1}{N} \ee[ \del_z m_N (z) ( m_N (w) - \ee[ m_N (w) ] ) ] \notag \\
=& - \frac{1}{N} \ee[ \del_z (m_N (z)- \msc (z) ) ( m_N (w) - \ee[ m_N (w) ] ) ] = \O ( N^{\eps} N^{-3} y^{-2} v^{-1} )
\end{align}
where we used the local law \eqref{eqn:ll} and the Cauchy integral formula to obtain the last estimate.  For the term \eqref{eqn:cov-1o3} we have,
\begin{align}
- & \frac{2 }{N^3} \sum_{ija} \ee[ G_{ia} (z) G_{ja} (w) G_{ji} (w) ] = - 2 \del_w \frac{1}{N^3} \ee[ \tr G(z) G(w) ].
\end{align}
We have the identity,
\beq
\frac{1}{N} \tr G(z) G(w) = \frac{m_N (z) - m_N (w) }{z-w} .
\eeq
If $z$ and $w$ are in opposite half-planes, then since $|z-w| \geq y + v$ we have by the local law \eqref{eqn:ll} with overwhelming probability,
\begin{align}
 & \del_w \frac{ m_N (z) - m_N (w) }{z-w} = \frac{ m_N(z) - m_N (w) }{ (z-w)^2} - \frac{ \del_w m_N (w) }{ (z-w)} \notag \\
= & \frac{ \msc (z) - \msc (w) }{ (z-w)^2} - \frac{ \del_w \msc (w) }{ z-w} + \O ( N^{\eps-1}(  y^{-2} v^{-1} + (y+v)^{-1}  v^{-2} )
\end{align}
and the main term in the last line equals $\del_w \frac{ \msc (z) - \msc (w) }{ z-w}$.  If $z$ and $w$ are in the same half-plane we consider first the case $ 10|z-w| \geq \max \{ y, v\} $. In this case we use the same method as above to derive the same estimate. In the case that $10 |z-w| \leq \max \{ y, v\}$ we see that $v \leq 2 y \leq 4 v$. We have the identity,
\beq
\frac{ m_N (z) - m_N (w) }{ z-w } = - \int_0^1 m_N' (z + s (w-z) ) \d s .
\eeq
For the integrand we have, $m_N'(z+ s (w-z) ) = \msc'(z + s (w-z) ) + \O (N^{\eps-1} y^{-2} )$ by the local law with overwhelming probability.  Summarizing for \eqref{eqn:cov-1o3} we have derived that with overwhelming probability,
\begin{align} \label{eqn:GG-est}
-\frac{2 }{N^3} \sum_{ija}  G_{ia} (z) G_{ja} (w) G_{ji} (w)   =& - 2\del_w \frac{ \msc (z) - \msc (w) }{ z -w } \notag \\
+& \O(N^{\eps-3}( y^{-2} v^{-1} + (y+v)^{-1}  v^{-2}) ) .
\end{align}
This yields the claim. \qed

\subsubsection{Third order term}

In the proof below, where possible, we will simplify errors that arise using
\beq \label{eqn:Psi-simple}
\frac{1}{N} \leq \Psi (E \pm \i \eta ) \leq C \frac{1}{ (N \eta )^{1/2}}
\eeq
for $E \pm \i \eta \in \D$, to absorb them into error terms that have arisen already in the proof.

\bel \label{lem:cov-3o}
We have,
\begin{align}
& \frac{1}{N^{3+1/2} } \sum_{ija}  (1+ \delta_{ia} )^2 \ee[ \del_{ia}^{2} (G_{ai} (z) ( G_{jj} (w) - \ee[ G_{jj} (w) ] ) )] \notag\\
= & 8\frac{  \msc(z)^2 \msc'(w) + \msc(z) \msc(w) \msc'(w) }{N^{5/2}} \notag \\
+ & N^{\eps} \O (N^{-3} y^{-2} v^{-1} + N^{-3} v^{-2} (y+v)^{-1} ) \notag \\
+ & N^{\eps} \O ( N^{-3/2} v^{-1} ( \Psi (z)^2 \Psi (w) + \Psi(z) \Psi(w)^2 ) +N^{-2} \Psi(z) \Psi (w) v^{-1} ).
\end{align}
\eel
\proof We have, using \eqref{eqn:del2Gia}, \eqref{eqn:delGiadelGjj} and \eqref{eqn:del2Gjj},
\begin{align}
& \frac{1}{N^{3+1/2} } \sum_{ija}  (1+ \delta_{ia} )^2 \ee[ \del_{ia}^{2} (G_{ai} (z) ( G_{jj} (w) - \ee[ G_{jj} (w) ] ) )]  \notag\\
= & \frac{ 2}{N^{2+1/2}} \sum_{ia} \ee[ G_{ia} (z)^3 ( m_N (w) - \ee[ m_N (w) ] ) ] \label{eqn:cov-3o1} \\
+ & \frac{6}{N^{2+1/2}} \sum_{ia} \ee[ G_{ia} (z) G_{ii} (z) G_{aa} (z) (m_N (w) - \ee[ m_N (w) ]  ) ] \label{eqn:cov-3o2} \\
+ & \frac{4}{N^{3+1/2} } \sum_{ia} \ee[ G_{ii} (z) G_{aa} (z) \del_w G_{ia} (w) ] \label{eqn:cov-3o3} \\
+ & \frac{4}{N^{3+1/2} } \sum_{ia} \ee[ G_{ia} (z)^2 \del_w G_{ia} (w) ] \label{eqn:cov-3o4} \\
+ & \frac{4}{N^{3+1/2}} \sum_{ia} \ee[ G_{ia} (z) G_{ai} (w) \del_w G_{ai} (w) ] \label{eqn:cov-3o5} \\
+ & \frac{4}{N^{3+1/2}} \sum_{ia} \ee[ G_{ia} (z) G_{ii} (w) \del_w G_{aa} (w) ] \label{eqn:cov-3o6} .
\end{align}
The first term \eqref{eqn:cov-3o1} is,
\begin{align}
&\frac{ 2}{N^{2+1/2}} \sum_{ia} \ee[ G_{ia} (z)^3 ( m_N (w) - \ee[ m_N (w) ] ) ] \notag\\
= & \frac{2}{N^{2+1/2} } \sum_{i} \ee[ ( G_{ii}(z)^3 - \msc(z)^3 )  ( m_N (w) - \ee[ m_N (w) ] ) ] + \O (N^{\eps} N^{-3/2} \Psi(z)^3 v^{-1} ) \notag \\
= & N^{\eps} \O ( N^{-7/2} y^{-1} v^{-1} +N^{-3/2} \Psi(z)^3 v^{-1} )
\end{align}
In the second line we used \eqref{eqn:entry-ll} to estimate the off-diagonal Green's function elements and also applied the local law \eqref{eqn:ll}. In the last line we expanded out $(G_{ii} - \msc + \msc )^3$ and used the local law \eqref{eqn:ll} to estimate the term linear in $G_{ii}$ and then the entry-wise law \eqref{eqn:entry-ll} to estimate the higher powers.

For the term \eqref{eqn:cov-3o2} we first introduce the notation,
\beq
\be = N^{-1/2} (1, \dots, 1)^T \in \rr^N ,
\eeq
i.e., the constant unit vector in $\rr^N$.  Then, writing $G_{ii} = (G_{ii} - \msc) + \msc$ and expanding (and doing the same for $G_{aa}$) we find for \eqref{eqn:cov-3o2}, 
\begin{align}
& \frac{1}{N^{2+1/2}} \sum_{ia} \ee[ G_{ia} (z) G_{ii} (z) G_{aa} (z) (m_N (w) - \ee[ m_N (w) ]  ) ] \notag \\
=& \frac{1}{N^{3/2}} \msc(z)^2\ee[ \be^T G(z) \be (m_N (w) - \ee[ m_N (w) ] ) ] \notag \\
+ &\frac{ 2 \msc (z)}{N^2} \sum_{i} \ee[ (G_{ii}(z) - \msc (z) ) ( G (z) \be )_i (m_N (w) - \ee[ m_N (w) ] ) ] \notag\\
+ & \O (N^{\eps} N^{-3/2} \Psi(z)^3 v^{-1} ) .
\end{align}
By the isotropic local law \eqref{eqn:iso} we see
\beq
 \frac{1}{N^{3/2}} \msc(z)^2\ee[ \be^T G(z) \be (m_N (w) - \ee[ m_N (w) ] ) ] = N^{\eps} \O ( N^{-5/2} \Psi(z) v^{-1} ).
\eeq
By the isotropic local law \eqref{eqn:iso} we see with overwhelming probability that,
\beq \label{eqn:cov-3o-iso}
(G(z) \be )_i = N^{-1/2} \msc (z) + N^{\eps} \O ( \Psi(z))
\eeq
and so,
\begin{align}
&\frac{ 2 \msc (z)}{N^2} \sum_{i} \ee[ (G_{ii}(z) - \msc (z) ) ( G (z) \be )_i (m_N (w) - \ee[ m_N (w) ] ) ] \notag \\
=& \frac{2 \msc(z)^2}{N^{3/2} } \ee[ (m_N (z) - \msc(z) ) ( m_N (w) - \ee[ m_N (w)] ) ] + N^{\eps} \O ( N^{-2} \Psi(z)^2 v^{-1})\notag\\
= & \O ( N^{-7/2} y^{-1} v^{-1} +N^{-2} \Psi(z)^2 v^{-1} ).
\end{align}
This completes estimating the term \eqref{eqn:cov-3o2}. We write \eqref{eqn:cov-3o3} as,
\begin{align}
 & \frac{1}{N^{3+1/2} } \sum_{ia} \ee[ G_{ii} (z) G_{aa} (z) \del_w G_{ia} (w) ]  \notag\\
 = & \frac{\msc(z)^2}{N^{5/2}} \ee[ \del_w \be^T G (w) \be ] + \frac{2 \msc (z) }{N^{3}} \sum_i \ee[ (G_{ii} (z)-\msc(z) ) \del_w (G (w) \be )_i ] \notag\\
+& \frac{1}{N^{7/2}} \sum_{ia} \ee[ (G_{ii}(z) - \msc (z) ) ( G_{aa} (z) - \msc (z) ) \del_w G_{ia} (w) ].
\end{align}
By the isotropic local law \eqref{eqn:iso} and the Cauchy integral formula, we have
\beq
\frac{\msc(z)^2}{N^{5/2}} \ee[ \del_w \be^T G (w) \be ] = \frac{ \msc(z)^2 \msc'(w) }{N^{5/2}} + N^{\eps} \O ( N^{-5/2} \Psi(w) v^{-1} ).
\eeq
By \eqref{eqn:cov-3o-iso} with $z$ replaced by $w$ and the entry-wise local law \eqref{eqn:entry-ll} we have,
\beq
\frac{2 \msc (z) }{N^{3}} \sum_i \ee[ (G_{ii} (z)-\msc(z) ) \del_w (G (w) \be )_i ] = N^{\eps} \O (N^{-3} y^{-1} v^{-1/2} + N^{-2} \Psi(z) \Psi(w) v^{-1} ),
\eeq
where we used a crude bound $| \msc'(w) | \leq Cv^{-1/2}$. Finally, after separating into diagonal and off-diagonal terms we find,
\begin{align}
&\frac{1}{N^{7/2}} \sum_{ia} \ee[ (G_{ii}(z) - \msc (z) ) ( G_{aa} (z) - \msc (z) ) \del_w G_{ia} (w) ] \notag\\
=& N^{\eps} \O( N^{-3/2} \Psi(z)^2 \Psi(w) v^{-1}  + N^{-5/2} \Psi(z)^2 ( \Psi(w)v^{-1} + v^{-1/2} ) ) ,
\end{align}
again using a crude bound $|\msc'(w) | \leq C v^{-1/2}$. 
This completes \eqref{eqn:cov-3o3}. We now turn to \eqref{eqn:cov-3o4}. Clearly,
\begin{align}
\frac{4}{N^{3+1/2} } \sum_{ia} \ee[ G_{ia} (z)^2 \del_w G_{ia} (w) ] = \frac{4}{N^{7/2}} \sum_{i} \ee[ G_{ii}(z)^2 \del_w G_{ii} (w) ] + N^{\eps} \O ( N^{-3/2} \Psi(z)^2 \Psi(w) v^{-1} ),
\end{align}
and 
\begin{align} \label{eqn:cov-3o-diag-1}
& \frac{1}{N^{7/2}} \sum_{i} \ee[ G_{ii}(z)^2 \del_w G_{ii} (w) ] \notag\\
=& \frac{ \msc(z)^2 \msc'(w) }{N^{5/2} } + N^{\eps} \O ( N^{-7/2} ( y^{-1} v^{-1/2} + v^{-2} ) + N^{-5/2} ( \Psi(z) \Psi(w) v^{-1} ),
\end{align}
by expanding $G_{ii}(z) = (G_{ii} (z) - \msc(z) ) + \msc(z)$ (and similarly for $G_{ii}(w)$), using the local law \eqref{eqn:ll} to estimate the terms linear in $G_{ii}$, and the entrywise local law \eqref{eqn:entry-ll} for the terms quadratic and cubic in $G_{ii}$. We also used the crude bound $|\msc'(w) | \leq v^{-1/2}$. 

This completes the treatment of \eqref{eqn:cov-3o4}. For \eqref{eqn:cov-3o5}, the diagonal terms are handled similarly to \eqref{eqn:cov-3o-diag-1} and we get
\begin{align} \label{eqn:cov-3o-diag-2}
& \frac{1}{N^{7/2}} \sum_{i} \ee[ G_{ii}(z) G_{ii} (w) \del_w G_{ii} (w) ] \notag\\
=& \frac{ \msc(z) \msc (w) \msc'(w) }{N^{5/2} } + N^{\eps} \O ( N^{-7/2} ( y^{-1} v^{-1/2} + v^{-2} )+ N^{-5/2} \Psi (z) \Psi(w) v^{-1} ).
\end{align}
The off-diagonal terms are $N^{\eps} \O (N^{-3/2} \Psi(z) \Psi(w)^2 v^{-1} )$. We write \eqref{eqn:cov-3o6} as,
\begin{align}
& \frac{1}{N^{7/2}} \sum_{ia} \ee[ G_{ia} (z) G_{ii} (w) \del_w G_{aa} (w) ] \notag \\
= & \frac{\msc'(w) \msc (w)}{N^{5/2} } \ee[ \be^T G(z) \be] \notag\\
+ & \frac{ \msc' (w) }{N^3} \sum_i \ee[ ( G(z) \be )_i ( G_{ii} (z) - \msc (z) ) ] + \frac{ \msc(z) }{ N^3} \ee[ (G(z) \be )_i  \del_w (G_{ii } (w) - \msc (w) ) ]  \notag\\
+ & N^{\eps} \O ( N^{-3/2} v^{-1} \Psi(w)^2 \Psi (z) ).
\end{align}
By the isotropic local law \eqref{eqn:iso} the first term is,
\beq
\frac{\msc'(w)  \msc (w)}{N^{5/2} } \ee[ \be^T G(z) \be] = \frac{ \msc'(w) \msc(z) \msc (w) }{N^{5/2}} + N^{\eps} \O (N^{-5/2} \Psi(z) v^{-1/2} ).
\eeq
Using \eqref{eqn:cov-3o-iso} as well as the local law \eqref{eqn:ll} we obtain easily that,
\begin{align}
& \frac{ \msc' (w) }{N^3} \sum_i \ee[ ( G(z) \be )_i ( G_{ii} (z) - \msc (z) ) ] + \sum_i \frac{ \msc(z) }{ N^3} \ee[ (G(z) \be )_i  \del_w (G_{ii } (w) - \msc (w) ) ]  \notag\\
= & N^{\eps} \O ( N^{-7/2} y^{-1} v^{-1/2} + N^{-2} \Psi(z)^2 v^{-1/2} + N^{-7/2} v^{-2} + N^{-2} \Psi(z) \Psi(w) v^{-1} )
\end{align}
Adding up all the errors and simplifying yields the claim. \qed 

\subsubsection{Fourth order term}

\bel \label{lem:cov-4o}
We have
\begin{align} \label{eqn:cov-4o-result}
& \frac{1}{N^4} \sum_{ija} (1 + \delta_{ia} )^3 \ee[ \del_{ia}^3 ( G_{ai} (z) ( G_{jj} (w) - \ee[ G_{jj} (w) ] ) ) ] \notag \\
= &   - \frac{12 \msc'(w) \msc (w) \msc(z)^2}{N^2} + N^{\eps} \O (N^{-3} ( v^{-1} y^{-1} + v^{-2} ) + N^{-2} \Psi(z) \Psi(w) v^{-1})
\end{align}
\eel
\proof We break up the left side of \eqref{eqn:cov-4o-result} into the four terms based on how the three derivatives  are distributed amongst the $G_{ia} (z)$ and $G_{jj} (w)$. 
 We first calculate, using \eqref{eqn:del3Gia},
\begin{align}
& -\frac{1}{N^3} \sum_{ia} (1+ \delta_{ia} )^3 \ee[ ( \del_{ia}^3 G_{ai} (z) ) (m_N (w) - \ee[ m_N (w) ] ) ] \notag \\
= &\frac{1}{N^3} \sum_{ia} \ee[ (6 G_{ia}^4 (z) + 36 G_{ia} (z)^2 G_{ii} (z) G_{aa} (z) ) (m_N (w) - \ee[ m_N (w) ] ) ] \notag\\
+ & \frac{6}{N^3} \sum_{ia} \ee [ G_{ii}^2(z) G_{aa} (z)^2 (m_N (w) - \ee[ m_N (w) ] ) ]. \label{eqn:cov-4o1}
\end{align}
Clearly for the second line,
\begin{align}
 & \frac{1}{N^3} \sum_{ia} \ee[ (6 G_{ia}^4 (z) + 36 G_{ia} (z)^2 G_{ii} (z) G_{aa} (z) ) (m_N (w) - \ee[ m_N (w) ] ) ] \notag \\
= & N^{\eps} \O (N^{-2} v^{-1} \Psi(z)^2 + N^{-3} v^{-1} ).
\end{align}
For the last term of \eqref{eqn:cov-4o1} we have,
\beq
\frac{1}{N^3} \sum_{ia} \ee [ G_{ii}^2(z) G_{aa} (z)^2 (m_N (w) - \ee[ m_N (w) ] ) ] = N^{\eps} \O ( N^{-3} y^{-1} v^{-1} )
\eeq
where we again expanded the resolvent entries around $\msc$ and used the local law \eqref{eqn:ll} to estimate the terms linear in $G_{ii}$ or $G_{aa}$.

The next term is, using \eqref{eqn:del2Gia} and \eqref{eqn:delGjj},
\begin{align} \label{eqn:cov-4o2}
-& \frac{1}{N^4} \sum_{ija} (1+\delta_{ia} )^3 \ee[ ( \del_{ia}^2 G_{ia} (z) ) ( \del_{ia} G_{jj} (w) ) ] \notag\\
= & \frac{4}{N^4} \sum_{ia} \ee[  ( \del_{w} G_{ia} (w) ) G_{ia} (z)^3 ] \notag\\
+ & \frac{12}{N^4} \sum_{ia} \ee[  ( \del_{w} G_{ia} (w) ) G_{ii} (z) G_{aa} (z) G_{ia} (z) ] 
\end{align}
For the first term, we have
\begin{align}
 \frac{1}{N^4} \sum_{ia} \ee[  ( \del_{w} G_{ia} (w) ) G_{ia} (z)^3 ] = N^{\eps} \O ( N^{-2} \Psi(z)^3 \Psi (w) v^{-1} + N^{-3} v^{-1} )
\end{align}
where we used a crude bound $| \del_w G_{ii} (w) | \leq C v^{-1}$ for the diagonal entries, which holds with overwhelming probability. For the term on the last line of \eqref{eqn:cov-4o2} we have by the same argument, 
\beq \label{eqn:cov-4o2a1}
\frac{1}{N^4} \sum_{ia} \ee[  ( \del_{w} G_{ia} (w) ) G_{ii} (z) G_{aa} (z) G_{ia} (z) ] = N^{\eps} \O ( N^{-2} \Psi(z) \Psi(w) v^{-1} + N^{-3} v^{-1} ) .
\eeq
For the next term contributing to the LHS of \eqref{eqn:cov-4o-result} we have, using \eqref{eqn:delGia} and \eqref{eqn:del2Gjj},
\begin{align} \label{eqn:cov-4o3}
 & \frac{3}{N^4} \sum_{ija} (1+\delta_{ia})^3 \ee[ (\del_{ia} G_{ia} (z) ) ( \del_{ia}^2 G_{jj} (w) ) ] \notag\\
= & - \frac{12}{N^4} \sum_{ia} \ee[ ( \del_w G_{ii} (w) ) G_{aa} (w) G_{ii} (z) G_{aa} (z) ] \notag\\
- & \frac{12}{N^4} \sum_{ia} \ee[ ( \del_w G_{ii} (w) ) G_{aa} (w) G_{ia}(z)^2 ] \notag\\
- & \frac{12}{N^4} \sum_{ia} \ee[ ( \del_w G_{ia} (w) ) G_{ia} (w) G_{ii} (z) G_{aa} (z)   ] \notag\\
- &\frac{12}{N^4} \sum_{ia} \ee[ ( \del_w G_{ia} (w) ) G_{ia} (w) G_{ia}(z)^2  ] 
\end{align}
We begin with the first term on the right side of \eqref{eqn:cov-4o3}. By expanding each resolvent entry around $\msc(z)$ or $\msc(w)$ and using the local law \eqref{eqn:ll} to estimate the terms that are linear in the $G_{ii}$ and the entrywise law \eqref{eqn:entry-ll} to estimate the terms quadratic or higher in the $G_{ii}$ we obtain,
\begin{align}
 &- \frac{12}{N^4} \sum_{ia} \ee[ ( \del_w G_{ii} (w) ) G_{aa} (w) G_{ii} (z) G_{aa} (z) ]  \notag \\
 = & - \frac{12 \msc'(w) \msc (w) \msc(z)^2}{N^2} + N^{\eps} \O (N^{-3} ( y^{-1} v^{-1/2} + v^{-2} ) + N^{-2}(v^{-1} \Psi(w) \Psi (z)  )).
\end{align}
For the second term on the right side of \eqref{eqn:cov-4o3} we have,
\beq
\frac{1}{N^4} \sum_{ia} \ee[ ( \del_w G_{ii} (w) ) G_{aa} (w) G_{ia}(z)^2 ]  = N^{\eps} \O ( N^{-3} y^{-1} v^{-1} ).
\eeq
The third term is,
\begin{align}
& \frac{1}{N^4} \sum_{ia} \ee[ ( \del_w G_{ia} (w) ) G_{ia} (w) G_{ii} (z) G_{aa} (z)   ] \notag\\
= & \frac{1}{N^3} \msc(z)^2 \ee[ \del_w^2 m_N (w) ] + \frac{2 \msc (z) }{N^4} \sum_i \ee[ ( \del_w^2 G_{ii} (w) ) (G_{ii} (z) - \msc (z) ) ]  \notag\\
+ & N^{\eps} \O ( N^{-2} \Psi(z)^2 \Psi(w)^2 v^{-1} + N^{-3} v^{-1}) .
\end{align}
Clearly,
\beq
\frac{1}{N^3} \msc(z)^2 \ee[ \del_w^2 m_N (w) ] = N^{\eps} \O (N^{-3} v^{-2} ),
\eeq
and 
\begin{align}
 \frac{1 }{N^4} \sum_i \ee[ ( \del_w^2 G_{ii} (w) ) (G_{ii} (z) - \msc (z) ) ] = N^{\eps} \O ( N^{-3} v^{-2} ).
\end{align}
The last term on the right side of \eqref{eqn:cov-4o3} is,
\begin{align}
\frac{1}{N^4} \sum_{ia} \ee[ ( \del_w G_{ia} (w) ) G_{ia} (w) G_{ia}(z)^2  ] = N^{\eps} \O ( N^{-2} \Psi(z)^2 \Psi(w)^2 v^{-1} + N^{-3} v^{-1} )
\end{align}
The final contribution to the left side of \eqref{eqn:cov-4o-result} is, using \eqref{eqn:del3Gjj},
\begin{align}
 -& \frac{1}{N^4} \sum_{ija} \ee[ G_{ia} (z) (\del_{ia}^3 G_{jj} (w) ) ] \notag \\
= & \frac{24}{N^4} \sum_{ia} \ee[ G_{ia} (z) G_{ai} (w) G_{ii} (w) \del_w G_{aa} (w) ]  \notag \\
+ & \frac{12}{N^4} \sum_{ia} \ee[  ( \del_{w} G_{ia} (w) ) G_{ia} (z) ( G_{aa} (w) G_{ii} (w) + G_{ia} (w)^2 ) ] .
\end{align}
By splitting into off-diagonal and diagonal terms, we see that the first term is
\beq
\frac{24}{N^4} \sum_{ia} \ee[ G_{ia} (z) G_{ai} (w) G_{ii} (w) \del_w G_{aa} (w) ]  = N^{\eps} \O ( N^{-2} \Psi(z) \Psi(w) v^{-1}  + N^{-3} v^{-1} )
\eeq
and the second term is
\beq
\frac{1}{N^4} \sum_{ia} \ee[  ( \del_{w} G_{ia} (w) ) G_{ia} (z) ( G_{aa} (w) G_{ii} (w) + G_{ia} (w)^2 ) ]  =  N^{\eps} \O ( N^{-2} \Psi(z) \Psi(w) v^{-1}  + N^{-3} v^{-1} ).
\eeq
This completes the proof. \qed

\subsubsection{Fifth order}

\bel \label{lem:cov-5o}
We have,
\begin{align}
 & \frac{1}{N^{9/2}} \sum_{i j a} (1+\delta_{ia} )^4 \ee[ \del_{ia}^4 (G_{ai} (z) (G_{jj} (w) - \ee[ G_{jj} (w) ] ) ] \notag \\
= & N^{\eps} \O (N^{-5/2} \Psi(z) v^{-1} + N^{-5/2} \Psi(w) v^{-1/2}).
\end{align}
\eel
\proof We write the term as,
\begin{align}  \label{eqn:cov-5o}
&\frac{1}{N^{9/2} } \sum_{i j a} (1+\delta_{ia} )^4 \ee[ \del_{ia}^4 (G_{ai} (z) (G_{jj} (w) - \ee[ G_{jj} (w) ] ) ]  \notag \\
= & \frac{1}{N^{7/2} } \sum_{ia} (1+\delta_{ia} )^4 \ee[ ( \del^4_{ia} G_{ai} (z) ( m_N (w) - \ee[ m_N (w) ] ) ] \notag \\
- & \frac{1}{N^{9/2}} \sum_{n=0}^3 c_n \sum_{ia} ( 1 + \delta_{ia} )^4 \ee[ ( \del^{3-n}_{ia} G_{ai} (z) ) ( \del_w \del_{ia}^n G_{ia} (w) ) ]
\end{align}
for some combinatorial constants $c_n$. Consider the first term on the right side of \eqref{eqn:cov-5o}. The expression $(1+\delta_{ia})^4 \del_{ia}^4 G_{ia}$ is a linear combination of products of five resolvent entries.  Among five resolvent entries there are 10 indices, 5 of which must be $i$ and 5 of which must be $a$. Therefore, each term contains at least one off-diagonal entry $G_{ia}$. Therefore, with overwhelming probability we have,
\beq \label{eqn:cov-5o-parity}
| \del_{ia}^4 G_{ia} | \leq C\delta_{ia} + N^{\eps} \Psi(z).
\eeq
It follows that the second line of \eqref{eqn:cov-5o} is $\O (N^{-5/2} \Psi(z) v^{-1})$. By the same parity argument, we find that for any derivative of even order $n = 2m$ there exists $C$ such that
\beq
|\del_{ia}^{2m} G_{ia} |\leq C\delta_{ia} + N^{\eps} \Psi(z),
\eeq
with overwhelming probability. Therefore, for any $n=0, 1, 2, 3$ we have
\beq
| \del^{3-n}_{ia} G_{ai} (z) ) ( \del_w \del_{ia}^n G_{ia} (w) | \leq N^{\eps} v^{-1} ( \delta_{ia} + \Psi(z) + \Psi (w) ),
\eeq
with overwhelming probability. It follows that the last line of \eqref{eqn:cov-5o} is $N^{\eps} \O (N^{-5/2} \Psi(z) v^{-1} + N^{-5/2} \Psi(w) v^{-1/2})$. This yields the claim. \qed

\subsubsection{Sixth and seventh order} 

\bel
We have,
\begin{align}
 & \frac{1}{N^{5}} \sum_{i j a} (1+\delta_{ia} )^5 \ee[ \del_{ia}^5 (G_{ai} (z) (G_{jj} (w) - \ee[ G_{jj} (w) ] ) ] \notag  \\
= & N^{\eps} \O ( N^{-3} v^{-1} ).
\end{align}
and
\begin{align}
 & \frac{1}{N^{11/2}} \sum_{i j a} (1+\delta_{ia} )^6 \ee[ \del_{ia}^6 (G_{ai} (z) (G_{jj} (w) - \ee[ G_{jj} (w) ] ) ] \notag  \\
= & N^{\eps} \O ( N^{-3} v^{-1} ).
\end{align}
\eel
\proof We consider, 
\begin{align}
&\frac{1}{N^{5} } \sum_{i j a} (1+\delta_{ia} )^5 \ee[ \del_{ia}^5 (G_{ai} (z) (G_{jj} (w) - \ee[ G_{jj} (w) ] ) ]  \notag \\
= & \frac{1}{N^{4} } \sum_{ia} (1+\delta_{ia} )^5 \ee[ ( \del^5_{ia} G_{ai} (z) ( m_N (w) - \ee[ m_N (w) ] ) ] \notag \\
- & \frac{1}{N^{5}} \sum_{n=0}^4 c_n \sum_{ia} ( 1 + \delta_{ia} )^5 \ee[ ( \del^{4-n}_{ia} G_{ai} (z) ) ( \del_w \del_{ia}^n G_{ia} (w) ) ]
\end{align}
for some combinatorial constants $c_n$. 
We easily conclude the proof from, the fact that for any $n$ we have
\beq
| \del_{ia}^k G_{ia} | \leq C, \qquad | \del_{ia}^k \del_w G_{ia} (w) | \leq C v^{-1},
\eeq
with overwhelming probability. The result for the seventh order term is proven in the same way. \qed

\subsubsection{Proof of Proposition \ref{prop:cov-self}}

This follows from the expansion \eqref{eqn:cov-cum-exp}; in each of the previous few subsections, the lemmas provide the required estimates for  each of the terms in the expansion. \qed

\subsection{Covariance estimates} \label{sec:cov-ests}

Recall the notation,
\beq
\kappa (E) = | |E| - 2|.
\eeq
For $z = x + \i y$ we define $\kappa(z) = \kappa(x)$. 
Since
\beq
|z + 2 \msc (z) | \asymp \sqrt{ \kappa + y },
\eeq
we see from Proposition \ref{prop:cov-self} and Lemma \ref{lem:a-cov-2} that,
\begin{align} \label{eqn:cov-complex}
& N^2 \ee[ (m_N(z) - \ee[ m_N(z) ] )( m_N(w) - \ee[ m_N(w) ] ) ] \notag\\
= &2 \frac{ \msc'(z) \msc'(w)}{( 1- \msc (z) \msc (w) )^2} - 4 s_3 N^{-1/2} \msc'(z) \msc'(w) ( \msc(z) + \msc(w) ) \notag  \\
+ & 2 s_4 \msc'(w) \msc (w) \msc'(z) \msc(z) \notag\\
+ & \frac{ N^{\eps}}{ \sqrt{ \kappa(x) + y }} \O ( N^{-3} ( y^{-2} v^{-1} + (y+v)^{-1} v^{-2} )) \notag \\
+ & \frac{ N^{\eps}}{\sqrt{ \kappa(x) + y }} \O ( N^{-3/2} v^{-1} ( \Psi (z)^2 \Psi (w) + \Psi(z) \Psi(w)^2 ) +N^{-2} \Psi(z) \Psi (w) v^{-1} )
\end{align}
From this, we conclude the following proof.

\noindent{\bf Proof of Theorem \ref{thm:global-cov}}.  We apply \eqref{eqn:cov-complex} assuming that $\sqrt{ \kappa(x) + y } \geq \sqrt{ \kappa(u) + v }$. Otherwise, we reverse the roles of $z$ and $w$. Since $v, y \geq N^{-1/2}$ we have,
\beq
\Psi(z) \leq C  \frac{ ( \kappa(x) + y)^{1/4}}{(N y)^{1/2} },
\eeq
and a similar estimate for $\Psi (w)$. 
This yields the claim. \qed

Let $\varphi$ be a function supported in $[-5, 5]$ and recall the following quasi-analytic extension of $\varphi$,
\beq \label{eqn:quasi-analytic-def}
\tilphi (x + \i y ) = ( \varphi (x) + \i y \varphi' (x) ) \chi (y),
\eeq
where $\chi$ is smooth even bump function that is on $[-1, 1]$ and $0$ outside of $[-2, 2]$. Recall the formula
\beq
\del_{\bar{z}} \tilphi (x+ \i y ) = \frac{1}{2} \i y \chi(y) \varphi''(x) + \i ( \varphi (x) + \i \varphi'(x) y ) \chi'(y),
\eeq
as well as the Helffer-Sj\"ostrand (HS) formula \cite[Appendix C]{benaych2019lectures}
\beq
\tr \varphi (H) - \ee[ \varphi (H) ] = \frac{1}{ \pi} \int_{\rr^2} \del_{\bar{z}} \tilphi (x + \i y )N (m_N (z) - \ee[ m_N (z) ] ) \d x \d y.
\eeq
Recall the definition of the norm $\| \varphi \|_{1, w}$ in \eqref{eqn:weighted-norm}. 
The following is a simple consequence of the HS formula.
\bel
Assume that $\| \varphi''\|_1 \leq N$ and let $\eps >0$. With overwhelming probability we have,
\beq
\left| \tr \varphi (H) - \ee[ \tr \varphi (H) ] \right| \leq N^{\eps} ( 1 + \| \varphi'\|_1 + \| \varphi''\|_1 N^{-1} )
\eeq
\eel
\proof We have by a direct application of the HS formula and the local law \eqref{eqn:ll},
\begin{align}
\left| \tr \varphi (H) - \ee[ \tr \varphi (H) ] \right| &\leq C \left|  \int_{\rr^2} y \varphi'' (x) \chi (y) N( m_N ( x + \i y ) - \ee[ m_N ( x + \i y ) ] ) \d x \d y \right|  \notag\\
&+ C N^{\eps} ( \| \varphi \|_1 + \| \varphi' \|_1 ).
\end{align}
For the term on the first line, by the symmetry of $\chi (y)$ we have,
\begin{align}
&\int_{\rr^2} y \varphi'' (x) \chi (y) N( m_N ( x + \i y ) - \ee[ m_N ( x + \i y ) ] ) \d x \d y \notag\\
= & 2 \int_{y>0} y \varphi'' (x) \chi (y) N \Im[ m_N ( x + \i y ) - \ee[ m_N ( x + \i y ) ]  ]\d x \d y \notag
\end{align}
We have that $ y \to y \Im [ m_N ( x + \i y )]$ is increasing and so using the local law \eqref{eqn:ll} it follows that for $0 < y < N^{\eps-1}$ that,
\beq
yN \Im [ m_N ( x + \i y ) ] \leq N^{\eps}
\eeq
with overwhelming probability and so, 
\beq
\left| \int_{y < N^{\eps-1}} y \varphi'' (x) \chi (y) N \Im[ m_N ( x + \i y ) - \ee[ m_N ( x + \i y ) ]  ]\d x \d y \right|  \leq C N^{2 \eps-1} \| \varphi'' \|_1
\eeq
 Let $ t  = N^{\eps} ( 1  + \| \varphi' \|_1 ) / ( \| \varphi'' \| )$. Then,
 \begin{align}
& \left| \int_{ N^{\eps-1}<y<t} y \varphi'' (x) \chi (y) N \Im[ m_N ( x + \i y ) - \ee[ m_N ( x + \i y ) ]  ]\d x \d y \right|  \notag\\
\leq & \int_{ y < t } C N^{\eps} \| \varphi''\|_1 \d y \leq C N^{\eps} ( 1 + \| \varphi' \|_1 ).
\end{align}
For larger $y$, we have by integration by parts and the compact support of $\varphi$ that
\begin{align}
& \left| \int_{ y > t} y \varphi'' (x) \chi (y) N \Im[ m_N ( x + \i y ) - \ee[ m_N ( x + \i y ) ]  ]\d x \d y \right| \notag\\
= & \left| \int_{ y > t} y \varphi' (x) \chi (y) N \Im[ m'_N ( x + \i y ) - \ee[ m'_N ( x + \i y ) ]  ]\d x \d y \right| \notag\\
\leq & N^{\eps} \int_{3 >|y|>t} \| \varphi' \|_1 y^{-1} \d y \leq  CN^{2\eps} \| \varphi'\|_1 
\end{align}
where we used the Cauchy integral formula. \qed

Following \cite{meso}, define
\beq \label{eqn:Oma-def}
\Oma := \{ (x, y) \in \rr^2 : |y| > N^{ \mfa -1 } \}
\eeq
where $0 < \mfa < 1 $. 

\bel \label{lem:HS-cutoff}
Let $\eps >0, \mfa >0$ and $H$ be a Wigner matrix. With overwhelming probability,
\begin{align}
\left| ( \tr \varphi (H) - \ee[ \tr \varphi (H) ] ) -  \frac{1}{ \pi }\int_{\Oma} \del_{\bar{z}} \tilphi (z) N ( m_N (z) - \ee[ m_N(z) ] ) \d x \d y \right| \leq C N^{\eps} N^{\mfa-1} \| \varphi'' \|_1  .
\end{align} 
\eel
\proof This is a straightforward conseqeunce of the local law. The argument of the previous lemma shows that the contribution from $|y| < N^{\eps-1}$ is $\O ( N^{2\eps-1} \| \varphi'' \|_1)$. The contribution from $N^{\eps-1} < |y| < N^{\mfa-1}$ is $\O ( N^{\eps+\mfa-1} \| \varphi'' \|_1 )$ due to the local law \eqref{eqn:ll} and the explicit form of the quasi-analytic extension $\varphi$.  \qed

\bet \label{thm:var-global-2} Assume $\| \varphi''\|_{1, w} \leq N$. 
We have for any $\eps >0$,
\begin{align}
\left| \Var ( \tr \varphi (H) ) - V ( \varphi ) \right| \leq  N^{\eps-2} \| \varphi'' \|_{1, w}^2 + N^{\eps-1}(1+  \| \varphi' \|_{1, w}  ) \| \varphi'' \|_{1, w}
\end{align}
where,
\begin{align} \label{eqn:varf-def}
V( \varphi ) &= \frac{1}{ 2 \pi^2} \int_{-2}^2 \int_{-2}^2 \frac{ (\varphi(x) - \varphi(u) )^2}{ (x-u)^2} \frac{ 4- xu}{ \sqrt{ 4 - x^2} \sqrt{ 4 - u^2} } \d x \d y \notag\\
&+ \frac{s_4}{2 \pi^2} \left( \int_{-2}^2 \varphi (x) \frac{ 2 -x^2}{ \sqrt{4-x^2} } \d x \right)^2 \notag\\
&- \frac{2s_3}{  \pi^2 N^{1/2}} \left( \int_{-2}^2 \varphi (x) \frac{ 2 -x^2}{ \sqrt{4-x^2} } \d x \right) \left( \int_{-2}^2 \varphi (x) \frac{x}{ \sqrt{4-x^2} } \d x \right)
\end{align}
\eet
\proof  By the previous two lemmas we have,
\beq
\Var( \tr \varphi (H) ) = \frac{1}{ \pi^2}\ee \left[ \left( \int_{\Oma} \del_{\bar{z}} \tilphi N (m_N (z) - \ee[ m_N (z) ] \d x \d y ) \right)^2 \right] + \O ( N^{2\mfa+ \eps -1} \| \varphi'' \|_1 (1+ \| \varphi' \|_1 ) )
\eeq
The first term on the right side is,
\begin{align}
 & \ee \left[ \left( \int_{\Oma} \del_{\bar{z}} \tilphi N (m_N (z) - \ee[ m_N (z) ]) \d x \d y ) \right)^2 \right] \notag \\
= & \int_{\Oma^2} \del_{\bar{z}} \tilphi (z) \del_{\bar{w}} \tilphi (w) N^2 \ee\big[ (m_N (z) - \ee[ m_N (z) ] ) ( m_N (w) - \ee[ m_N (w) ])\big] \d x \d y \d u \d v
\end{align}
where we use the notation $z = x + \i y$ and $w = u + \i v$. 
Let $F(z, w)$ be the function defined by the second and third lines of \eqref{eqn:cov-complex} so that,
\begin{align}
 |H_N (z, w) | :=& \left| N^2 \ee[ (m_N (z) - \ee[ m_N (z) ] ) ( m_N (w) - \ee[ m_N (w) ] ) ] - F(z, w) \right| \notag\\
\leq & \frac{N^{\eps}}{ \sqrt{ \kappa(x) + |y| } } \left( \frac{1}{N y^2 |v|} + \frac{1}{ N v^2 |y|} \right).
\end{align}
We simplified the errors in \eqref{eqn:cov-complex} using the upper bound of \eqref{eqn:Psi-simple}. 
We now use this estimate to replace the terms $m_N(z), m_N (w)$ with $F(z, w)$ in the above integral expression for the variance. 
Clearly,
\beq
\int_{\Oma^2} | \chi'(y) \chi'(v) (\varphi (x) + \i y \varphi' (x) )(\varphi (u ) + \i v \varphi'(u) )H_N (z, w) | \d x \d y \d u \d v \leq CN^{\eps-1} ( 1 + \| \varphi' \|_{1, w} )^2
\eeq
and
\begin{align}
& \int_{\Oma^2} | \chi (y) \chi' (v) y \varphi''(x) (\varphi (u ) + \i v \varphi'(u) )H_N (z, w) | \d x \d y \d u \d v \notag\\
\leq &  C N^{\eps-1} \| \varphi'' (x) \|_{1, w} (1 + \| \varphi' \|_{1, w} ) \int_{N^{-1} < |y| < 1 } |y|^{-1} \d y \leq CN^{2 \eps-1} \| \varphi'' \|_{1, w} (1 + \| \varphi' \|_{1, w} ).
\end{align}
Define $t = N^{\eps} (1 + \| \varphi' \|_{1, w} ) \| \varphi'' \|_{1, w}^{-1} \wedge 1$.  We have,
\begin{align}
\int_{ |y|,|v| < t} | \varphi''(x)  \varphi''(v) y v H_N(z, w) | \d x \d y \d u \d v \leq & \| \varphi''\|_{1, w}^2 N^{\eps-1} \int_{|y|, |v| < t } |y|^{-1} \d v \d y \notag\\
\leq & C N^{3 \eps-1} \| \varphi'' \|_{1, w} ( 1 + \| \varphi' \|_{1, w} ) .
\end{align}
By integration by parts in the $x$ variable, and the fact that $z \to H_N (z, w)$ is holomorphic away from the real axis and the Cauchy integral formula,
\begin{align}
& \left| \int_{|y| >t , |v| < t} \varphi''(x) y \chi'(y) \varphi''(u) v \chi' (v ) H_N(z, w) \d x \d y \d u \d v \right| \notag\\
= & \left| \int_{|y| >t , |v| < t} \varphi'(x) y \chi'(y) \varphi''(u) v \chi' (v ) \del_z H_N(z, w) \d x \d y \d u \d v \right| \notag\\
\leq & N^{\eps} \| \varphi' \|_{1, w} \| \varphi'' \|_{1, w} \int_{|y| > t , |v| < t } (y^{-2} + |yv|^{-1} ) \d y \d v \leq C N^{2 \eps} \| \varphi'\|_{1, w} \| \varphi'' \|_{1, w} .
\end{align}
Finally,
\begin{align}
& \left| \int_{|y| >t , |v| > t} \varphi''(x) y \chi'(y) \varphi''(u) v \chi' (v ) H_N(z, w) \d x \d y \d u \d v \right| \notag\\
= & \left| \int_{|y| >t , |v| > t} \varphi'(x) y \chi'(y) \varphi'(u) v \chi' (v ) \del_z \del_w H_N(z, w) \d x \d y \d u \d v \right| \notag\\
\leq & N^{\eps-1} \| \varphi'\|_{1, w}^2 \int_{|y|, |v| > t } y^{-2} |v|^{-1} \d v \d y \leq N^{2 \eps-1} \| \varphi' \|_{1, w} \| \varphi'' \|_{1, w}
\end{align}
The above estimates show that,
\begin{align}
& \int_{\Oma^2} \del_{\bar{z}} \tilphi (z) \del_{\bar{w}} \tilphi (w) N^2 \ee[ (m_N (z) - \ee[ m_N (z) ] ) ( m_N (w) - \ee[ m_N (w) ) ] \notag \\
= & \int_{\Oma^2} \del_{\bar{z}} \tilphi (z) \del_{\bar{w}} \tilphi (w)  F(z, w)  + N^{3\eps-1} \O ( \| \varphi'' \|_{1,w} (1+ \| \varphi' \|_{1, w} )) .
\end{align}
Using the estimate,
\beq \label{eqn:Fzwbd}
|F(z, w) | \leq C \frac{1}{ \sqrt{ \kappa(x) + y} \sqrt{ \kappa(u) + v } } \frac{1}{ y^2 + v^2}
\eeq
it is straightforward to conclude that,
\begin{align}
 \int_{\Oma^2} \del_{\bar{z}} \tilphi (z) \del_{\bar{w}} \tilphi (w)  F(z, w) &= \int_{\cc^2}  \del_{\bar{z}} \tilphi (z) \del_{\bar{w}} \tilphi (w)  F(z, w) \notag \\
  +&  N^{\eps} \O ( N^{2\mfa-2} \| \varphi''\|_{1, w}^2 + N^{\mfa-1} \| \varphi'\|_{1, w}\| \varphi''\|_{1, w} )
\end{align}
The integral on the right side may then be calculated by Green's theorem which in complex notation states that
\beq
\int_\Omega \del_{\bar{z}} F (z) \d x \d y = - \frac{\i}{2} \int_{\del \Omega} F(z) \d z.
\eeq
Therefore,
\begin{align}
 & \frac{1}{ \pi^2} \int_{\cc^2} \del_{\bar{z}} \tilphi (z) \del_{\bar{w}} \tilphi (w) F(z, w) \d x \d y \d u \d v \notag\\
= &- \frac{1}{ 4 \pi^2} \int_{\rr^2} \varphi (x) \varphi (y) ( F(x + \i 0 , y + \i 0 ) + F ( x - \i 0 , y  - \i 0) - 2 F ( x + \i 0, y - \i 0 ) ) \d x \d y.
\end{align}
The contributions from the terms that do not involve $s_3$ were calculated in the proof of \cite[Lemma 4.2]{meso}. The term that involves $s_3$ may be calculated as
\begin{align}
 & -  \frac{8 s_3}{N^{1/2} \pi^2 } \int_{\cc^2} \del_{\bar{z}} \tilphi (z) \del_{\bar{w} } \tilphi (w) \msc'(z) \msc'(w) \msc(z) \d x \d y \d u \d v \notag\\
= &- \frac{8 s_3}{ N^{1/2} \pi^2} \left( \int_{\rr} \varphi (x) \Im [ \msc' (x + \i 0) \msc ( x + \i 0) ] \d x \right) \left( \int_{\rr} \varphi (y) \Im [ \msc' (y + \i 0) ] \d y \right) \notag\\
= &  - \frac{8 s_3}{ N^{1/2} \pi^2} \left( \int_{-2}^2 \varphi (x) \frac{ x^2-2}{2 \sqrt{4- x^2 }} \d x \right) \left( \int_{-2}^2 \varphi (y) \frac{ - y}{2 \sqrt{ 4 - y^2} } \d y \right)
\end{align}
where the last line follows from the explicit formula of $\msc$. 
We conclude the proof. \qed

\subsection{Expectation}

\bel \label{lem:expect-self}
Let $H$ be a Wigner matrix. Let $z  = x + \i y $.  Then,
\begin{align}
(z  + 2 \msc (z) ) ( \ee[ m_N (z)] - \msc (z) ) &= - \frac{1}{N} \msc' (z) \notag\\
+& \frac{ 4 s_3}{N^{3/2} } \msc(z)^3 - \frac{s_4}{N} \msc(z)^4 \\
+& N^{\eps} \O ( (N y)^{-2} ).
\end{align}
\eel
\proof By the cumulant expansion, 
\begin{align} \label{eqn:expect-expand}
\frac{1}{N} \sum_i z \ee[ G_{ii} (z) ] +1 =& \frac{1}{N} \sum_{ia} \ee[ H_{ia} G_{ia} (z) ] \notag\\
= & - \ee[ m_N (z)^2] - \del_z \ee[ m_N (z) ] \notag\\ 
+ & \sum_{n=3}^5 \frac{1}{N^{1+n/2}} \sum_{ia} \frac{ (1+\delta_{ia})^{n-1} s_n}{(n-1)!} \ee[ \del_{ia}^n G_{ia} ] \notag \\
+ & N^{\eps} \O (N^{-2} ).
\end{align}
By the local law \eqref{eqn:ll} we have,
\begin{align}
- \ee[ m_N (z)^2] - \frac{1}{N}\del_z \ee[ m_N (z) ]  &= - \msc(z)^2 - 2 \msc(z) \ee[ m_N (z) - \msc(z) ] - \frac{1}{N} \del_z \msc(z) \notag\\
&+ N^{\eps} \O ((Ny)^{-2} ).
\end{align}
By \eqref{eqn:del2Gia} and similar arguments as in the proof of Lemma \ref{lem:cov-3o}, we have that
\beq
\frac{1}{N^{5/2}} \sum_{ia} ( 1+\delta_{ia} )^2\ee[ \del_{ia}^2 G_{ia} ] = \frac{8}{N^{3/2}} \msc(z)^3 + N^{\eps} \O ( (N y)^{-2} ).
\eeq
From \eqref{eqn:del3Gia} and the entrywise local law \eqref{eqn:entry-ll} we see,
\begin{align}
\frac{1}{N^3} \sum_{ia} (1+ \delta_{ia})^3 \ee[ \del_{ia}^3 G_{ia} ] = - \frac{6}{N^3} \sum_{ia} \ee[ G_{ii}^2 G_{aa}^2 ] + N^{\eps} \O ((N y)^{-2} ).
\end{align}
Expanding $G_{ii} = (G_{ii} - \msc ) + \msc$ (and similarly for $G_{aa}$) and estimating the terms linear in $G_{ii}$ using the local law \eqref{eqn:ll} and the higher order terms using the entry-wise local law \eqref{eqn:entry-ll} we see that,
\begin{align}
\frac{1}{N^3} \sum_{ia} \ee[ G_{ii}^2 G_{aa}^2 ]  = \frac{ \msc(z)^4}{N} + N^{\eps}\O ( (N y)^{-2} ).
\end{align}
Finally, the term with $n=5$ in \eqref{eqn:expect-expand} is estimate using the parity estimate \eqref{eqn:cov-5o-parity}. \qed

\bel \label{lem:wig-expect}
Let $f$ be a function, and $H$ a Wigner matrix. Then,
\begin{align}
\ee[ \tr f (H) ] &= N \int f(x) \rhosc (x) \d x - \frac{1}{ 2 \pi }\int f(x) \frac{1}{\sqrt{4-x^2} } \d x \notag \\
&+ \frac{ f(2) + f(-2) }{4} + \frac{s_4}{2 \pi} \int_{-2}^2 f(x) \frac{ x^4 - 4 x^2 + 2}{\sqrt{4-x^2}} \d x \notag\\
&+ \frac{2 s_3}{N^{1/2}} \frac{1}{ \pi} \int_{-2}^2 f(x) \frac{ 3x - x^3}{ \sqrt{4 -x^2}} + N^{\eps} \O ( N^{-1} \| f'' \|_{1, w}  )
\end{align}
\eel 
\proof The following is proven in a straightforward manner from the HS formula and Lemma \ref{lem:expect-self}, similar to the argument proving Theorem \ref{thm:var-global-2},
\begin{align}
\ee [ \tr f (H ) ] - N \int f(x ) \rhosc (x) \d x &= \frac{1}{ \pi} \int_{\Oma} \del_{\bar{z}} \tilf (z) \frac{ - \msc' (z) -s_4 \msc(z)^4 + 4s_3N^{-1/2} \msc(z)^3 }{ z + 2 \msc (z) } \d x \d y \notag\\ 
&+ N^{\eps} \O ( N^{\mfa-1} \| f'' \|_{1, w} )
\end{align}
where $\tilf$ is as in \eqref{eqn:quasi-analytic-def} and $\Oma$ is as in \eqref{eqn:Oma-def}. From the estimate,
\beq
\left|  \frac{ - \msc' (z) }{ z + 2 \msc (z) } \right| C \leq \frac{1}{ \kappa(x) + y },
\eeq 
we see that
\begin{align}
 &\int_{\Oma} \del_{\bar{z}} \tilf (z) \frac{ - \msc' (z) -s_4 \msc(z)^4 + 4s_3N^{-1/2} \msc(z)^3 }{ z + 2 \msc (z) } \d x \d y \notag\\
 = & \int_{\cc} \del_{\bar{z}} \tilf (z) \frac{  - \msc' (z) -s_4 \msc(z)^4 + 4s_3N^{-1/2} \msc(z)^3 }{ z + 2 \msc (z) } \d x \d y + N^{\eps} \O ( N^{\mfa-1} \| f'' \|_{1, w} ).
\end{align}
The integral without the $s_3 N^{-1/2}$ term as calculated in \cite[Lemma 4.3]{meso}. Calculating the additional term is similar to what appears there, and the result follows.
 \qed

\section{Mesoscopic estimates for Dyson Brownian Motion} \label{sec:char}

In this section, we will consider Dyson Brownian Motion (DBM),
\beq
\d \lambda_i (t) = \sqrt{ \frac{2}{N \beta} } \d B_i (t) + \frac{1}{N} \sum_{j \neq i } \frac{1}{ \lambda_i (t) - \lambda_j (t) } \d t - \frac{ \lambda_i (t) }{2} \d t
\eeq
where the initial data are the eigenvalues of a Wigner matrix. Compared to \eqref{eqn:int-dbm} we introduced the parameter $\beta  \in \rr$ above which equals $\beta=1$ in the real symmetric case and $\beta=2$ in the complex Hermitian case. The parameter does not play a substantial role but we introduce it for completeness.

 We will consider the evolution of
\beq
s_N(z, t) := \frac{1}{N} \sum_i \log ( \lambda_i (t) - z )
\eeq
for $z$ in the upper half-plane. By the It{\^{o}} formula,
\begin{align}
\d s_N = \frac{ - m_N^2}{2} + \frac{1}{ 2N } \left( 1- \frac{2}{ \beta} \right) m_N' - \frac{1}{2} - \frac{ z m_N (z)}{2} + \frac{1}{N^{3/2}}\sqrt{ \frac{2}{ \beta}} \sum_i \frac{ \d B_i }{ \lambda_i - z}.
\end{align}
Define now,
\begin{align} \label{eqn:char-hh1}
s(z) := \int \log (x -z ) \rhosc (x) \d x.
\end{align}
The characteristics as defined in \eqref{eqn:char-def} go backwards in time. We need the characteristics that go forward in time. For this, we fix a final time $T$ and define
\beq
\tilz_{t} := z_{T-t}
\eeq
for some final condition $\tilz_{T} = z_{0} = z$ which we suppress from the notation.  Note that due to \eqref{eqn:delz-bwd}, $\tilz_s$ satisfies
\beq \label{eqn:char-a2}
\del_s \tilz_s = - \msc ( \tilz_s) - \frac{ \tilz_s}{2}.
\eeq
Then, with $s (z)$ as in \eqref{eqn:char-hh1},
\beq
\del_s s ( \tilz_s ) = \msc(\tilz_s)^2 + \tilz_s \frac{ \msc(\tilz_s ) }{2} = \frac{ \msc ( \tilz_s)^2 -1}{2},
\eeq
where we used $\msc^2 + z \msc + 1 = 0$. Plugging in the characteristic, we see that
\begin{align} \label{eqn:char-a1}
\d s_N (\tilz_s) &= -m_N^2/2 + \frac{1}{2N} \left( 1- \frac{2}{\beta} \right) m_N'  - \frac{1}{2} - \tilz_s m_N /2 + \frac{1}{N^{3/2}}\sqrt{\frac{2}{ \beta} } \sum_i \frac{ \d B_i}{ \lambda_i - \tilz_s} \\
&+ m_N \msc + m_N \tilz_s /2 .
\end{align}
In Section \ref{sec:char-ests} we estimate the difference of $\d s_N ( \tilz_s) - \del_s s ( \tilz_s)$ and use this to find an estimate of $s_N (\tilz_T, T)$ in terms of $s_N (\tilz_0, 0)$ and an explicit Gaussian random variable. In Section \ref{sec:char-cov} we calculate the variance of this Gaussian random variable as well as the covariance between two of these variables for different characteristics that are well-separated. 

Before beginning, we record the following easy consequence of rigidity.
\bep \label{prop:dbm-ll}
Let $\eps >0$ and $\tau >0$. With overwhelming probability we have simultaneously for all $t \in [0, 1]$ and $z \in \D_{\tau,1}$ that,
\beq
| m_N (z) - \msc (z) | \leq \frac{ N^{\eps}}{ N \Im [z]}
\eeq
as well as
\beq
| \lambda_i (t) - \gamma_i | \leq \frac{ N^{\eps}}{N^{2/3} \min \{ i^{1/3}, (N+1-i)^{1/3}  \} }.
\eeq
\eep

\subsection{Estimates along characteristics} \label{sec:char-ests}

\bep \label{prop:char-est}
Let $\tilz_s$, $s \in [0, T]$ be a characteristic, and assume that the final condition $z = \tilz_T$ satisfies $\Re[z] \in (-2+\kappa, 2 - \kappa)$ and $\Im [z] \geq N^{\tau-1}$ for some $\tau, \kappa >0$. Assume $T < (\log(N))^{-1}$.   Let $\eps >0$. With overwhelming probability,
\begin{align}
s_N (\tilz_T, T) - s (\tilz_T) & = s_N ( \tilz_0, 0) - s ( \tilz_0) \notag\\
&+ \frac{1}{N^{3/2}} \sqrt{ \frac{2}{\beta}} \sum_i \int_0^T \frac{ \d B_i (s) }{ \gamma_i - \tilz_s } \notag\\
&+ \frac{1}{2N} \left( 1 - \frac{2}{\beta} \right) \int_0^T \msc' ( \tilz_s) \d s \notag\\
&+ N^{\eps} \O \left( \frac{1}{ N^2 \Im [z_T] } \right)
\end{align}
\eep
\proof Using \eqref{eqn:char-a2} and \eqref{eqn:char-a1} we have,
\begin{align}
\d ( s_N ( \tilz_s ) - s( \tilz_s ) ) &= - \frac{1}{2} ( \msc ( \tilz_s) - m_N ( \tilz_s) )^2 \label{eqn:char-err-1} \\
&+ \frac{1}{2N} \left( 1- \frac{2}{ \beta} \right) (m_N'( \tilz_s) - \msc' ( \tilz_s) ) \label{eqn:char-err-2} \\
&+ \frac{1}{2N} \left( 1- \frac{2}{ \beta} \right)  \msc' ( \tilz_s) \label{eqn:char-err-3} \\
&\frac{1}{N^{3/2}}\sqrt{\frac{2}{ \beta} } \sum_i \d B_i \left( \frac{1}{ \lambda_i -\tilz_s} - \frac{1}{ \gamma_i - \tilz_s } \right)\label{eqn:char-err-4}   \\
&+ \frac{1}{N^{3/2}}\sqrt{\frac{2}{ \beta} } \sum_i \frac{ \d B_i}{ \gamma_i - \tilz_s} \label{eqn:char-mart}.
\end{align}
We now estimate the integrals in $s$ of the various terms on the right.  By the local law estimates of Proposition \ref{prop:dbm-ll} we have, with overwhelming probability,
\begin{align}
\int_0^T \left| m_N ( \tilz_s) - \msc ( \tilz_s ) \right|^2 \d s \leq N^{\eps} \int_0^T \frac{1}{N^2 ( \Im [\tilz_s ] )^2 }\d s.
\end{align}
From the explicit form of the characteristics we see that,
\beq
c \leq - \del_s \Im[ \tilz_s ] \leq c^{-1}
\eeq
for some $c>0$ depending on $\kappa >0$ and our assumption that $T  < ( \log(N) )^{-1}$. Therefore, for the term \eqref{eqn:char-err-1} we have the estimate,
\beq
\int_0^T \left| m_N ( \tilz_s) - \msc ( \tilz_s ) \right|^2 \d s \leq N^{\eps} \frac{1}{N^2 \Im [\tilz_T] }.
\eeq
By the Cauchy integral formula we have for all $ s \in [0, T]$,
\beq
\left| m_N' ( \tilz_s ) - \msc' ( \tilz_s ) \right| \leq \frac{N^{\eps}}{N  \Im [ \tilz_s]^2}
\eeq
and so we arrive at the same estimate for the time integral of \eqref{eqn:char-err-2} as we did for \eqref{eqn:char-err-1}.  For the term \eqref{eqn:char-err-4}, we have by the rigidity estimates of Proposition \ref{prop:dbm-ll} that,
\beq
 \sum_i \left| \frac{1}{ \lambda_i - \tilz_s } - \frac{1}{ \gamma_i - \tilz_s} \right|^2 \leq \frac{N^{\eps}}{N ( \Im [\tilz_s ] )^3}
\eeq
with overwhelming probability for all $ s\in [0, T]$. Therefore, by the Burkholder-Davis-Gundy inequality:
\begin{align}
\left| \int_0^T \frac{1}{N^{3/2}}\sqrt{\frac{2}{ \beta} } \sum_i \d B_i \left( \frac{1}{ \lambda_i -\tilz_s} - \frac{1}{ \gamma_i - \tilz_s } \right) \right| \leq \frac{N^{\eps}}{N^2 \Im [ \tilz_T] }
\end{align}
with overwhelming probability. The claim follows. \qed

\subsection{Covariance calculation} \label{sec:char-cov}

We now calculate the variance of the Gaussian random variables appearing in Proposition \ref{prop:char-est}. 
\bel
Let $\tilz_s$ be a characteristic whose final condition $\tilz_T$ satisfies $\Re[ \tilz_T] \in (-2+\kappa, 2- \kappa)$ and $\Im [ \tilz_T ] \geq N^{\tau-1}$. Assume that $T < \log(N)^{-1}$. Then,
\begin{align}
 2\frac{1}{N} \sum_i \int_0^T \frac{ \Im [ \tilz_s ]^2}{ | \gamma_i - \tilz_s |^4} \d s 
 =&  \log(  \Im [ \tilz_0 ] / \Im[ \tilz_T]  )  \notag\\
 -&\frac{T}{2} + \Re[ \log ( 1- \msc ( \tilz_T)^2 ) - \log ( 1- \msc ( \tilz_0)^2 ) ] \notag\\
 + & \O( (N \Im [ \tilz_T] )^{-1} ).
\end{align}
\eel
\remark Note that as long as $z$ is separated from $\pm 2$ by a constant $c$, then $|1- \msc^2 (z) | \geq c'$ for some $c' >0$. 

\proof First, it is easy to see that,
\beq
\left| \frac{1}{N} \sum_{i} \frac{1}{ | \gamma_i - \tilz_s |^4} - \int \frac{1}{ |x - \tilz_s |^4} \rhosc (x) \d x \right| \leq \frac{C}{N \Im[ \tilz_s]^4}.
\eeq
Therefore,
\beq
\left| \frac{1}{N} \sum_i \int_0^T \frac{ \Im [ \tilz_s ]^2}{ | \gamma_i - \tilz_s |^4}  \d s  - \int_0^T \int \frac{ \Im [ \tilz_s ]^2}{ | x - \tilz_s |^4} \rhosc (x) \d x \d s\right| \leq \frac{C}{ N \Im [\tilz_T] }.
\eeq
The integral against $\rhosc$ admits an exact calculation. First, observe that
\beq
 2 \int \frac{ \Im [ \tilz_s ]^2}{ | x - \tilz_s |^4} \rhosc (x) \d x = - \Re[ \msc' ( \tilz_s ) ] + \frac{ \Im [ \msc ( \tilz_s ) ]}{ \Im[\tilz_s]}. 
\eeq
Since,
\beq
\del_s \tilz_s = - \msc ( \tilz_s) - \frac{ \tilz_s}{2}
\eeq
we have
\beq
\frac{ \Im [ \msc ( \tilz_s) ] }{ \Im [ \tilz_s] } = \frac{ - \del_s \Im [ \tilz_s]}{ \Im [z_s] } - \frac{1}{2}.
\eeq
Next,
\beq
\del_s \log (1 - \msc^2 ( \tilz_s ) ) = \frac{  - 2 \msc \msc' ( \tilz_s ) }{ 1 - \msc^2 ( \tilz_s ) } \del_s \tilz_s = - \msc' ( \tilz_s)
\eeq
where we used,
\beq
\del_s \tilz_s = - \msc ( \tilz_s) - \tilz_s/2 = \frac{ 1 - \msc^2 ( \tilz_s) }{ 2 \msc ( \tilz_s ) }
\eeq
the second identity following from $\msc^2 + z \msc +1 = 0$. Therefore,
\beq
2 \int \frac{ \Im [ \tilz_s ]^2}{ | x - \tilz_s |^4} \rhosc (x) \d x =  \del_s \Re[ \log ( 1- \msc^2 ( \tilz_s ) ) ] -  \frac{ \del_s \Im [ \tilz_s]}{\Im [ \tilz_s] } - \frac{1}{2}.
\eeq
The claim follows after integration in $s$. \qed

\bel
Let $\tilz_s$ and $\tilw_s$ be two characteristics such that the final conditions satisfy $\Re[ \tilz_T], \Re [ \tilw_T] \in (-2+\kappa, 2 - \kappa)$ and $\Im [ \tilz_T], \Im[\tilw_T] \geq N^{\tau -1}$ for some $\tau, \kappa >0$. Assume also that,
\beq
\kappa \Im [ \tilz_T] \leq \Im [ \tilw_T] \leq \frac{1}{ \kappa} \Im [ \tilz_T].
\eeq
Assume that
\beq
T < \frac{ | \Re[ \tilz_T] - \Re[ \tilw_T] |}{ \log(N)}.
\eeq
Then,
\begin{align}
\left| \frac{1}{N}  \int_0^T \sum_i \frac{  \Im [ \tilz_s] \Im [ \tilw_s] }{ | \gamma_i - \tilz_s|^2 | \gamma_i - \tilw_s |^2} \d s \right| \leq C \frac{T^2 + ( \Im [ \tilz_T] )^2}{ ( \Re[ \tilz_T] - \Re[ \tilw_T])^2} + C \frac{1}{N \Im[ \tilz_T] }
\end{align}
\eel
\proof Due to the explicit form of the characteristics, we see that
\beq
c \Im [ \tilz_s] \leq \Im [ \tilw_s] \leq \frac{1}{c} \Im [ \tilz_s].
\eeq
for some $c>0$ for all $ s \in [0, T]$. 
Similar to the proof of the previous lemma we have,
\beq
\left| \frac{1}{N} \int_0^T \sum_i \frac{  \Im [ \tilz_s] \Im [ \tilw_s] }{ | \gamma_i - \tilz_s|^2 | \gamma_i - \tilw_s |^2} \d s  -\int_0^T \int \frac{ \Im [ \tilz_s ] \Im [ \tilw_s]}{ | x- \tilz_s|^2 | x- \tilw_s|^2} \rhosc (x) \d x \d s \right| \leq \frac{C}{N \Im [\tilz_T] }.
\eeq
The integral in $x$ is bounded via convolution of two Poisson kernels,
\beq
\int \frac{ \Im [ \tilz_s ] \Im [ \tilw_s]}{ | x- \tilz_s|^2 | x- \tilw_s|^2} \rhosc (x) \d x  \leq C \frac{ \Im [ \tilz_s]}{ ( \Re[ \tilz_s] - \Re[ \tilw_s] )^2 + ( \Im [ \tilz_s ] )^2 }.
\eeq
Under our assumptions we have,
\beq
| \Re[ \tilz_s] - \Re[ \tilw_s]  | \geq c | \Re[ \tilz_T] - \Re[ \tilw_T] |
\eeq
for all $ s \in [0, T]$, using the explicit form of the characteristics. Hence, using the explicit form of the characteristics to estimate $\Im [ \tilz_s ] \leq C ( \Im [ \tilz_T ] + T )$,
\begin{align}
\int_0^T \int \frac{ \Im [ \tilz_s ] \Im [ \tilw_s]}{ | x- \tilz_s|^2 | x- \tilw_s|^2} \rhosc (x) \d x  \d s \leq &  C  \int_0^T \frac{  \Im [ \tilz_T] + T }{ | \Re[ \tilz_T ] - \Re[ \tilw_T ] |^2} \d s \notag\\
\leq &  C \frac{T^2 + ( \Im [ \tilz_T])^2}{ ( \Re[ \tilz_T] - \Re[ \tilw_T])^2}
\end{align}
This yields the claim. \qed

The above two lemmas immediately imply the following.
\bep \label{prop:char-cov}
Let $\tilz_s$ and $\tilw_s$ be two characteristics such that the final conditions satisfy $\Re[ \tilz_T], \Re [ \tilw_T] \in (-2+\kappa, 2 - \kappa)$ and $\Im [ \tilz_T], \Im[\tilw_T] \geq N^{\tau -1}$ for some $\tau, \kappa >0$.  Assume that there is a $t = N^{-\omega}$ such that
\beq
t \asymp \Im [ \tilz_T] \asymp \Im [ \tilw_T],
\eeq
and a parameter $\delta >0$ such that,
\beq
t \leq N^{-\delta} | \Re[ \tilz_T ] - \Re[ \tilw_T] |.
\eeq
Assume that there is $0 < \delta_1 < \delta$ such that
\beq
T \asymp N^{\delta_1} t.
\eeq
Then there is a $c>0$ so that the variances of the Gaussian random variables
\beq
Z_z :=\frac{1}{N^{1/2}} \sum_i \int_0^T \frac{ \d B_i}{ \gamma_i - \tilz_s} \d s, \qquad Z_w := \frac{1}{N^{1/2}} \sum_i \int_0^T \frac{ \d B_i}{ \gamma_i - \tilw_s} \d s
\eeq
satisfy
\beq
\Var ( Z_z) \geq c \delta_1 \log(N), \qquad \Var (Z_w) \geq c \delta_1 \log(N)
\eeq
and their covariance satisfies
\beq
| \Cov (Z_z, Z_w ) | \leq C N^{2 ( \delta_1 -\delta)} + (Nt)^{-1} \leq \frac{C}{ \log(N)}.
\eeq
\eep
\proof The estimate for the covariance follows immediately. For the variance, note that $\Im [ \tilz_0 ] \asymp \Im [ \tilz_T ] +T \asymp T$ by the explicit form of the characteristics and the assumption that $T \gg t$. The lower bound follows. \qed

\section{Estimates for the Gaussian ensembles} \label{sec:gaussians}

If $\rho_k$ are the $k$-point functions of a point process (in our case the eigenvalue of a random matrix), i.e., 
\beq
\sum_i \ee[  f( \lambda_i ) ] = \int f (x) N \rho_1 (x) \d x, \qquad \sum_{ i \neq j } \ee[ f (\lambda_i, \lambda_j ) ] = N (N-1) \int f (x, y) \rho_2 (x, y) \d x \d y
\eeq
then for two functions $f, g$ we have
\beq
\Cov \left( \tr f (H) , \tr g (H) \right) = \frac{1}{2} \int \left( f(x) - f(y) \right) \left( g (x) - g (y) \right) T_N (x, y) \d x \d y
\eeq
where the cluster function is,
\beq
T_N (x, y) = N^2 \rho_1 (x) \rho_1 (y) - N (N-1) \rho_2 (x, y).
\eeq
Introduce the function,
\beq \label{eqn:geta-def}
g_{E, \eta} (x) := \frac{ \eta}{ ( x - E)^2 + \eta^2}
\eeq
The un-normalized Hermite polynomials are
\beq \label{eqn:unnorm-poly}
H_n (x) = (-1)^n \e^{ x^2} \left( \frac{ \d }{ \d x } \right)^n \e^{ - x^2}.
\eeq
They satisfy,
\beq
\int_\rr \e^{ - x^2} H_n (x) H_m(x) = \delta_{nm} 2^n n! \sqrt{ \pi}.
\eeq
The normalized Hermite polyomials $h_n$ are,
\beq \label{eqn:norm-poly}
h_n(x) := \frac{1}{ \sqrt{ 2^n n! \sqrt{ \pi } } } H_n (x)
\eeq
which satisfy,
\beq
\int_{\rr} h_i (x) h_j (x) \e^{-x^2} \d x = \delta_{ij}.
\eeq

\subsection{Covariance for GUE}

With the above notation, define,
\beq \label{eqn:hermite-unscaled}
\psi_n (x) := \e^{-x^2/4} h_n ( x / \sqrt{2} ) 2^{-1/4}.
\eeq
For the GUE, the Christoffel-Darboux formula is (see, e.g., Chapter 5 of \cite{pastur2011eigenvalue}),
\beq
T_N (x, y) = N\left( \frac{ \psi_N ( \sqrt{N} x ) \psi_{N-1} ( \sqrt{N} y ) - \psi_N ( \sqrt{N} y ) \psi_{N-1} ( \sqrt{N} x ) }{x-y} \right)^2 =: \frac{\Psi_N(x, y)^2}{(x-y)^2}.
\eeq
The estimates in Appendix \ref{a:hermite} imply that,
\beq \label{eqn:psiN-bd}
| \Psi_N (x, y) |^2 \leq C \left( \frac{\e^{- c N (|x| - 2)_+^{3/2} }}{|x-2|^{1/2} + N^{-1/3} } \right) \left( \frac{\e^{- c N (|y| - 2)_+^{3/2} }}{|y-2|^{1/2} + N^{-1/3} } \right)
\eeq
\bep \label{prop:gue-cov-calc}
Let $0 < \eta, \delta < 10^{-1}$. There is a $C>0$ depending on $\delta$ but not on $\eta >0$ so that the following holds. For any $E_1, E_2 $ with $|E_1|, |E_2| < 2 - 2 \delta$ we have,
\beq
\left| \Cov ( \tr g_{E_1, \eta } (H) , \tr g_{E_2, \eta } (H) ) \right| \leq C \frac{1}{ (E_1 - E_2)^2 + \eta^2 }
\eeq
\eep
\proof We see from the exponential tails in \eqref{eqn:psiN-bd} that,
\begin{align}
\Cov ( \tr g_{E_1, \eta } (H) , \tr g_{E_2, \eta } (H) ) &= \frac{1}{2}\int_{|x|, |y| < 2 + \delta} \frac{ (g_1 (x) - g_1 (y) )( g_2 (x) - g_2 (y) ) }{(x-y)^2} \Psi_N(x, y)^2 \d x \d y  \notag\\
&+ \O \left( \eta^2 \e^{- c N } \right),
\end{align}
where we abbreviate $g_i  = g_{E_i, \eta}$.  If $||x|-2| < \delta$ and $|y| < 2 + \delta$ it is not hard to see that,
\beq
\left| \frac{ (g_1 (x) - g_1 (y) )( g_2 (x) - g_2 (y) ) }{(x-y)^2}  \right| \leq \frac{C \eta^2}{ ( (y -E_1)^2 + \eta^2) ( (y-E_2 )^2 + \eta^2 )}
\eeq
Therefore,
\begin{align}
|\Cov ( \tr g_{1 } (H) , \tr g_{2 } (H) ) | \leq & \left|  \frac{1}{2}\int_{|x|, |y| < 2 - \delta} \frac{ (g_1 (x) - g_1 (y) )( g_2 (x) - g_2 (y) ) }{(x-y)^2} \Psi_N(x, y)^2 \d x \d y  \right| \notag\\
+& C \frac{ \eta}{ ((E_1 - E_2)^2 + \eta^2 }
\end{align}
by the convolution property of Poisson kernels. 
Now,
\begin{align}
& \frac{ (g_1 (x) - g_1 (y) )( g_2 (x) - g_2 (y) ) }{(x-y)^2}   \notag\\
=&  \eta^2 \frac{ (x - E_1 + y-E_1 )(x-E_2 + y-E_2)}{ ( (x-E_1)^2 + \eta^2 )( (y-E_1)^2 + \eta^2 ) ( (x-E_2)^2 + \eta^2 )( (y-E_2)^2 + \eta^2 )}.
\end{align}
The claim then follows from Lemma \ref{lem:cal-bd-1} below and \eqref{eqn:psiN-bd} above. \qed

\bel \label{lem:cal-bd-1}
For $0 \leq E < 5$ and $0 < \eta < 10^{-1}$ we have,
\beq
\eta^2 \int_{|x|, |y| < 10} \frac{ |(x+y)(x+y- 2 E)|}{ (x^2 + \eta^2)(y^2 + \eta^2)( (x-E)^2+\eta^2)( (y-E)^2+\eta^2 ) } \d x \d y  \leq C \frac{1}{ E^2 + \eta^2 }.
\eeq
\eel
\proof We write the numerator as $(x+y)(x+y-2 E) = x(x-E)  + y(y-E) + x (y-E) + y (x-E)$.  By symmetry it suffices to bound the $x(x-E)$ term and the $x(y-E)$ term.  For the  first we have,
\begin{align}
& \int_{|x|, |y| < 10} \frac{| x(x-E) |}{ (x^2 + \eta^2)(y^2 + \eta^2)( (x-E)^2+\eta^2)( (y-E)^2+\eta^2 ) } \d x \d y \notag\\
 =&  \left( \int_{|x|<10} \frac{ | x (x-E)|}{ (x^2 + \eta^2)( (x-E)^2 + \eta^2 ) } \d x \right) 
\left( \int_{|y|<10} \frac{ 1}{ (y^2 + \eta^2)( (y-E)^2 + \eta^2 )} \d y \right).
\end{align}
The $y$ integration is bounded above by $\eta^{-1} ( E^2 + \eta^2 )^{-1}$ by the convolution property of the Poisson kernel. The integrand in the $x$ integral can be bounded by $( x^2 + \eta^2)^{-1} + ((x-E)^2  + \eta^2 )^{-1}$ by Cauchy-Schwarz, and so yields $\O (\eta^{-1})$.

We now estimate the integral corresponding to $|x(y-E)|$. We have,
\begin{align}
\int_{|x| < 10 } \frac{ |x|}{(x^2+\eta^2 )( (x-E)^2 + \eta^2 )}  \d x &= \int_{|x| < E/2} \frac{ |x|}{(x^2+\eta^2 )( (x-E)^2 + \eta^2 )}  \d x \notag\\
&+ \int_{2 E >|x| > E/2} \frac{ |x|}{(x^2+\eta^2 )( (x-E)^2 + \eta^2 )}  \d x \notag\\
&+ \int_{10>|x| > 2 E} \frac{ |x|}{(x^2+\eta^2 )( (x-E)^2 + \eta^2 )}  \d x .
\end{align}
The last integral is bounded by,
\beq
\int_{|x| > 2 E} \frac{ |x|}{(x^2+\eta^2 )( (x-E)^2 + \eta^2 )}  \d x \leq C \int_\rr \frac{1}{ x^3 + (E + \eta)^3} \d x \leq \frac{C}{ E^2 + \eta^2}.
\eeq
The second last integral is bounded by,
\beq
\int_{2 E >|x| > E/2} \frac{ |x|}{(x^2+\eta^2 )( (x-E)^2 + \eta^2 )}  \d x \leq \frac{C E}{ E^2 + \eta^2} \int_\rr \frac{1}{ (x-E)^2  + \eta^2 } \d x \leq C \frac{E}{ \eta ( E^2 + \eta^2 ) }.
\eeq
The first integral is bounded by,
\beq
\int_{|x| < E/2} \frac{ |x|}{(x^2+\eta^2 )( (x-E)^2 + \eta^2 )}  \d x  \leq C \frac{ |E|}{ E^2 + \eta^2 } \int_\rr \frac{1}{ x^2 + \eta^2} \d x \leq C \frac{E}{ \eta ( E^2 + \eta^2 ) }.
\eeq
Therefore,
\beq
\int_{|x| < 10 } \frac{ |x|}{(x^2+\eta^2 )( (x-E)^2 + \eta^2 )}  \d x \leq C \frac{1}{ \eta (E+\eta)}.
\eeq
By making the change of variable $x = y-E$ we see by the same argument as above,
\beq
\int_{|y| < 10 } \frac{ |y-E|}{(y^2+\eta^2 )( (y-E)^2 + \eta^2 )}  \d y \leq C \frac{1}{ \eta (E+\eta)}.
\eeq
This completes the proof. \qed

\subsection{GOE preliminaries}

We first make some simplifications that will allow us to avoid considering the behavior of the cluster function near the spectral edges.  Write $R_N (x, y)$ for twice the cluster function of the GOE, i.e.,
\beq \label{eqn:GOE-cluster-def}
\frac{1}{2} R_N(x, y) = N^2 p_1 (x) p_1 (y) - N (N-1 ) p_2 (x, y)
\eeq
Let now $g_i$ be as in the GUE case above, i.e., a Poisson kernel centered at $|E_i| \leq 2 - 6 \delta$, some $\delta >0$. Let $\chi$ be a bump function that is $1$ on $[-2 +4 \delta, 2 - 4\delta]$ and $0$ outside of $[-2+3\delta, 2 -3 \delta]$.  Define $f_i = \chi g_i$.

Let $\chi_1$ be a bump function that is $0$ outside of $[-2+\delta, 2 - \delta]$ and $1$ on $[-2+2 \delta, 2 - 2 \delta]$.
\bel \label{lem:goe-edge-simplification}
Let $\eps >0$ and $D>0$. Then, for $E_i$, $f_i$ as $0 < \eta < 1$ as above we have
\begin{align}
\left|\Cov ( \tr g_1 (H) , \tr g_2 (H) ) \right|  & \leq \left| \frac{1}{4} \int_{\rr^2} (f_1(x) - f_1 (y) ) ( f_2 (x) - f_2 (y) ) \chi_1 (x) \chi_1 (y) R_N (x, y) \d x \d y \right| \notag \\
& + C N^{\eps} \frac{1}{ (E_1 - E_2 )^2 + \eta^2 } + \eta^{-2} N^{-D}. \label{eqn:goe-edge-simplification}
\end{align}
\eel
\proof 
We have, (writing $\tr h (H) = h$ for notational simplicity),
\begin{align}
\Cov ( g_1, g_2 ) =& \Cov ( \chi g_1, \chi g_2 ) \notag\\
 +& \Cov ( (1- \chi ) g_1, g_2 ) + \Cov ( g_1, (1- \chi ) g_2 ) 
- \Cov ( (1- \chi ) g_1, (1- \chi ) g_2 ).
\end{align} 
Since $(1- \chi) g_i$ equals $\eta$ times a smooth function having derivatives bounded uniformly in $N$ it is straightforward to conclude by Theorem \ref{thm:var-global-2} or rigidity that,
\beq
\Var ( (1- \chi) g_i ) \leq C \eta^2 N^{\eps}.
\eeq
Moreover, since the function $\eta \to  \eta \Im [ m_N ( x + \i \eta)]$ is increasing, it follows that
\beq
\Var (g_i ) \leq \eta^{-2} \ee[ N^{\eps} \Im [ m_N ( E_i + \i N^{\eps-1} ) ]^2 ] \leq C N^{3 \eps} \eta^{-2},
\eeq
where the last inequality follows by the local law \eqref{eqn:ll}. 
It follows by Cauchy-Schwarz that,
\beq
| \Cov (g_1, g_2 ) | \leq | \Cov (f_1, f_2 ) | + C N^{\eps}
\eeq
where we defined $f_i = \chi g_i$.  

   We now show that we can restrict the domain of integration in the formula for the covariance of the $f_i$'s. Let $\chi_1$ be a bump function as defined before the statement of the lemma that is $0$ outside of $[-2+\delta, 2 - \delta]$ and $1$ on $[-2+2 \delta, 2 - 2 \delta]$.  Recall that $R_N$ is symmetric and so,
\begin{align}
\Cov (\tr f_1 (H), \tr f_2 (H) ) &= \frac{1}{4} \int_{\rr^2} ( f_1 (x) - f_1 (y) ) ( f_2 (x) - f_2 (y) ) R_N (x, y) \d x \d y \notag\\
&= \frac{1}{2} \int_{\rr^2} f_1(x) f_2(x) R_N (x, y) \d x \d y - \frac{1}{2} \int_{\rr^2} f_1 (x) f_2 (y) R_N (x, y) \d x \d y.
\end{align}
On the other hand, consider,
\begin{align}
& \frac{1}{4} \int_{\rr^2} (f_1 (x) - f_1 (y) ) ( f_2 (x) - f_2 (y) ) \chi_1 (x) \chi_1 (y) R_N (x, y) \d x \d y \notag\\
= & \frac{1}{2} \int_{\rr^2} f_1 (x) f_2 (x)  \chi_1 (x) \chi_1 (y) R_N (x, y) \d x \d y \notag\\
- & \frac{1}{2} \int_{\rr^2} f_1 (x) f_2 (y) \chi_1 (x) \chi_1 (y) R_N(x, y) \d x \d y.
\end{align}
Since the $f_i$ are supported in $[-2+3 \delta, 2-3 \delta]$ and $\chi_1$ is $1$ on this interval, we see that $f_1 (x) f_2 (y) \chi_1 (x) \chi_1 (y) = f_1 (x) f_2 (y)$ and $f_1 (x) f_2 (x)  \chi_1 (x) \chi_1 (y)= f_1(x) f_2 (x) \chi_1 (y)$. Therefore,
\begin{align}
& \Cov (\tr f_1 (H), \tr f_2 (H) )  -\frac{1}{4} \int_{\rr^2} (f_1 (x) - f_1 (y) ) ( f_2 (x) - f_2 (y) ) \chi_1 (x) \chi_1 (y) R_N (x, y) \d x \d y  \notag\\
= & \frac{1}{2} \int_{\rr^2} f_1 (x) f_2 (x) (1 - \chi_1 (y) ) R_N (x, y) \d x \d y.
\end{align}
Recall that $R_N (x, y) = 2 (N^2 p_1 (x) p_1 (y) - N(N-1) p_2 (x, y) )$ and so, letting $F(x) = f_1 (x) f_2 (x)$ and $G(y) = 1- \chi_1 (y)$ we have,
\begin{align}
- &\frac{1}{2} \int_{\rr^2} f_1 (x) f_2 (x) (1 - \chi_1 (y) ) R_N (x, y) \d x \d y \notag\\
= & \sum_{i \neq j } \ee[ F ( \lambda_i ) G( \lambda_j ) ] - \sum_{i, j} \ee[ F ( \lambda_i ) ] \ee[ G( \lambda_j ) ] \notag\\
= &\sum_{i, j } \ee[ F ( \lambda_i ) G( \lambda_j ) ] - \sum_{i, j} \ee[ F ( \lambda_i ) ] \ee[ G( \lambda_j ) ] \notag\\
= & \Cov \left( \tr F(H), \tr G(H) \right)
\end{align}
where in the second equality we used that since $F$ and $G$ have disjoint support, $F( \lambda_i ) G( \lambda_i ) = 0$.  To bound the last term we treat the cases of $\eta \lesssim N^{-1}$ and $\eta \gtrsim N^{-1}$ separately.

Let $\eps >0$, and consider $\eta \leq N^{\eps-1}$. By rigidity we have with overwhelming probability that $\tr G ( H ) - \ee[ \tr G (H) ] = \O ( N^{\eps} )$. Therefore, (using that $F$ is positive),
\beq
| \Cov \left( \tr F(H), \tr G(H) \right) | = | \ee[ \tr F(H) ( \tr G(H) - \ee[ \tr G(H) ] ) ] | \leq N^{\eps} \ee[ \tr F(H) ] + C N^{-D} \eta^{-2}
\eeq
for any $D>0$.  On the other hand, since the support of $F$ is restricted to the bulk, we have that, (since  the 1-point function or density of states of the GOE is bounded in the bulk, see Theorem \ref{thm:gaussian-bulk})
\beq
\ee[ \tr F(H) ] \leq C N \int_\rr f_1 (x) f_2 (x) \d x\leq C N \frac{\eta}{ (E_1 - E_2)^2 + \eta^2 }\leq \frac{ CN^{\eps}}{(E_1 - E_2)^2 + \eta^2 } 
\eeq
where we used the convolution property for Poisson kernels as well as the fact that we are assuming $N\eta \leq N^{\eps}$.

Assume now that $\eta \geq N^{\eps-1}$. Then a straightforward calculation using rigidity shows that with overwhelming probability,
\beq
| \tr F(H) - \ee[ \tr F(H) ] | \leq C N^{\eps} \| F' \|_1 \leq C \frac{ N^{\eps}}{ \eta} \| F \|_1 \leq CN^{\eps} \frac{1}{ (E_1 - E_2)^2 + \eta^2 },
\eeq
again using the convolution property of the Poisson kernels. 
This yields the claim. \qed

\subsection{Covariance for GOE, $N$ even}

In this section, we prove the following.
\bep \label{prop:goe-even}
Let $g_{E, \eta}$ be as in \eqref{eqn:geta-def}. Let $\kappa >0$ and $\mfc >0$ and assume that $|E_1|, |E_2| \leq 2 - \kappa$. Let $\eps>0$.  For any $N^{-\mfc}< \eta < \frac{1}{2}$ we have,
\beq
\left| \Cov \left(  \tr g_{E_1, \eta} (H) , \tr g_{E_2, \eta } (H) \right)  \right| \leq  C N^{\eps} \left( \frac{|E_1-E_2| }{ \eta ( (E_1 - E_2 )^2 + \eta^2 ) } + \frac{1}{ (E_1 - E_2 )^2 + \eta^2 } \right)
\eeq
where $H$ is a GOE matrix of even dimension.
\eep
In order to prove the above Proposition, we have by Lemma \ref{lem:goe-edge-simplification} that it suffices to estimate the first term on the right side of \eqref{eqn:goe-edge-simplification}. We now state the exact formula for the cluster function found in \cite{pastur2011eigenvalue}. 
With $\psi_n$ as in \eqref{eqn:hermite-unscaled}, introduce
\beq
\phiN (x) := N^{1/4} \psi_N ( \sqrt{N} x), \qquad \phiNm (x) := N^{1/4} \psi_{N-1} ( \sqrt{N } x )
\eeq
By the asymptotics in Appendix \ref{a:hermite}, we have
\beq \label{eqn:phiN-asymp}
| \phiN (x) | + |\phiNm (x) | \leq C \frac{ \e^{ - c N (|x| - 2 )_+^{3/2} }}{ ||x|-2 |^{1/4} + N^{-1/6} }
\eeq
Define the operator, for $ x\in \rr$ and $f : \rr \to \cc$,
\beq
\sg (x) := \frac{1}{2} \mathrm{sign} (x), \qquad  ( \sg f ) (t) := \int_{\rr} \sg (t -s ) f (s) \d s.
\eeq
Note that for sufficiently well-behaved $f$ we have,
\beq
\frac{\d}{\d t} \sg f = f .
\eeq
From (6.1.19) and Theorem 6.1.6 and (6.2.10) of \cite{pastur2011eigenvalue} we have for the cluster function defined as in \eqref{eqn:GOE-cluster-def} and $N$ even, 
\begin{align}
R_N (x, y) =& 2 S_N (x, y) S_N(y, x) \notag\\
 +& D_N(x, y) ( I_N (y, x) - \sg (y-x) ) + D_N (y, x) ( I_N (x, y) - \sg (x-y) )
\end{align}
and
\begin{align}
S_N (x, y) = K_N (x, y) + \frac{N}{2} \phiNm (x) \sg \phiN (y)
\end{align}
where
\begin{align}
K_N (x, y) = \frac{ \phiN (x) \phiNm (y) - \phiN (y) \phiNm (x) }{(x-y) }
\end{align}
(see \cite[(4.2.17), (5.1.3)]{pastur2011eigenvalue}} for the Christoffel-Darboux formula in this scaling)
and,
\begin{align}
D_N (x, y) = - \partial_y S_N (x, y), \qquad I_N (x, y) = \int_y^x S_N (s, y) \d s.
\end{align}
We apply the exact formula above, and decompose the covariance of $f_i$ into several terms corresponding to the various terms in $R_N$ above.  To facilitate this we further introduce,
\beq
E_{N,1} (x, y ) = \frac{N}{2} \phiNm (x) \sg \phiN (y), \qquad \tilde{I}_N (x, y ) = \int_{y}^x K_N (s, y) \d s
\eeq
and
\beq
F_{N, 1} (x, y) = \int_y^x E_{N, 1} (s, y) \d s.
\eeq
Importantly, note that
\beq
\del_x \tilde{I}_N (x, y) = K_N (x, y) , \qquad \del_x F_{N, 1} (x, y ) = E_{N, 1} (x, y).
\eeq
\noindent{\bf Proof of Proposition \ref{prop:goe-even}}. By Lemma \ref{lem:goe-edge-simplification} it suffices to estimate the first term on the right side of \eqref{eqn:goe-edge-simplification}. 
After changing switching the roles of $x$ and $y$ in a few terms we find the following formula for the first term on the RHS of \eqref{eqn:goe-edge-simplification}.
\begin{align}
& \frac{1}{4} \int_{\rr^2} (f_1 (x) - f_1 (y) ) (f_2 (x) - f_2 (y) ) R_N (x, y) \chi_1 (x) \chi_1 (y) \d x \d y \notag\\
= &  \frac{1}{2} \int (f_1 (x) - f_1 (y) )(f_2 (x) - f_2 (y) ) K_N (x, y)^2 \chi_1 (x) \chi_1 (y)  \d x \d y \label{eqn:even-cov-1} \\
- &  \frac{1}{2} \int (f_1 (x) - f_1 (y) )(f_2 (x) - f_2 (y) ) (  \partial_y K_N (x, y) ) \tilde{I}_N (y, x ) \chi_1 (x) \chi_1 (y)   \d x \d y \label{eqn:even-cov-2} \\
+&  \int (f_1 (x) - f_1 (y) )(f_2 (x) - f_2 (y) ) E_{N, 1} (x, y) K_N (y, x)  \chi_1 (x) \chi_1 (y)   \d x \d y \label{eqn:even-cov-3} \\
+ & \frac{1}{2} \int(f_1 (x) - f_1 (y) )(f_2 (x) - f_2 (y) ) E_{N, 1} (x, y) E_{N, 1} (y, x) \chi_1 (x) \chi_1 (y)  \d x \d y \label{eqn:even-cov-4} \\
- & \frac{1}{2} \int(f_1 (x) - f_1 (y) )(f_2 (x) - f_2 (y) ) ( \partial_y E_{N, 1} (x, y) ) \tilde{I}_N (y, x)   \chi_1 (x) \chi_1 (y) \d x \d y \label{eqn:even-cov-5} \\
- & \frac{1}{2} \int(f_1 (x) - f_1 (y) )(f_2 (x) - f_2 (y) ) ( \partial_y  K_N (x, y) ) F_{N, 1} (y, x)  \chi_1 (x) \chi_1 (y)  \d x \d y \label{eqn:even-cov-6} \\
- & \frac{1}{2} \int(f_1 (x) - f_1 (y) )(f_2 (x) - f_2 (y) ) ( \partial_y E_{N, 1} (x, y) ) F_{N, 1} (y, x) \chi_1 (x) \chi_1 (y)   \d x \d y \label{eqn:even-cov-7} \\
+ & \frac{1}{2} \int(f_1 (x) - f_1 (y) )(f_2 (x) - f_2 (y) )  ( \partial_y S_N (x, y) )  \sg (y-x ) \chi_1 (x) \chi_1 (y)  \d x \d y \label{eqn:even-cov-8}
\end{align}
Over the next sequence of lemmas, we bound each of the terms appearing above, which will conclude the Proposition. \qed 
\bel
\label{lem:even-1} For the term \eqref{eqn:even-cov-1} we have,
\begin{align}
\left|  \int (f_1 (x) - f_1 (y) )(f_2 (x) - f_2 (y) ) K_N (x, y)^2 \chi_1 (x) \chi_1 (y)  \d x \d y  \right| \leq \frac{C}{ (E_1 - E_2)^2 + \eta^2 }
\end{align}
\eel
\proof 
The region of integration is restricted to $|x|, |y| \leq 2 - \delta$ on which we have,
\beq \label{eqn:cov-KNest}
| K_N (x, y) | \leq \frac{C}{ |x-y|},
\eeq
which follows from the explicit form of $K_N$ above and the Hermite polynomial asymptotics in Appendix \ref{a:hermite}. 
When $|x|, |y| < 2 - 4 \delta$ we have that $f_i = g_i$ and so by the proof of Proposition \ref{prop:gue-cov-calc} (using Lemma \ref{lem:cal-bd-1}) we easily see that,
\beq
\int_{|x|, |y| < 2 - 4 \delta}\left|  \frac{ (f_1 (x) - f_1 (y) )(f_2 (x) - f_2 (y) ) }{(x-y)^2} \right| \d x \d y \leq \frac{C}{ (E_1 - E_2)^2 + \eta^2 }.
\eeq
When $|x|, |y| > 2 - 5 \delta$ we see from the explicit form of $f_i$ that,
\beq
\left|  \frac{ (f_1 (x) - f_1 (y) )(f_2 (x) - f_2 (y) ) }{(x-y)^2} \right| \leq C \eta^2 ,
\eeq
and so this portion of the integral contributes $\O (1)$.  This leaves the cross-terms where $|x|  > 2 -4  \delta$ and $|y| < 2 - 5 \delta$, or vice versa. By symmetry, it suffices to deal with the former region. When $|x| > 2 - 4 \delta$ and $|y| < 2 - 5 \delta$ we have that,
\beq
\left|  \frac{ (f_1 (x) - f_1 (y) )(f_2 (x) - f_2 (y) ) }{(x-y)^2} \right| \leq \frac{C \eta^2}{ (y-E_1)^2 + \eta^2}  \frac{1}{ (y-E_2)^2 + \eta^2}
\eeq
and so this region contributes $\O ( \eta/  ( (E_1-E_2)^2 + \eta^2 ) $. The claim follows. \qed


\bel \label{lem:cov-Itildeest}
We have for $|x|, |y| < 2 - \delta$ that
\beq \label{eqn:cov-Itildeest}
| \tilde{I}_N (x, y) | \leq  C \log(N) ( 1 \wedge ( N |x-y| ) ) .
\eeq
\eel
\proof We may assume $ x> y$; the proof for $x < y$ is similar.    We write, for some $y < A < x$,
\beq
\tilde{I}_N (x, y) = \int_y^{A} K_N (s, y) \d s + \int_A^{x} K_N(s, y) \d s.
\eeq
By \eqref{eqn:cov-KNest} we have,
\beq
\left| \int_A^{x} K_N(s, y) \d s \right| \leq C \log \left( \frac{x-y}{A-y} \right).
\eeq
For the other term we have by the explicit form of $K_N$ and the fact that the $\phiN$ and $\phiNm$ are bounded in the domain of integration,
\beq
\left|  \int_y^{A} K_N (s, y) \d s \right| \leq C \int_y^A \frac{| \phiN (s) - \phiN (y) | + | \phiNm (s) - \phiNm (y) |}{ s -y }  \d s.
\eeq
By \cite[(5.2.4)]{pastur2011eigenvalue}, we have that $| \partial_s \phiN (s) | + | \partial_s \phiNm (s) | \leq CN$ for $|s| < 2 - \delta$ and so this integral is bounded by $CN|A-y|$.  If $|x-y| \leq N^{-1}$ we take $A = x$. If $|x-y| \geq N^{-1}$ we take $A = y + N^{-1}$. The estimate follows. \qed


\bel \label{lem:even-2}
For the term \eqref{eqn:even-cov-2} we have,
\begin{align}
 & \left| \int (f_1 (x) - f_1 (y) )(f_2 (x) - f_2 (y) ) ( \partial_y K_N (x, y) ) \tilde{I}_N (y, x )  \chi_1 (x) \chi_1 (y)  \d x \d y  \right| \notag\\
\leq & C \log(N)  \left( \frac{|E_1 - E_2 | }{ \eta  ( ( E_1 - E_2)^2 + \eta^2 )} + \frac{1}{ ( E_1 - E_2)^2 + \eta^2 }  \right)
\end{align}
\eel
\proof 
We will treat the integral by integrating by parts in $y$. Note that the integrand is of compact support so there is no boundary term. For the term arising from the partial integration where the derivative hits $\tilde{I}_N$ we have, $\partial_y \tilde{I}_N (y, x) = K_N (y, x)$, and so this term was dealt with in Lemma \ref{lem:even-1}. Consider the term where where the partial derivative hits $\chi_1 (y)$. Note that
\beq
(f_1(x) - f_1 (y) ) (f_2 (x) - f_2 (y) ) \chi_1' (y) = f_1(x) f_2(x) \chi_1' (y)
\eeq
and that this term is non-zero only if $|x-y| \geq \delta$. Therefore,
\begin{align}
& \int_{\rr^2} |(f_1(x) - f_1 (y) ) (f_2 (x) - f_2 (y) ) \chi_1' (y) \chi_1 (x) K_N (x, y) \tilde{I}_N (x, y) | \d x \d y \notag\\
\leq & C \log(N) \int_{\rr} |f_1 (x) f_2 (x) \chi_1 (x) | \d x \leq C \log(N) \frac{ \eta}{ (E_1 - E_2)^2 + \eta^2 }.
\end{align}
It remains to consider the term,
\begin{align}
 & \int_{\rr^2} | ( \del_y [(f_1 (x) - f_1 (y) )(f_2 (x) - f_2 (y) ) ] ) K_N (x, y) \tilde{I}_N (y, x) |\chi_1 (x) \chi_1 (y) \d x \d y \notag\\
\leq &  C \log(N) \int_{\rr^2} \left| \frac{ \del_y [(f_1 (x) - f_1 (y) )(f_2 (x) - f_2 (y) ) ] }{x-y} \right| \chi_1 (x) \chi_1 (y) \d x \d y 
\end{align}
where we used \eqref{eqn:cov-KNest} and \eqref{eqn:cov-Itildeest}. 
%
If both $|x|, |y| > 2 - 5 \delta$ then,
\beq
\left| \frac{ \del_y [(f_1 (x) - f_1 (y) )(f_2 (x) - f_2 (y) ) ] }{x-y} \right| \leq C \eta^2
\eeq
by the explicit form of the $f_i$. If $|x| > 2 - 4 \delta$ and $|y| < 2 - 5 \delta$ then one can estimate easily,
\begin{align}
\left| \frac{ \del_y [(f_1 (x) - f_1 (y) )(f_2 (x) - f_2 (y) ) ] }{x-y} \right| \leq C \frac{ \eta}{ ( ( y - E_1 )^2 + \eta^2)( (y-E_2 )^2 + \eta^2 ) }.
\end{align}
The integral over $x$ and $y$ contributes an acceptable error by the convolution property of Poisson kernels.  For  the region $ |x| < 2 - 5 \delta$ and $|y| > 2 - 4 \delta$ we have,
\beq
\left| \frac{ \del_y [(f_1 (x) - f_1 (y) )(f_2 (x) - f_2 (y) ) ] }{x-y} \right| \leq C \frac{ \eta^2}{ ( (x-E_1)^2 + \eta^2 )( (x - E_2 )^2 + \eta^2 ) }
\eeq
and we obtain a similar (in fact, better) estimate. For $|x|, |y| < 2 - 4 \delta$ we have that $f_i = g_i$. By symmetry it will suffice to consider the case that the derivative hits $g_1(y)$.  In this case, we can estimate,
\begin{align}
\left| (\del_y g_1 (y) ) \frac{ g_2 (x) - g_2 (y) }{ (x-y)} \right| \leq C \eta^2 \frac{|y - E_1|}{ (y- E_1 )^4 + \eta^4 } \frac{ |x-E_2| + |y-E_2 | }{ (  ( x- E_2)^2 + \eta^2 ) (  (  y -E_2 )^2 + \eta^2 ) }.
\end{align}
By estimating the integral of the right side over $|x|, |y| < 10$ it suffices to check the case $E_1 = 0$ and $E_2 = E$ for some $0 < E < 5$.  There are two terms, corresponding to the two terms $|x-E_2|$ and $|y-E_2|$ on the right side above. The $|x-E_2|$ term is bounded by,
\begin{align}
 & \eta^2 \int_{|y| < 10 } \frac{|y|}{ (y^4 + \eta^4 )((y-E)^2 + \eta^2 ) } \d y \int_{ |x| < 10 } \frac{ |x-E|}{ ( x - E)^2 + \eta^2 } \d x \notag\\
\leq & C \log(N) \int_{|y| < 10} \frac{ \eta}{( y^2 + \eta^2 ) ( (y-E)^2 + \eta^2 ) } \d y \notag\\
\leq & C \log(N) \frac{1}{ E^2 + \eta^2}.
\end{align}
In the first inequality we did the $x$ integration which contributes $\O ( | \log \eta | )  = \O ( \log(N))$ and cleared a power of $y^2 + \eta^2$ in the denominator using the $|y\eta|$ in the numerator. The last estimate follows from the convolution property of Poisson kernels. 

The final term to consider is,
\begin{align}
 & \eta^2 \int_{|y| < 10 } \frac{ |y| |y-E|}{(y^4 + \eta^4 ) ( (y-E)^2 + \eta^2 )} \d y \int_{ |x| < 10 } \frac{1}{ (x-E)^2 + \eta ^2 } \d x \notag\\
\leq & C \eta \int_{|y| < 10 } \frac{ |y| |y-E|}{(y^4 + \eta^4 ) ( (y-E)^2 + \eta^2 )} \d y  
\end{align}
For this integral, first,
\beq
 \eta \int_{|y| < E/2 } \frac{ |y| |y-E|}{(y^4 + \eta^4 ) ( (y-E)^2 + \eta^2 )} \d y \leq C \frac{E}{ E^2 + \eta^2 } \int \frac{ \eta |y|}{ y^4 + \eta^4 } \d y \leq C \frac{E}{ \eta ( E^2 + \eta^2 ) }.
\eeq
Second,
\beq
\eta \int_{|y| > E/2} \frac{ |y| |y-E|}{(y^4 + \eta^4 ) ( (y-E)^2 + \eta^2 )} \d  y \leq \frac{C}{ E^2 + \eta^2 } \int_{|y| < 10 } \frac{ |y-E|}{ (y-E)^2 + \eta^2 } \d y \leq \frac{C \log(N) }{ E^2 + \eta^2 }.
\eeq
Summarizing, we have shown
\begin{align}
 & \int_{|x|, |y | < 2 -4 \delta } \left| (\del_y g_1 (y) ) \frac{ g_2 (x) - g_2 (y) }{ (x-y)} \right| \d x \d y  \notag\\
\leq &  C \log(N) \left( \frac{|E_1 - E_2 | }{ \eta  ( ( E_1 - E_2)^2 + \eta^2 )} + \frac{1}{ ( E_1 - E_2)^2 + \eta^2 } \right).
\end{align}
As stated above, the same argument applies to when the derivative hits $g_2$ instead. This completes the proof. \qed

\bel \label{lem:even-cov-EN}
We have,  for $|x|, |y| < 2 - \delta$,
\beq
| E_{N, 1} (x, y) | \leq C, \qquad |F_{N, 1} (x, y) | \leq C .
\eeq
\eel
\proof The  estimate follows from \cite[(6.2.9)]{pastur2011eigenvalue}.  \qed 

\bel \label{lem:even-3}
We have for the terms \eqref{eqn:even-cov-3} and \eqref{eqn:even-cov-4},
\begin{align}
 & \left|  \int (f_1 (x) - f_1 (y) )(f_2 (x) - f_2 (y) )  E_{N, 1} (x, y) K_N (y, x) \chi_1 (x) \chi_1 (y)  \d x \d y \right| \notag\\
 \leq &  \frac{C}{ (E_1 - E_2)^2 + \eta^2 }
\end{align}
and
\begin{align}
\left|  \int(f_1 (x) - f_1 (y) )(f_2 (x) - f_2 (y) ) E_{N, 1} (x, y) E_{N, 1} (y, x) \chi_1 (x) \chi_1 (y) \d x \d y  \right|  \leq \frac{C}{ (E_1 - E_2)^2 + \eta^2 }
\end{align}
\eel
\proof By Lemma \ref{lem:even-cov-EN} we have that $|E_{N, 1} (x , y) | \leq C$ for all $x$ and $y$ in the region where the integrand is non-zero. Therefore, we can crudely estimate,
\beq
|E_{N, 1} (s, t) | \leq \frac{C}{|s-t|}.
\eeq
The integral of $|(f_1 (x) - f_1 (y) )(f_2 (x) - f_2 (y) )(x-y)^{-2}|$ over $|x|, |y| < 2-\delta$ was estimated in the proof of Lemma \ref{lem:even-1}. \qed

\bel
For the term \eqref{eqn:even-cov-5} we have,
\begin{align}
& \left|  \int(f_1 (x) - f_1 (y) )(f_2 (x) - f_2 (y) ) ( \partial_y E_{N, 1} (x, y) )\tilde{I}_N (y, x) \chi_1 (x) \chi_1 (y) \d x \d y  \right| \notag\\
\leq & C \log(N) \left( \frac{|E_1 - E_2 | }{ \eta  ( ( E_1 - E_2)^2 + \eta^2 )} + \frac{1}{ ( E_1 - E_2)^2 + \eta^2 } \right)
\end{align}
\eel
\proof 
We integrate by parts in $y$. There is no boundary term due to the compact support of the integrand. When this derivative hits $\tilde{I}_N (y, x)$ one gets $K_N (y, x)$. This term was dealt with in Lemma \ref{lem:even-3}. When the derivative hits the $f_i$'s, one can, as in Lemma \ref{lem:even-3} simply bound $|E_{N, 1} (s, t) | \leq C |s-t|^{-1}$ and so one is left with estimating the integral of
\beq
\left|  ( \del_y  [ (f_1 (x) - f_1 (y) )( f_2 (x) - f_2 (y)  ] ) |x-y|^{-1}  \tilde{I}_N (x, y) \right| 
\eeq
over $|x|, |y| < 2 - \delta$. This was completed in Lemma \ref{lem:even-2}.  The final term is when the derivative hits $\chi_1 (y)$, which is treated similarly to the analogous term in the proof of Lemma \ref{lem:even-2}. That is, due to the support properties of the $f_i$ we have that  the integrand is bounded by,
\beq
| (f_1 (x) - f_1 (y) ) ( f_2 (x) - f_2 (y) ) \chi_1'(y) E_{N, 1} (x, y) \tilde{I}_N (x, y) | \leq C \log(N) f_1(x) f_2(x)
\eeq
which contributes an acceptable error after integration in $x$ due the convolution property of the Poisson kernel.
 \qed

\bel
For the terms \eqref{eqn:even-cov-6} and \eqref{eqn:even-cov-7} we have,
\begin{align}
&\left|  \int(f_1 (x) - f_1 (y) )(f_2 (x) - f_2 (y) ) ( \partial_y  K_N (x, y) ) F_{N, 1} (y, x) \chi_1 (x) \chi_1 (y) \d x \d y \right| \notag  \\
\leq &   C \log(N) \left( \frac{|E_1 - E_2 | }{ \eta  ( ( E_1 - E_2)^2 + \eta^2 )} + \frac{1}{ ( E_1 - E_2)^2 + \eta^2 } \right)
\end{align}
and
\begin{align}
 &  \left|  \int(f_1 (x) - f_1 (y) )(f_2 (x) - f_2 (y) ) [ \partial_y E_{N, 1} (x, y) F_{N, 1} (y, x)  ] \d x \d y  \right| \notag \\
 \leq &   C \log(N) \left( \frac{|E_1 - E_2 | }{ \eta  ( ( E_1 - E_2)^2 + \eta^2 )} + \frac{1}{ ( E_1 - E_2)^2 + \eta^2 } \right)
 \end{align}
\eel
\proof We start with the first estimate. Partial integration in $y$ leads to a term with $K_N (x, y) E_{N, 1} (y, x)$  as well as one where the derivative hits the $f_i$'s and $\chi_1$. The former term already appeared in Lemma \ref{lem:even-3} and was estimated there. The  term where the derivative hits the $f_i$'s can be estimated in the same manner as Lemma \ref{lem:even-2} using the estimate of Lemma \ref{lem:even-cov-EN} for $F_{N, 1}$ (instead of the estimate $| \tilde{I}_N (x, y) | \leq C  \log(N)$ which was used in that proof). The integrand that arises when the derivative hits $\chi_1 (y)$ is bounded by $C f_1(x) f_2(x) $ which again contributes an acceptable error after integration in $x$.

For the second estimate stated in the lemma, we again integrate by parts. The term arising when the derivative hits $F_{N, 1}$ is  $E_{N, 1} (x, y) E_{N, 1} (y, x)$, which was estimated  in Lemma \ref{lem:even-3}. When the derivative hits the $f_i$'s we can proceed similarly to Lemma \ref{lem:even-2}. \qed

\bel
For the term \eqref{eqn:even-cov-8} we have,
\begin{align}
& \left| \int(f_1 (x) - f_1 (y) )(f_2 (x) - f_2 (y) ) [ \partial_y S_N (x, y)  \sg (y-x )  ] \d x \d y  \right| \notag\\
\leq & C \log(N) \left( \frac{|E_1 - E_2 | }{ \eta  ( ( E_1 - E_2)^2 + \eta^2 )} + \frac{1}{ ( E_1 - E_2)^2 + \eta^2 } \right)
\end{align}
\eel 
\proof We again integrate by parts. Since we have $\sg (y-x)$ there is a boundary term at $y=x$. However, this boundary term is $0$ due to the fact that $S_N$ is a continuous function, and the factors $(f_i(x) - f_i(y) )$ which vanish when $x=y$. There are then terms that arise when the derivative hits the $f_i$'s and when it hits $\chi_1$. Since $|S_N(x, y) | \leq C\log(N) |x-y|^{-1}$ and $| \sg (y-x ) | \leq 1$, these kinds of terms have been estimated in the previous lemmas.  \qed

\subsection{Covariance for GOE, $N$ odd}
In this section we prove the following.
\bep \label{prop:goe-odd}
The conclusions of Proposition \ref{prop:goe-even} hold also for GOE matrices of odd dimension.
\eep
We again use Lemma \ref{lem:goe-edge-simplification} and so in this section we estimate the first term on the right side of \eqref{eqn:goe-edge-simplification}. The exact formula for $R_N (x, y)$ in the case of odd dimension is similar to the case of even dimension, except for a few terms. After introducing the formulas, we just show how to handle the extra terms that arise, as the rest is the same as the case of even dimension.

The work \cite{pastur2011eigenvalue} does not contain formulas for $N$ odd, so we use Mehta's book  \cite{mehta2004random}. We will list the formulas there, rescale them to match our convention and then relate them to the quantities listed in the $N$-even case. We will need a few new estimates for some of the terms that appear.

The formulas in \cite{mehta2004random} are given in terms of the Hermite functions,
\beq
\varphi_n (x) = \frac{1}{( 2^n n! \sqrt{ \pi } )^{1/2}} \e^{- x^2/2} H_n (x) = \e^{-x^2/2} h_n (x),
\eeq
with $H_n$ and $h_n$ as in \eqref{eqn:unnorm-poly} and \eqref{eqn:norm-poly}, and
 are orthonormal.  Compared with \eqref{eqn:hermite-unscaled} we have,
\beq
\psi_n (x) = 2^{-1/4} \varphi_n ( x 2^{-1/2} ).
\eeq
From \cite[(7.2.20)-(7.2.27)]{mehta2004random}  (which defines $\tilde{K}$ below), and \cite[(7.2.31)]{mehta2004random} (which defines the cluster function in terms of $\tilde{K}$, but see also \cite[(7.2.30)]{mehta2004random} for the correlation functions as well as Chapter 5.1 for the definition of the determinant of quarternions, and the formula for the covariance in terms of correlation functions above) and \cite[(A.10.1)]{mehta2004random} for the Christoffel-Darboux formula in this scaling, we have,
\begin{align}
\Cov ( \tr f(H), \tr g (H) ) = \frac{1}{4} \int_{\rr^2} & ( f(x (2/N)^{1/2}) - f(y(2/N)^{1/2}) ) ( g(x(2/N)^{1/2}) - g(y(2/N)^{1/2}) ) \notag\\
& \times \tr \left(  \tilde{K}_N (x, y) \tilde{K}_N (y, x) \right) \d x \d y
\end{align}
where,
\begin{align}
\tilde{K}_N (x, y) := \left( \begin{matrix} \tilde{S}_N (x, y) & \tilde{D}_N (x, y) \\ \tilde{J}_N (x, y)  & \tilde{S}_N (y, x) \end{matrix}\right)
\end{align}
where
\begin{align}
\tilde{S}_N (x,y) &:= \left( \frac{N}{2} \right)^{1/2} \frac{ \varphi_N (x ) \varphi_{N-1} (y) - \varphi_N (y) \varphi_{N-1} (x) }{ x-y }  \notag\\
&+ \varphi_{N-1} (x) \left(  \left( \frac{N}{2} \right)^{1/2} \sg \varphi_{N} (y) + \frac{1}{ \int \varphi_{N-1} (t) \d t } \right),  \notag\\
\tilde{D}_N (x, y) &:= - \frac{\del}{ \del y } \tilde{S}_N (x, y), \notag\\
\tilde{J}_N (x, y) &:= \int \sg ( x - t ) \tilde{S}_N (t , y ) \d t \notag\\
&- \sg (x-y) - \sg \varphi_{N-1} (y) \left(  \int \varphi_{N-1} (t) \d t \right)^{-1} .
\end{align}
So, rescaling, we obtain 
\beq
\Cov ( \tr f (H) , \tr g (H) ) = \frac{1}{4} \int_{\rr^2} ( f(x) - f(y) ) ( g (x) - g (y) ) \tr \hat{K}_N (x, y) \hat{K}_N (y, x) \d x \d y
\eeq
where
\beq
\hat{K}_N (x, y) = \left( \begin{matrix} \hat{S}_N (x, y) &  (2/N)^{1/2} \hat{D}_N (x, y) \\ (N/2)^{1/2} \hat{J}_N (x, y) & \hat{S}_N (y, x) \end{matrix} \right) 
\eeq
and
\begin{align}
\hat{S}_N (x, y) & := K_N (x, y) + \hat{E}_{N, 1} (x, y), \notag\\
\hat{E}_{N, 1} (x, y) &:= \frac{N}{2} \phiNm (x) \sg \phiN (y) + \phiNm (x) ( \int \phiNm (t) \d t )^{-1}, \notag\\
\hat{D}_N (x, y) &:= - \frac{ \del}{ \del y} \hat{S}_N (x, y), \notag\\
\hat{J}_N (x, y) &:= \int \sg (x - t) \hat{S}_N (t, y) \d t \notag\\
& - \eps (x-y) - \sg \phiNm (y) \left( \int \phiNm (t) \d t ) \right)^{-1}.
\end{align}
So we have in this case that twice the cluster function is given by,
\beq
R_N (x, y) = 2 \hat{S}_N (x, y) \hat{S}_N (y, x) + \hat{J}_N (x, y) \hat{D}_N (y, x) + \hat{J}_N (y, x) \hat{D}_N (x, y).
\eeq
We now prove some estimates on the terms arising above.
\bel \label{lem:cov-goe-odd-1}
Let $\delta >0$. For $x, y$ satisfying $|x|, |y| < 2 - \delta$ we have,
\beq
| \hat{E}_{N, 1} (x, y) | \leq C.
\eeq
\eel
\proof Let $N = 2m + 1$. By \cite[(6.2.8)]{pastur2011eigenvalue} we have,
\beq 
| \hat{E}_{N, 1} (x, y) - \hat{E}_{N, 1} (x, 0) | \leq C.
\eeq
For the latter term, we have
\beq \label{eqn:EN-0}
\hat{E}_{N, 1} (x, 0) =  \phiNm (x)  \left(  \frac{N}{2} \sg \phiN (0)+( \int \phiNm (t) \d t )^{-1} \right),
\eeq
by definition. We now estimate the factor appearing on the RHS, showing it is $\O (1)$. 
We have,
\begin{align}
\frac{N}{2} \sg \phiN (0) &= - \frac{N}{2} \int_0^\infty \phiN (t) \d t \notag\\
&= - \left( \frac{N}{2} \right)^{3/4} \int_0^\infty \varphi_{2m+1} (t) \d t \notag\\
&= - \frac{m^{1/2}}{2^{1/2}} + \O \left( 1 \right) .
\end{align}
In the first line we used the fact that $\phiN$ is an odd function. In the second line we used the definition of $\phiN$ in terms of the Hermite functions. In the last line we used Lemma \ref{lem:herm-int-2}. 
For the other term we have,
\begin{align}
\int \phiNm (t) \d t &= \left( \frac{N}{2} \right)^{-1/4} \int \varphi_{2m} (t) \d t \notag \\
&= \frac{2^{1/2}}{m^{1/2}} \left( 1 + \O ( m^{-1} ) \right)
\end{align}
We used Lemma \ref{lem:herm-int-1} in the second line. 
The claim follows after using the fact that the Hermite functions are bounded in the domain in question. \qed

\bel
Let $\delta >0$. We have for $x, y $ satisfying $|x|, |y| < 2 - \delta$ that.
\beq
\left| \int \sg (x -t ) \hat{E}_{N, 1} (t, y) \right| \leq C
\eeq
and
\beq
\left| \sg \phiNm (y)  \left( \int \phiNm (t) \d t ) \right)^{-1} \right| \leq C
\eeq
\eel
\proof By \cite[(6.2.8)]{pastur2011eigenvalue} we have
\beq
| \sg \phiN (y) - \sg \phiN (0) | \leq C N^{-1}
\eeq
and so it suffices to prove the first estimate with $y=0$ due to the fact that scaled Hermite functions $\phiN$ and $\phiNm$ have bounded $L^1$ norms. For $y=0$, we have the explicit form \eqref{eqn:EN-0} for which the factor on the RHS is $\O(1)$, as shown in the previous proof. The first estimate follows.

For the  second estimate, note that the first estimate of \cite[(6.2.9)]{pastur2011eigenvalue} (this applies since here $N$ is odd and there $N$ is even) implies that
$| \sg \phiNm (y) | \leq C N^{-1}$. The previous proof shows that the other factor is $\O (N^{1/2})$ and so we conclude the proof. \qed

\bel \label{lem:cov-goe-odd-2}
Let $|x|, |y| < 2 - \delta$. Then,
\beq
\left| \int \sg (x-t) \tilde{S}_N (t, y) \d t \right| \leq C \log(N) ,
\eeq
and so,
\beq
| \hat{J}_N (x, y) | \leq C \log(N).
\eeq
\eel
\proof First,
\beq
\int_{ t : |t-y| > \delta/2} | \tilde{S}_N (t, y) | \d t \leq C ( \| \phiN \|_1 + \| \phiNm \|_1 ) \leq C ,
\eeq
by the asymptotics  \eqref{eqn:phiN-asymp}.  Next, the proof of Lemma \ref{lem:cov-Itildeest} shows that,
\beq
\int_{ t : |t-y| \leq \delta/2} | K_N (t, y) | \d t \leq C \log(N) ,
\eeq
and the claim follows. \qed

\noindent{\bf Proof of Proposition \ref{prop:goe-odd}.} If suffices to estimate the first term on the right side of \eqref{eqn:goe-edge-simplification}. For this term we see from the above that we have the formula,
\begin{align}
 & \frac{1}{4} \int_{\rr^2} (f_1 (x) - f_1 (y) ) (f_2 (x) - f_2 (y) ) \chi_1 (x) \chi_1 (y) R_N (x, y) \d x \d y \notag\\
= &  \frac{1}{2} \int_{\rr^2}  (f_1 (x) - f_1 (y) ) (f_2 (x) - f_2 (y) ) \hat{S}_N (x, y) \hat{S}_N (y, x) \chi_1 (x) \chi_1 (y) \d x \d y \label{eqn:cov-goe-odd-1} \\
- & \frac{1}{2} \int_{\rr^2} (f_1 (x) - f_1 (y) ) (f_2 (x) - f_2 (y) )  ( \del_y \hat{S}_N (x, y) ) \hat{J}_N (y, x) \chi_1 (x ) \chi_1 (y) \d x \d y. \label{eqn:cov-goe-odd-2}
\end{align}
Given the estimates that we have established above, estimating the above terms is very similar to the proof of Proposition \ref{prop:goe-even}. Given Lemma \ref{lem:cov-goe-odd-1}, the term \eqref{eqn:cov-goe-odd-1} is estimated in the exact same fashion as the terms \eqref{eqn:even-cov-1}, \eqref{eqn:even-cov-3} and \eqref{eqn:even-cov-4}. For the term \eqref{eqn:cov-goe-odd-2}, we integrate by parts in $y$. When the derivative hits $\hat{J}_N (y, x)$ we use $\del_y \hat{J}_N (y, x) = \hat{S}_N (y, x)$ and see that this term is identical to \eqref{eqn:cov-goe-odd-1}. When the derivative hits the other terms we see from Lemma \ref{lem:cov-goe-odd-2} that the integrand is bounded by
\beq
C \log(N)^2 | ( \del_y [ \chi_1 (y) (f_1 (y) - f_1 (x) ) ( f_2 (x) - f_2 (y) ) ] ) (x-y)^{-1} |.
\eeq
Estimation of the integral of this quantity was handled in the proof of Lemma \ref{lem:even-2}. \qed

\subsection{Expectation} \label{sec:gauss-expect}

\bel
Let $f$ be a function, and $H$ a matrix from the GOE. Then,
\begin{align}
\ee[ \tr f (H) ] &= N \int f(x) \rhosc (x) \d x - \frac{1}{ 2 \pi }\int f(x) \frac{1}{\sqrt{4-x^2} } \d x \notag \\
&+ \frac{ f(2) + f(-2) }{4} + N^{\eps} \O ( N^{-1} \| f'' \|_{1, w}  )
\end{align}
\eel 
\proof This follows immediately from Lemma \ref{lem:wig-expect}.
 \qed

\section{Proof of main theorems} \label{sec:main}

\bet \label{thm:GDE-tech}
 Let $ \varphi \in H^{1/2+s}$ for $s >0$. Let $H$ be a Gaussian divisible ensemble of the form
 \beq
 H = \e^{-T/2} W + \sqrt{1 - \e^{-T} } G
 \eeq
 where $T = N^{-\omega}$ and $\omega>0$ satisfies,
 \beq
 \omega < \frac{s}{100}.
 \eeq
 Then,
 \beq
 \Var ( \tr ( \varphi (H) ) ) \leq C \| \varphi \|_{H^{1/2+s}}^2 .
 \eeq
\eet
\proof Fix a small $\mfa >0$. Let $\chi$ be a smooth bump function that is 1 on the interval $[-1, 1]$ and $0$ outside of $[-2, 2]$. Define
\beq
\hat{\psi}_1 ( \xi) = \chi( \xi / N^{\mfa} ) \hat{\varphi} ( \xi),
\eeq
and 
\beq
\hat{\psi}_2 ( \xi ) = ( 1- \chi ( \xi / N^{\mfa} ) ) \sum_{|k - \mfa \log_2 (N) | < 100 } \hat{\varphi}_k ( \xi )
\eeq
and
\beq
\hat{\psi}_3 ( \xi ) = \sum_{k \geq 50 + \mfa \log_2 (N) } \hat{\varphi}_k ( \xi).
\eeq
By the properties of the Littlewood-Paley decomposition,
\beq
\varphi = \psi_1 + \psi_2 + \psi_3.
\eeq
By \cite[Lemma 2.1, p. 52]{bahouri2011fourier}, we have that
\beq
\| \psi_1 \|_{C^2  } + \| \psi_2 \|_{C^2} \leq C N^{3 \mfa/2} \| \varphi \|_2.
\eeq
We choose $\mfa = 1/10$.  Now, if $p$ is the inverse Fourier transform of $\chi$ we have,
\beq
\psi_1 (x) = N^{\mfa} \int p ( N^{\mfa} (x - y) ) \varphi (y) \d x.
\eeq
If $x \notin (-2+\kappa/2 , 2 - \kappa/2)$ we see from the fact that,
\beq
|p^{(m)} (x) | \leq \frac{C_{n, m}}{ 1 + |x|^n}
\eeq
that,
\beq
| \psi^{(m)}_1 (x) | \leq C N^{-100} \frac{1}{ 1 + |x|^{10}} \| \varphi\|_2.
\eeq
for $0 \leq m \leq 5$ and $x \notin (-2+\kappa/2, 2 - \kappa/2)$. By the properties of the Littlewood-Paley decomposition, we see that 
\beq
\hat{\psi}_2 ( \xi ) = (1 - \chi ( \xi / N^{\mfa}  ) ) q ( \xi / N^{\mfa} ) \hat{ \varphi} ( \xi)
\eeq
for some smooth $q$ that is of compact support.  Explicitly, 
\beq
q (\xi ) = \sum_{|k| \leq 100 } \hat{ \omega} ( 2^{-k } \xi ).
\eeq
 By the same argument as above, we see that
\beq
| \psi^{(m)}_2 (x) | \leq C N^{-100} \frac{1}{ 1 + |x|^{10}} \| \varphi\|_2
\eeq
for $0 \leq m \leq 5$ and $x \notin (-2+\kappa/2, 2 - \kappa/2)$ as well. 
Therefore by Theorem \ref{thm:var-global-2} we have that,
\beq
| \Var ( \tr \psi_1 (H) ) - V ( \psi_1 ) | + | \Var ( \tr \psi_2 (H) ) - V ( \psi_2 ) |\leq N^{-1} N^{3 \mfa + \eps} \| \varphi \|_2^2.
\eeq
Clearly
\beq
\| \psi_1 \|_2 + \| \psi_2 \|_2 \leq C \| \varphi \|_2.
\eeq
We see from the above estimates that
\beq
| V ( \psi_i ) | \leq C \| \varphi \|_2^2 + C \int_{ ( -2 + \kappa/4, 2 - \kappa/4)^2} \frac{ ( \psi_i (x) - \psi_i (y) )^2}{ (x-y)^2} \d x \d y.
\eeq
The last term is bounded above by the integral over $\rr^2$ which is, up to constants, the norm $ \| \psi_i \|_{\dot{H}^{1/2}}$. Clearly,
\beq
\| \psi_i \|_{\dot{H}^{1/2}} \leq C \| \varphi \|_{H^{1/2+s}}
\eeq
and so 
\beq
 \Var ( \tr ( \psi_1 (H) ) )  +\Var ( \tr ( \psi_2 (H) ) )  \leq C \| \varphi \|_{H^{1/2+s}} .
\eeq
Let
\beq
M = 10s^{-1} \log_2(N).
\eeq
and
\beq
\psi_4 := \sum_{k \geq M} \varphi_k.
\eeq
By Theorem \ref{thm:main-uv} we have,
\beq
\Var ( \tr ( \psi_4 (H) ) ) \leq C N^{-1} \| \varphi \|_{H^{1/2+s}}. 
\eeq
Consider now $k$ such that
\beq
\mfa \log_2(N) -100 \leq k \leq M.
\eeq
By Theorems \ref{thm:meso-var} and \ref{thm:sub} we have,
\beq
\Var ( \tr \varphi_k (H) ) \leq C N^{\eps+\omega} 2^{k} \| \varphi_k \|_2^2 + N^{-100} \| \varphi \|_2^2
\eeq
for any $\eps >0$. We may now choose $\eps >0$ so that,
\beq
\eps + \omega <  \frac{\mfa s}{2}.
\eeq
Then
\beq
N^{\eps+\omega} 2^{k} \| \varphi_k \|_2^2 \leq N^{-\mfa s} \| \varphi \|^2_{H^{1/2+s}}.
\eeq
Hence
\beq
\Var ( \tr ( \psi_3 - \psi_4) ( H) ) \leq C \log(N) N^{-\mfa s }  \| \varphi \|^2_{H^{1/2+s}} .
\eeq
The 
claim follows.\qed


\subsection{Proof of Theorem \ref{thm:main-1}}

We recall the following, \cite[Proposition 4.1]{erdHos2011universality}.
\bep \label{prop:reverse-heat}
Let $u$ be a Wigner-smooth probability density as in Definition \ref{def:smooth-density}. For any $K>0$ there is an $\alpha_K >0$ so that for any positive $t < \alpha_K$ there exists a probability density $p_t$ with mean $0$ and variance $1$ so that
\beq
\int | \e^{t A} p_t - u | \d x \leq C t^K,
\eeq
where $A = \frac{1}{2} \del_{x}^2 - \frac{x}{2} \del-x$ is the generator of the Ornstein-Uhlenbeck flow. Moreover, $g_t$ is Wigner-smooth with uniform constants depending on those for $u$.
\eep
The only difference from the above and the construction in \cite{erdHos2011universality} is the slightly weaker assumption on the decay, but that this holds can be checked by the explicit form of the short construction given there. 

Now consider a smooth Wigner matrix $H$ and a function $ \varphi \in H^{1/2+s}$. Denote the density of its  off-diagonal entries by $u_0$ and of the diagonal entries by $u_d$. Let,
\beq
M = 10s^{-1} \log_2(N),
\eeq
and
\beq
T = N^{-\omega}, \qquad \omega = \frac{s}{100}.
\eeq
Let,
\beq
\psi_1 = \sum_{k \geq -1}^M \varphi_k, \qquad \psi_2 = \sum_{k > M} \varphi_k.
\eeq
We see from Theorem \ref{thm:main-uv} that
\beq
\Var ( \tr \psi_2 ( H) ) \leq C N^{-1} \| \varphi \|_{H^{1/2+s}}^2.
\eeq
We choose $K$ sufficiently large so that with $t = T$, Proposition \ref{prop:reverse-heat} guarantees existence of a densities $p_o,p_d$ so that,
\beq \label{eqn:smooth-match}
\int | \e^{t A} p_o - u_o | \d x +\int | \e^{t A} p_d - u_d | \d x \leq C \frac{1}{N^{100}} \frac{1}{2^{10M}}.
\eeq
Let $W_1$ be the Wigner ensemble with off-diagonal densities $p_0$ and diagonal densities $p_d$ and $W = \e^{-t/2} W_1 + \sqrt{1 - \e^{-t}} G$ with $G$ GOE so that $W$ has densities $\e^{t A} p_d$ and $\e^{tA} p_o$.  We apply Theorem \ref{thm:GDE-tech} to $W$.$^1${\let\thefootnote\relax\footnotetext{$1$. Note that strictly speaking, the third and fourth cumulants of $W_1$ may not satisfy \eqref{eqn:rs-cumu} but, by \eqref{eqn:smooth-match}, and since $H$ does satisfy \eqref{eqn:rs-cumu}, $W_1$ satisfies it up to, say, an error of $\O (N^{-100})$ which does not affect the proof.}}
For $\eta_1, \eta_2 \geq 2^{-M}$ it follows from the Lipschitz continuity of the resolvent that,
\begin{align}
 \bigg| &\Cov \left( \Im \left[ \tr (H - x - \i \eta_1 )^{-1} \right], \Im \left[ \tr (H - y - \i \eta_2 )^{-1} \right] \right) \notag\\
 & - \Cov \left( \Im \left[ \tr (W - x - \i \eta_1 )^{-1} \right], \Im \left[ \tr (W - y - \i \eta_2 )^{-1} \right] \right) \bigg| \leq C N^{-10}.
\end{align}
From the decomposition,
\begin{align}
&\Var ( \sum_{k \leq M} \tr \varphi_k (H ) )  \notag\\
= & \sum_{l, k \leq M} \int_{\rr^2} g_l(x) g_k (y)  \Cov \left( \Im \left[ \tr (H - x - \i \eta_l )^{-1} \right], \Im \left[ \tr (H - y - \i \eta_k )^{-1} \right] \right) \d x \d y ,
\end{align}
the estimate
\beq
\left| \Var \left(  \sum_{k \leq M} \tr \varphi_k (H ) \right) - \Var \left(  \sum_{k \leq M} \tr \varphi_k (H ) \right) \right| \leq N^{-2} \| \varphi \|^2_{2}
\eeq
follows. This yields the estimate for functions in $H^{1/2+s}$. \qed 

\section{Fine estimates of characteristic functions} \label{sec:fine-lss}

Let $f$ be a smooth function supported in $[-5, 5]$ and $H$ a Wigner matrix.  Recall the quasi-analytic extension of $f$ given by,
\beq
\tilf (x + \i y ) := (f (x) + \i y f'(x) ) \chi (y)
\eeq
where $\chi$ is an even bump function that is $1$ for $|y| < 1$ and $0$ for $|y| > 2$. 
We fix $\mfa >0$ for the rest of the section and define
\beq
\ea ( \xi ) := \exp \left[ \i \xi \left(  \frac{1}{ \pi }\int_{\Oma} \del_{\bar{z}} \tilf (z) N ( m_N (z) - \ee[ m_N(z) ] ) \d x \d y \right) \right] ,
\eeq
and
\beq
\psia ( \xi ) := \ee[ \ea ( \xi ) ],
\eeq
where $\Oma$ was defined in \eqref{eqn:Oma-def}. Lemma \ref{lem:HS-cutoff} implies that $\psia$ is a good approximation for the characteristic function of the linear statistic.

\subsection{Self-consistent equation for Stein's method}

We will apply Stein's method to calculate $\psia ( \xi)$. The main tool technical estimate that we will use in our application of Stein's method is the following proposition.
\bep \label{prop:stein-self} Let $\eps  >0$ and $\tau >0$.  For any $z = x + \i y $ with $z \in \D_{\tau, 1}$ we have,
\begin{align}
&(z+2 \msc (z) ) N \ee[ \ea ( \xi ) (m_N (z) - \ee[ m_N (z) ] ) ] \notag\\
= &  - \frac{ 2 \i \xi}{N} \frac{1}{ \pi} \int_{\Oma} \del_{\bar{w}} \tilf (w) \psi_a \del_w \frac{ \msc (z) - \msc(w) }{ z- w } \d u \d v  \notag\\
 + & \psia \frac{ 4 s_3 \msc(z) \i \xi}{\pi N^{1/2}} \int_{\Oma} \del_{\bar{w}} \tilf (w)( \msc (z) + \msc(w)) \msc'(w) \d u \d v \notag \\
- &  \frac{ 2 \xi^2 s_3}{ \pi^2 N^{1/2}} \psia \msc(z) \left( \int_{\Oma} \del_{\bar{w}} \tilf (w) \msc' (w) \d u \d v \right)^2 \notag \\
 - & \frac{ - 2 \xi \i s_4}{ \pi}  \psia\int_{\Oma} \del_{\bar{w}} \tilf (w) \msc'(w) \msc(w) \msc(z)^2 \d u \d v  \notag\\
+& N^{\eps} (1 + |\xi|^6) ( 1+ \| f'\|_1^6 ) \O ( N^{-1} ( 1 + \| f''\|_1 )y^{-1} + N^{-1} y^{-2} )
\end{align}
Above, $w = u + \i v$. 
\eep
The proof of the above proposition is largely similar to the proof of Proposition \ref{prop:cov-self}, in that it will be seen to be a consequence of
 a cumulant expansion. First, we have that whenever the entry-wise local semicircle law holds,
\beq \label{eqn:delea-1}
| \del_{ia}^k \ea ( \xi ) | \leq C (1 + | \xi| )^k N^{\eps}(1+  \| f' \|_1)^k
\eeq
with overwhelming probability for any $k$ \cite{meso}. This allows for,
\begin{align} \label{eqn:stein-exp}
&z \sum_{i} \ee[ \ea (G_{ii} (z) - \ee[ G_{ii} (z) ] ) ] \notag\\
= & \sum_{ia} \ee[ \ea (G_{ia}(z) H_{ia} - \ee[ G_{ia} (z) H_{ia} ])] \notag\\
= & \sum_{n=2}^{5} \sum_{ia} \frac{(1+\delta_{ia} )^{n-1} s_n}{N^{n/2} (n-1)!} \ee[ \del_{ia}^{n-1} \ea (G_{ia}(z)  - \ee[ \del_{ia}^{n-1} G_{ia} (z)  ] )]  \notag \\
+&N^{\eps} \O( (1+ | \xi | )^6 N^{-1} (1+  \| f' \|_1)^6)
\end{align}
In the next few subsections we will deal with each of the terms on the second last line of \eqref{eqn:stein-exp}. We note that many of the terms that appear when the derivatives $\del_{ia}$ act on $\ea$ have been dealt with in Section \ref{sec:global-self} in the proof of Proposition \ref{prop:cov-self} and so we will refer to the results in that section as appropriate.  As in that proof, we will refer to the term in the second line of \eqref{eqn:stein-exp} with $n=k$ as the $k$th order term. 

Throughout the proof we will use the notation $z = x \pm \i y$ as well as $w = u \pm \i v$ and $w_i = u_i \pm \i v_i$. We will often use the fact that since $f$ is assumed to be supported in $[-5, 5]$ that $\| f \|_1  \leq C \| f' \|_1 \leq C^2 \| f'' \|_1$ to simplify some error terms.  In applying the entrywise local law \eqref{eqn:entry-ll} we will often use that $\Psi(z) \leq C (N y )^{-1/2}$ when $y \geq N^{-1}$.
\subsubsection{Second order term} 
\bel
For any $\eps >0$ we have,
\begin{align}
& \sum_{ia} \frac{ (1+ \delta_{ia} ) }{N} \ee[ \del_{ia}  \ea ( G_{ia}(z) - \ee[ \del_{ia} G_{ia} (z)]  \notag\\
= & -2 N \msc (z) \ee[ \ea (m_N(z) - \ee[ m_N (z) ] ) ]  \notag\\
+ & - \frac{ 2 \i \xi}{N} \frac{1}{ \pi} \int_{\Oma} \del_{\bar{w}} \tilf (w) \psi_a \del_w \frac{ \msc (z) - \msc(w) }{ z- w } \d u \d v  \notag\\
+ & N^{\eps} \O ( N^{-1} (1+ \| f'' \|_1) y^{-1}  + N^{-1} y^{-2} (1 + \| f' \|_1 ) ).
\end{align}
\eel
\proof We write,
\begin{align}
&\sum_{ia} \frac{ (1+ \delta_{ia} ) }{N} \ee[ \del_{ia}  \ea ( G_{ia}(z) - \ee[ \del_{ia} G_{ia} (z)] \notag\\
=& -N  \ee[ \ea (m_N (z)^2 - \ee[ m_N(z)^2]) ] \notag\\
- &  \ee[ \ea ( \del_z m_N (z) - \ee[ \del_z m_N (z) ] ) ]\notag\\
+& \frac{1+\delta_{ia}}{N} \sum_{ia} \ee[ (\del_{ia} \ea) G_{ia} (z) ]
\end{align}
Using the local law \eqref{eqn:ll} we have,
\begin{align}
 -N  \ee[ \ea (m_N (z)^2 - \ee[ m_N(z)^2]) ] &= -2 N \msc (z) \ee[ \ea (m_N(z) - \ee[ m_N (z) ] ) ]  \notag\\ 
 &+ N^{\eps} \O ( N^{-1} y^{-2} ),
\end{align}
and
\beq
\ee[ \ea ( \del_z m_N (z) - \ee[ \del_z m_N (z) ] ) ] =  N^{\eps} \O ( N^{-1} y^{-2} ).
\eeq
The last term is,
\begin{align}
 & \frac{1+\delta_{ia}}{N} \sum_{ia} \ee[ (\del_{ia} \ea) G_{ia} (z) ] \notag\\
 = &- \frac{ 2 \i \xi}{N} \frac{1}{ \pi} \int_{\Oma} \del_{\bar{w}} \tilf (w) \ee[ \ea \del_w \tr G(z) G(w) ] \d u \d v.
\end{align}
In order to deal with this term, we apply the estimate \eqref{eqn:GG-est} derived in the proof of Lemma \ref{lem:cov-2o}. We first have, denoting $w= u + \i v$,
\begin{align}
 & \frac{1}{N} \int_{\Oma} | (f(u) + \i f'(u) v ) \chi'(v)| \left| \ee \left[ \ea \left( \del_w \tr G(w) G(z) - \del_w  \frac{ \msc(z) - \msc (w) }{ z -w } \right) \right] \right| \d u \d v \notag\\
 \leq & C \int_{\Oma} | (f(u) + \i f'(u) v ) \chi'(v)| N^{-1} ( y^{-2} v^{-1} + v^{-2} ( y  + v)^{-1} ) \d u \d v \notag\\
\leq &  C N^{\eps} ( 1  + \| f' \|_1 ) N^{-1} y^{-2}
\end{align}
Let $ t = N^{\eps} ( 1 + \| f' \|_1 ) \| f'' \|_1^{-1} \wedge 1$. We have,
\begin{align}
& \int_{ N^{\mfa-1} < |v| < t } | f''(u ) v  \chi (v) | \left| \ee \left[ \ea \left( \del_w \tr G(w) G(z) - \del_w  \frac{ \msc(z) - \msc (w) }{ z -w } \right) \right] \right| \d u \d v  \d u \d v \notag\\
\leq & C N^{\eps} \| f'' \|_1 \int_{N^{-1} < v < t } v N^{-1} ( y^{-2} v^{-1} + y^{-1} v^{-2} ) \d v \leq C N^{3 \eps-1} ( (1 + \| f' \|_1 ) y^{-2} + \| f'' \|_1 y^{-1} ).
\end{align}
By integration by parts, the Cauchy-Riemann equations (i.e., $\del_z F = \del_x F$ for holomorphic $F$) and the Cauchy integral formula we have,
\begin{align}
& \left| \int_{ |v| > t } f''(u) v \chi (v) \ee \left[ \ea \left( \del_w \tr G(w) G(z) - \del_w  \frac{ \msc(z) - \msc (w) }{ z -w } \right) \right]  \d u \d v  \right| \notag\\
= & \left| \int_{ |v| > t } f'(u) v \chi (v) \ee \left[ \ea \left( \del^2_w \tr G(w) G(z) - \del^2_w  \frac{ \msc(z) - \msc (w) }{ z -w } \right) \right]  \d u \d v  \right| \notag\\
\leq & CN^{\eps} \| f'\|_1 \int_{t < v < 2} N^{-1}v  (y^{-2} v^{-2} + y^{-1} v^{-3} ) \d v \leq C N^{3 \eps-1}  ( 1 + \| f ' \|_1 )y^{-2} +  (1+ \| f'' \|_1) y^{-1} ).
\end{align}
Therefore,
\begin{align} \label{eqn:stein-ibp}
- & \frac{ 2\i \xi}{N} \frac{1}{ \pi} \int_{\Oma} \del_{\bar{w}} \tilf (w) \ee[ \ea \del_w \tr G(w) G(w) ] \d u \d v \notag\\
= & - \frac{ 2\i \xi}{N} \frac{1}{ \pi} \int_{\Oma} \del_{\bar{w}} \tilf (w) \psi_a \del_w \frac{ \msc (z) - \msc(w) }{ z- w } \d u \d v \notag\\
+ & N^{3\eps} \O ( ( 1 + \| f ' \|_1 )N^{-1} y^{-2} + N^{-1} (1+ \| f'' \|_1) y^{-1}).
\end{align}
This completes the proof. \qed

\subsubsection{Third order term}
\bel
We have,
\begin{align}
& \sum_{ia} \frac{ (1+ \delta_{ia} )^2 }{N^{3/2}} \ee[ \del^2_{ia}  \ea ( G_{ia}(z) - \ee[ \del^2_{ia} G_{ia} (z)]  \notag\\
= & \psia \frac{ 8 \msc(z) \i \xi}{\pi N^{1/2}} \int_{\Oma} \del_{\bar{w}} \tilf (w) (\msc(w) + \msc(z) ) \msc'(w) \d u \d v \notag \\
- &  \frac{ 4  \xi^2}{ \pi^2 N^{1/2}} \psia \msc(z) \left( \int_{\Oma} \del_{\bar{w}} \tilf (w) \msc' (w) \right)^2 \d u \d v \notag \\
+ & N^{\eps} (1+ |\xi| ) \O ( N^{-1} ( \| f'' \|_1 + 1 ) y^{-1} ) 
+  N^{\eps}|\xi|^2 \O ( N^{-1} y^{-1/2} (1+ \| f' \|_1) ( 1+  \| f'' \|_1) )  \notag\\
+& N^{\eps} \O ( N^{-1} y^{-3/2} )
\end{align}
\eel
\proof 
We write,
\begin{align}
&\sum_{ia} \frac{ (1+ \delta_{ia} )^2 }{N^{3/2}} \ee[ \del^2_{ia}  \ea ( G_{ia}(z) - \ee[ \del^2_{ia} G_{ia} (z)])] \notag\\ 
= &\sum_{ia}  \frac{ (1+ \delta_{ia} )^2 }{N^{3/2}}\ee[ \ea ( \del_{ia}^2 G_{ia} (z) - \ee[ \del_{ia}^2 G_{ia} (z) ] ) ]  \label{eqn:stein-3o1} \\
+ & \sum_{ia}  \frac{2 (1+ \delta_{ia} )^2 }{N^{3/2}}\ee[ \del_{ia} \ea \del_{ia} G_{ia} (z) ]  \label{eqn:stein-3o2}\\
+ & \sum_{ia} \frac{ (1+ \delta_{ia} )^2 }{N^{3/2}} \ee[ ( \del^2_{ia} \ea ) G_{ia} (z) ]. \label{eqn:stein-3o3}
\end{align}
We start with \eqref{eqn:stein-3o1}. Using \eqref{eqn:del2Gia} this is split into,
\begin{align}
& \sum_{ia}  \frac{ (1+ \delta_{ia} )^2 }{N^{3/2}}\ee[ \ea ( \del_{ia}^2 G_{ia} (z) - \ee[ \del_{ia}^2 G_{ia} (z) ] ) ]  \notag\\
= & \frac{1}{N^{3/2}} \sum_{ia} 2  \ee[ \ea (  G_{ia} (z)^3 - \ee[ G_{ia} (z)^3 ] ) ] \notag\\
+ & \frac{1}{N^{3/2}} \sum_{ia} 6 \ee[ \ea ( G_{ia} (z) G_{ii} (z) G_{aa} (z) - \ee[  G_{ia} (z) G_{ii} (z) G_{aa} (z)  ] ) ] .
\end{align}
By similar logic as in the estimation of the term \eqref{eqn:cov-3o1} we obtain,
\beq
\frac{1}{N^{3/2}} \sum_{ia} 2  \ee[ \ea (  G_{ia} (z)^3 - \ee[ G_{ia} (z)^3 ] ) ]  = N^{\eps} \O ( N^{-1} y^{-3/2} ).
\eeq
Similar logic as in the estimation of the term \eqref{eqn:cov-3o2} yields,
\begin{align}
\frac{1}{N^{3/2}} \sum_{ia} 6 \ee[ \ea ( G_{ia} (z) G_{ii} (z) G_{aa} (z) - \ee[  G_{ia} (z) G_{ii} (z) G_{aa} (z)  ] ) ] = N^{\eps} \O( N^{-1} y^{-3/2} ).
\end{align}
This completes estimation of the term \eqref{eqn:stein-3o1}.  We turn to \eqref{eqn:stein-3o2}. It is split into,
\begin{align}
&\sum_{ia}  \frac{2 (1+ \delta_{ia} )^2 }{ N^{3/2}}\ee[ \del_{ia} \ea \del_{ia} G_{ia} (z) ]  \notag\\
= & \i \xi \sum_{ia}  \frac{4  }{\pi N^{3/2}} \int_{\Oma} \del_{\bar{w}} \tilf (w)  \ee[ \ea \del_w G_{ia} (w) (G_{ii} (z) G_{aa} (z) ) ] \notag\\
+ &\i \xi \sum_{ia}  \frac{4  }{\pi N^{3/2}} \int_{\Oma} \del_{\bar{w}} \tilf (w)  \ee[ \ea \del_w G_{ia} (w) G_{ia}(z)^2 ] 
\end{align}
Similar logic  to the estimation of \eqref{eqn:cov-3o3} yields that the following estimate holds with overwhelming probability,
\begin{align}
\frac{1}{N^{3/2}} \sum_{ia} \del_w G_{ia} (w) G_{ii} (z) G_{aa} (z)  =  N^{-1/2} \msc(z)^2 \msc'(w)
+ N^{\eps} ( N^{-1} v^{-3/2} y^{-1} ).
\end{align}
By direct integration, we obtain the estimate
\begin{align}
& \i \xi \sum_{ia}  \frac{4  }{N^{3/2}} \int_{\Oma} \del_{\bar{w}} \tilf (w)  \ee[ \ea \del_w G_{ia} (w) (G_{ii} (z) G_{aa} (z) ) ] \notag\\
= & \i \xi 4 N^{-1/2} \frac{1}{ \pi} \int_{\Oma} \del_{\bar{w}} \tilf (w) \msc(z)^2 \msc'(w) \psia +N^{\eps} | \xi |\O ( N^{-1} y^{-1} (1+\| f'' \|_1) ).
\end{align}
In a similar fashion, we use the logic used to estimate \eqref{eqn:cov-3o4} to find,
\begin{align}
&i \xi \sum_{ia}  \frac{4  }{\pi N^{3/2}} \int_{\Oma} \del_{\bar{w}} \tilf (w)  \ee[ \ea \del_w G_{ia} (w) G_{ia}(z)^2 ]  \notag\\
= &\i \xi 4 N^{-1/2} \frac{1}{ \pi} \int_{\Oma} \del_{\bar{w}} \tilf (w) \msc(z)^2 \msc'(w) \psia +N^{\eps} | \xi | \O ( N^{-1} y^{-1}(1+ \| f'' \|_1 ) ).
\end{align}
This completes the estimation of \eqref{eqn:stein-3o2}.  We split the term \eqref{eqn:stein-3o3} into,
\begin{align}
& \sum_{ia} \frac{ (1+ \delta_{ia} )^2 }{N^{3/2}} \ee[ ( \del^2_{ia} \ea ) G_{ia} (z) ] \notag\\
=& \sum_{ia} \frac{2 \i \xi }{ \pi N^{3/2} } \int_{\Oma} \del_{\bar{w}} \tilf (w) \ee[ \ea \del_w( G_{ii} (w) G_{aa} (w) )G_{ia} (z) ] \notag\\
+ &  \sum_{ia} \frac{2 \i \xi }{ \pi N^{3/2} } \int_{\Oma} \del_{\bar{w}} \tilf (w) \ee[ \ea \del_w( G_{ia} (w)^2 )G_{ia} (z) ] \notag\\ 
+ & \sum_{ia} \frac{ - 4 \xi^2}{ \pi^2 N^{3/2}} \int_{\Oma^2} \del_{\bar{w}_1} \del_{\bar{w}_2} \tilf (w_1) \tilf (w_2) \ee[ \ea \del_{w_1} G_{ia} (w_1) \del_{w_2} G_{ia} (w_2) G_{ia} (z) ]
\end{align}
Similar logic to the estimation of \eqref{eqn:cov-3o5} and \eqref{eqn:cov-3o6} shows the estimates,
\begin{align}
& \frac{1}{N^{3/2}} \sum_{ia} G_{ia} (z) G_{ii} (w) \del_w G_{aa} (w)\notag \\
= & \msc (w) \msc' (w) \msc(z) + N^{\eps} \O ( N^{-1} y^{-1/2} v^{-2} )
\end{align}
and
\begin{align}
& \frac{1}{N^{3/2}} \sum_{ia} G_{ia} (z) G_{ia} (w) \del_w G_{ia} (w) \notag \\
= & \msc (w) \msc' (w) \msc(z) + N^{\eps} \O ( N^{-1} y^{-1/2} v^{-2} )
\end{align}
holds with overwhelming probability. Therefore, by direct integration, we have:
\begin{align}
& \sum_{ia} \frac{2 \i \xi }{ \pi N^{3/2} } \int_{\Oma} \del_{\bar{w}} \tilf (w) \ee[ \ea \del_w( G_{ia} (w)^2 +G_{ii} (w) G_{aa} (w) )G_{ia} (z) ] \notag\\
 &= \frac{ 8 \i \xi}{ \pi N^{1/2}}  \int_{\Oma} \del_{\bar{w}} \tilf (w) \msc'(w) \msc(w) \msc(z) \psia +N^{\eps} | \xi | \O(N^{-1} (1+ \| f'' \|_1 ) y^{-1/2} ).
\end{align}
Denote $w_i = u_i \pm v_i $. 
Appying the entry-wise local law \eqref{eqn:entry-ll} we see that with overwhelming probability,
\begin{align}
&\frac{1}{N^{3/2}} \sum_{ia} \del_{w_1} G_{ia} (w_1) \del_{w_2} G_{ia} (w_2) G_{ia} (z) \notag\\
= & \frac{\msc(z) \msc'(w_1) \msc'(w_2)}{N^{1/2}} + N^{\eps} \O (N^{-1} y^{-1/2} v_1^{-3/2}v_2^{-3/2} ),
\end{align} 
and so by a similar argument  to that leading to the estimate \eqref{eqn:stein-ibp} we have,
\begin{align}
 & \sum_{ia} \frac{ - 4 \xi^2}{ \pi^2 N^{3/2}} \int_{\Oma^2} \del_{\bar{w}_1} \del_{\bar{w}_2} \tilf (w_1) \tilf (w_2) \ee[ \ea \del_{w_1} G_{ia} (w_1) \del_{w_2} G_{ia} (w_2) G_{ia} (z) ] \notag\\
 = & \frac{ - 4 \xi^2}{ \pi^2 N^{1/2}} \psia \msc(z) \left( \int_{\Oma} \del_{\bar{w}} \tilf (w) \msc' (w) \d u \d v\right)^2 + N^{\eps} | \xi |^2 \O ( N^{-1} y^{-1/2} (1+ \| f' \|_1)(1+ \| f'' \|_1 ))
\end{align}
This yields the claim. \qed

\subsubsection{Fourth order term}
\bel
We have,
\begin{align}
& \frac{1}{N^2} \sum_{ia}\frac{(1+\delta_{ia})^3}{N^2} \ee[ \del_{ia}^3 \ea ( G_{ia} (z) - \ee[ \del_{ia}^3 G_{ia} (z) ] ) ] \notag\\
 = & \frac{ - 12 \xi \i}{ \pi} \int_{\Oma} \del_{\bar{w}} \tilf (w) \msc'(w) \msc(w) \msc(z)^2 \d u \d v \notag \\
 + & N^{\eps} ( 1 + | \xi |^3 ) ( 1 + \| f' \|_1^3 ) \O ( N^{-1} y^{-1} ( 1 + \| f'' \|_1 ) )
\end{align}
\eel
\proof We write,
\begin{align}
& \frac{1}{N^2} \sum_{ia}\frac{(1+\delta_{ia})^3}{N^2} \ee[ \del_{ia}^3 \ea ( G_{ia} (z) - \ee[ \del_{ia}^3 G_{ia} (z) ] ) ] \notag\\
 = & \sum_{ia} \frac{ ( 1 + \delta_{ia} )^3}{N^2} \ee[ \ea ( \del_{ia}^3 G_{ia} (z) - \ee[ \del_{ia}^3 G_{ia} (z) ] ) ] \label{eqn:stein-4o1} \\
 + & \sum_{ia} \frac{ 3 ( 1 + \delta_{ia} )^3}{N^2} \ee[ ( \del_{ia}\ea ) ( \del_{ia}^2 G_{ia} (z) ) ] \label{eqn:stein-4o2} \\
 + &  \sum_{ia} \frac{ 3 ( 1 + \delta_{ia} )^3}{N^2} \ee[ ( \del_{ia}^2\ea ) ( \del_{ia} G_{ia} (z) ) ] \label{eqn:stein-4o3} \\
  + &  \sum_{ia} \frac{  ( 1 + \delta_{ia} )^3}{N^2} \ee[ ( \del_{ia}^3\ea ) (  G_{ia} (z) ) ] \label{eqn:stein-4o4}
\end{align}
For \eqref{eqn:stein-4o1}, the same logic used to estimate \eqref{eqn:cov-4o1} shows that with overwhelming probability,
\beq
\frac{1}{N^2} \sum_{ia} ( 1 + \delta_{ia} )^3  ( \del_{ia}^3 G_{ia} (z) - \ee[ \del_{ia}^3 G_{ia} (z) ] ) = N^{\eps} \O ( N^{-1} y^{-1} ).
\eeq
For \eqref{eqn:stein-4o2},  we see from the entry-wise local law  \eqref{eqn:entry-ll} that
\beq
\frac{1}{N^2} \sum_{ia} \del_{w} G_{ia} (w) ( 2 G_{ia}^3 (z) + 6 G_{ii} (z) G_{aa} (z) G_{ia} (z) ) = N^{\eps } \O ( N^{-1} y^{-1/2} v^{-3/2} )
\eeq
holds with overwhelming probability. Therefore by direct integration,
\beq
\sum_{ia} \frac{ 3 ( 1 + \delta_{ia} )^3}{N^2} \ee[ ( \del_{ia}\ea ) ( \del_{ia}^2 G_{ia} (z) ) ] = N^{\eps} \O ( |\xi | N^{-1} y^{-1/2} \| f'' \|_1 ).
\eeq
For \eqref{eqn:stein-4o3}, we first write it as
\begin{align}
& \sum_{ia} \frac{  ( 1 + \delta_{ia} )^3}{N^2} \ee[ ( \del_{ia}^2\ea ) ( \del_{ia} G_{ia} (z) ) ] \notag \\
= & \frac{ - 2 \xi \i }{ \pi N^2} \sum_{ia} \int_{\Oma} \del_{\bar{w}} \tilf (w) \ee[ \ea \del_w ( G_{ii} (w) G_{aa} (w) + G_{ia}(w)^2 ) ( G_{ii} (z) G_{aa} (z) + G_{ia}(z)^2 )   ] \notag\\
 +& \frac{ 4 \xi^2}{ \pi^2 N^2} \sum_{ia} \int_{\Oma^2} \del_{\bar{w}_1} \del_{\bar{w}_2} \tilf (w_1) \tilf (w_2) \del_{w_1} \del_{w_2} \ee[ \ea G_{ia} (w_1) G_{ia} (w_2) ( G_{ii}(z) G_{aa} (z) + G_{ia}(z)^2 )] 
\end{align}
Using the same logic as leading to the estimate of \eqref{eqn:cov-4o3}, we see that with overwhelming probability,
\begin{align}
& \frac{1}{N^2} \sum_{ia} \del_w ( G_{ii} (w) G_{aa} (w) + G_{ia}(w)^2 ) ( G_{ii} (z) G_{aa} (z) + G_{ia}(z)^2 ) \notag \\
= & 2 \msc'(w) \msc(w) \msc(z)^2 + N^{\eps} \O (N^{-1} y^{-1} v^{-2} ),
\end{align}
and so
\begin{align}
 & \frac{ - 2 \xi \i }{ \pi N^2} \sum_{ia} \int_{\Oma} \del_{\bar{w}} \tilf (w) \ee[ \ea \del_w ( G_{ii} (w) G_{aa} (w) + G_{ia}(w)^2 ) ( G_{ii} (z) G_{aa} (z) + G_{ia}(z)^2 )   ] \notag \\
 = & \frac{ - 4 \xi \i}{ \pi} \int_{\Oma} \del_{\bar{w}} \tilf (w) \msc'(w) \msc(w) \msc(z)^2 \d u \d v + N^{\eps} | \xi| \O (N^{-1} \| f'' \|_1 y^{-1} ).
\end{align}
By the entrywise local law \eqref{eqn:entry-ll} we have,
\beq
\frac{1}{N^2} \sum_{ia} \del_{w_1} \del_{w_2} G_{ia} (w_1) G_{ia}(w_2) ( G_{ii} (z) G_{aa} (z) + G_{ia} (z)^2 ) = N^{\eps} \O (N^{-1} v_1^{-3/2} v_2^{-3/2} ),
\eeq
with overwhelming probability. Therefore, by a similar argument as to the one that leads to \eqref{eqn:stein-ibp} that uses integration by parts we find,
\begin{align}
&\frac{ \xi^2}{ \pi^2 N^2} \sum_{ia} \int_{\Oma^2} \del_{\bar{w}_1} \del_{\bar{w}_2} \tilf (w_1) \tilf (w_2) \del_{w_1} \del_{w_2} \ee[ \ea G_{ia} (w_1) G_{ia} (w_2) ( G_{ii}(z) G_{aa} (z) + G_{ia}(z)^2 )]  \notag\\
= &N^{\eps} | \xi|^2 \O (N^{-1} ( 1 +  \| f''\|_1) ( 1 +  \| f'\|_1) )
\end{align}
This completes \eqref{eqn:stein-4o3}. 
For the final term \eqref{eqn:stein-4o4}, we may simply use the estimate,
\beq \label{eqn:delea-2}
|\del_{ia}^3 \ea | \leq (1+ | \xi |^3 ) (1 + \| f' \|_1^3 ) \left(  \delta_{ia} + N^{-1/2}( 1 + \| f''\|_1 )^{1/2} \right)
\eeq
which follows from direct calculation, the entry-wise local law \eqref{eqn:entry-ll} and an argument using integration by parts similar to the one leading to \eqref{eqn:stein-ibp}. This completes the proof.
\qed

\subsubsection{Fifth order term}

\bel We have,
\begin{align}
& \frac{1}{N^{5/2}} \sum_{ia} (1+\delta_{ia} )^4 \ee[ \del_{ia}^4 ( \ea (G_{ia} - \ee[ \del_{ia}^4 G_{ia} ] ) ] \notag\\
= & N^{\eps} (1+ |\xi|^4 ) ( 1 + \| f' \|_1 )^4 \O (N^{-1} y^{-1} ( 1 + \| f'' \|_1 ) )
\end{align}
\eel
\proof We write the fifth order term as,
\begin{align}
& \frac{1}{N^{5/2}} \sum_{ia} (1+\delta_{ia} )^4 \ee[ \del_{ia}^4 ( \ea (G_{ia} - \ee[ \del_{ia}^4 G_{ia} ] ) ] \notag\\
= & \frac{1}{N^{5/2}} \sum_{ia} ( 1 + \delta_{ia} )^4 \ee[ \ea ( \del_{ia}^4 G_{ia} -  \ee[ \del_{ia}^4 G_{ia} ] ) ] \notag\\
+ & \sum_{n=1}^4 \frac{1}{N^{5/2}} \sum_{ia} (1 + \delta_{ia} )^4 c_n \ee[ ( \del_{ia}^n \ea ) \del_{ia}^{4-n} G_{ia} ]. \label{eqn:stein-5o}
\end{align}
for some combinatorial factor $c_n$. From the entry-wise local law \eqref{eqn:entry-ll} we easily see that,
\beq
\frac{1}{N^{5/2}} \sum_{ia}  ( \del_{ia}^4 G_{ia} -  \ee[ \del_{ia}^4 G_{ia} ] ) = N^{\eps} \O (N^{-1} y^{-1/2} )
\eeq
with overwhelming probability. Consider now each of the terms in the last line of \eqref{eqn:stein-5o}. Note that for $m = 2k$ even we have,
\beq
| \del_{ia}^{2k} G_{ia} (z) | \leq C (\delta_{ia} + N^{\eps} (N y)^{-1/2} )
\eeq
with overwhelming probability.  Together with the estimate \eqref{eqn:delea-1} we see that the terms with $n=2, 4$ of the last line of \eqref{eqn:stein-5o} are $N^{\eps} \O ( N^{-1} y^{-1/2} (1 + |\xi| )^4 (1 + \| f' \|_1 )^4 )$.  For the terms with $n=1, 3$ note that the estimate \eqref{eqn:delea-2} holds also with LHS $| \del_{ia} \ea|$. Therefore, these terms are $N^{\eps}  ( 1 + | \xi|^3 ) (1 + \| f' \|_1^3) \O ( N^{-1/2} ( 1 + \| f''\|_1 )^{1/2} )$. This completes the proof. \qed

\subsubsection{Proof of Proposition \ref{prop:stein-self}}

This follows immediately from the lemmas of the previous subsections.

\subsection{Proof of Theorem \ref{thm:wig-smooth-lss}} \label{sec:wig-smooth-lss-proof}

The estimates for the expectation and variance are Lemma \ref{lem:wig-expect} and Theorem \ref{thm:var-global-2}, respectively. We turn to the esimate of the characteristic function. 
We have
\beq
\frac{ \d }{ \d \xi } \psia ( \xi ) = \frac{ \i}{ \pi} \int_{\Oma} \del_{\bar{z}} \tilf (z) N \ee[ \ea (m_N(z) - \ee[ m_N (z) ] ) ]  \d x \d y.
\eeq
From Proposition \ref{prop:stein-self} we then see that,
\begin{align}
\frac{ \d }{ \d \xi } \psia ( \xi ) &= (- \xi V_1 (f) + \frac{ \i 3 \xi^2}{ N^{1/2} } B_1 (f)) \psia ( \xi ) \notag\\
&+N^{\eps} ( 1+ |\xi |^6) ( 1 + \| f' \|_{1, w}^7) \O ( N^{-1} ( 1 + \| f'' \|_{1, w} ) )
\end{align}
where
\beq
B_1 (f) := \frac{2}{ 3 \pi^3} \left( \int_{\Oma} \del_{\bar{z}} \tilf (z) \msc' (z) \d x \d y \right)^3
\eeq
and 
\beq
V_1 (f) := \int_{\Oma^2} \del_{\bar{z}} \tilf (z) \del_{\bar{w}} \tilf (w) F(z, w) \d x \d y \d u \d v
\eeq
where $F(z, w)$ is defined by the function on the second and third lines of \eqref{eqn:cov-complex}. 
 Now, with $V(f)$ defined by \eqref{eqn:varf-def} and $B(f)$ defined by,
\beq
B(f) := \frac{2}{ 3 \pi^3} \left( \int_{\rr^2} \del_{\bar{z}} \tilf (z) \msc' (z) \d x \d y \right)^3=- \frac{1}{ 12 \pi^3} \left( \int_{-2}^2  f(x) \frac{x}{\sqrt{4-x^2}} \d x \right)^3,
\eeq
it is not hard to derive (see the proof of Theorem \ref{thm:var-global-2} for the estimate of $V_1(f)-V(f)$),
\beq
| V_1 (f) - V(f) | + |B_1(f) - B (f) | \leq N^{\eps} (1+ \| f'\|_1 )^2 N^{\mfa-1} (1+ \|f''\|_{1, w} ).
\eeq
Therefore,
\begin{align}
\frac{ \d }{ \d \xi } \psia ( \xi ) &= (- \xi V (f) + \frac{ 3 \i \xi^2 s_3}{ N^{1/2} } B (f)) \psia ( \xi ) \notag\\
&+N^{\mfa+ \eps} ( 1+ |\xi |^6) ( 1 + \| f' \|_{1, w}^7) \O ( N^{-1} ( 1 + \| f'' \|_{1, w} ) )
\end{align}
Therefore, using the fact that $B(f)$ is real we have,
\begin{align} \label{eqn:stein-1}
 & \left| \psia ( \xi ) - \exp\left[ - \xi^2 V(f)/2 + \i \xi^3 N^{-1/2} B(f) \right] \right| \notag\\
\leq & C N^{\mfa+ \eps} ( 1+ |\xi |^6) ( 1 + \| f' \|_{1, w}^7)  N^{-1} (1 + \| f'' \|_{1, w} )  \int_{0}^{\xi} \e^{ (s^2 - \xi^2) V(f) } \d s.
\end{align}
By Lemma \ref{lem:var-positive} we have that $V(f) \geq 0$ for $N \geq 4$. 
Therefore, the integral on the last line of \eqref{eqn:stein-1} is bounded by $|\xi|$. Therefore, by Lemma \ref{lem:HS-cutoff} we have the following estimate,
\begin{align}
& \left| \ee[\exp ( \i \xi ( \tr f (H) - \ee[ \tr f (H) ] ) ) ] - \exp\left[ - \xi^2 V(f)/2 + \i \xi^3 N^{-1/2} B(f) \right] \right| \notag\\
\leq & N^{\eps} (1+ | \xi|^6 + \| f' \|_{1, w}^7 ) N^{-1} (1 + \| f'' \|_{1, w} ),
\end{align}
and the claim follows. \qed

\appendix

\section{Derivatives of resolvent entries with respect to matrix elements} \label{a:derivs}

In this section, we collect some identities for derivatives of resolvent entries with respect to matrix entries. We have the general formula
\beq 
(1+\delta_{ia} ) \del_{ia} G_{jk} = - G_{ij} G_{ka} - G_{aj} G_{ki}.
\eeq
For first derivatives, this specializes to
\beq \label{eqn:delGia}
(1+ \delta_{ia} ) \del_{ia} G_{ia} = - G_{ii} G_{aa} - G_{ia}^2
\eeq
and
\beq \label{eqn:delGjj}
(1+ \delta_{ia} ) \del_{ia} G_{jj} = - 2 G_{ij} G_{ja}.
\eeq
For second derivatives, we get
\beq
 \label{eqn:del2Gia}
(1 + \delta_{ia} )^2 \del^2_{ia} (G_{ia} ) = 2 G_{ia}^3 + 6 G_{ii} G_{aa} G_{ia} .
\eeq
and
\beq \label{eqn:delGiadelGjj}
(1+ \delta_{ia} )^2 ( \del_{ia} G_{ia}(z) ) ( \del_{ia} G_{jj}(w) ) = 2 G_{ii} (z) G_{aa} (z) G_{ji} (w) G_{ja} (w) + 2 G_{ia} (z)^2 G_{ji} (w) G_{ja} (w) 
\eeq
and 
\beq \label{eqn:del2Gjj}
( 1 + \delta_{ia})^2 \del^2_{ia} G_{jj} = 4 G_{ij} G_{ja} G_{ai} + 2 G_{ij} G_{ji} G_{aa} + 2 G_{aj} G_{ja} G_{ii}.
\eeq
For the third derivatives, we have
\begin{align} \label{eqn:del3Gia}
 -(1+ \delta_{ia} )^3 \del_{ia}^3 G_{ia} &= 6 G_{ia}^4 + 36 G_{ia}^2 G_{ii} G_{aa} + 6 G_{ii}^2 G_{aa}^2,
\end{align}
as well as
\begin{align} \label{eqn:del3Gjj}
- (1+\delta_{ia} )^3 G_{jj} = 12  G_{ai}^2 G_{ij} G_{ja} + 12 G_{ja}^2 G_{ii} G_{ai} + 12 G_{ij}^2 G_{aa} G_{ai} + 12 G_{aj} G_{ij} G_{aa} G_{ii}
\end{align}


\section{Semicircle calculations}

\bel
We have,
\begin{align} \label{eqn:a-cov-lead}
\frac{ \msc (z)  - \msc (w) }{ z- w} = \frac{ \msc (z) \msc (w) }{ 1 - \msc(z) \msc(w) }.
\end{align}
\eel
\proof 
Recall,
\beq
\msc(z) = \frac{-z+\sqrt{z^2-4}}{2}, \qquad \msc(z) + z = - \frac{1}{\msc(z)}.
\eeq
Consider,
\begin{align}
\msc (z) - \frac{1}{\msc(w) } &= \msc (z) + w + \msc (w) = \frac{1}{2} \left( -z + \sqrt{z^2-4} +w + \sqrt{w^2-4} \right) \\
&= \frac{1}{2} \left( w -z + \frac{w^2 - z^2}{ \sqrt{ w^2-4}- \sqrt{z^2-4} } \right) =  \frac{1}{2} \left( w-z + \frac{(w-z)(w+z)}{ \sqrt{ w^2-4}- \sqrt{z^2-4} }  \right) \\
&= \frac{w-z}{2} \left( \frac{ \sqrt{w^2-4} - \sqrt{z^2-4} +w + z }{ \sqrt{w^2-4} - \sqrt{z^2-4} } \right) \\
&= (w-z) \left( \frac{ \msc (w) + w - \msc (z) }{ \sqrt{w^2-4} - \sqrt{z^2-4} }  \right) \\
&= (w-z) \left( \frac{-\msc (z) - (\msc (w))^{-1} }{ \sqrt{w^2-4} - \sqrt{z^2-4} } \right)
\end{align}
where in the second line we used difference of squares factoring $(\sqrt{a}-\sqrt{b})(\sqrt{a} + \sqrt{b} ) = a-b$.  Hence,
\begin{align}
1- \msc (z) \msc (w) &= - \msc (w) \left( \msc (z) - \frac{1}{\msc (w) } \right) \notag\\
 =& (-\msc (w)) (w-z) \left( \frac{-\msc (z) - (\msc (w))^{-1} }{ \sqrt{w^2-4} - \sqrt{z^2-4} } \right)  \notag \\
&= \frac{w-z}{ \sqrt{w^2-4} - \sqrt{z^2-4} }  ( 1+ \msc (z) \msc (w) ).
\end{align}
Therefore:
\beq
\frac{ \sqrt{w^2-4} - \sqrt{z^2-4} } {w-z}  = \frac{1+\msc (z) \msc (w)}{1 - \msc (z) \msc (w) }.
\eeq
Continuing, we have:
\begin{align}
\frac{\msc (z) - \msc (w)}{z-w} &= \frac{ -z + w + \sqrt{z^2-4} - \sqrt{w^2-4}}{2(z-w)} \\
&= - \frac{1}{2} + \frac{1}{2} \frac{ \sqrt{w^2-4} - \sqrt{z^2-4} } {w-z} \\
&= - \frac{1}{2} + \frac{1}{2}\left( \frac{1+\msc (z) \msc (w)}{1 - \msc (z) \msc (w) }\right)\\
&= \frac{\msc (z) \msc (w) }{1 - \msc (z) \msc (w) }
\end{align}
This yields the claim. \qed
 
\bel \label{lem:a-cov-2}
We have,
\begin{align}
\frac{1}{ z + 2 \msc (z) } \del_w \frac{ \msc (z)  - \msc (w) }{ z- w}  = - \frac{ \msc'(z) \msc'(w)}{( 1- \msc (z) \msc (w) )^2}
\end{align}
\eel
\proof This follows from the previous lemma and the fact that,
\beq \label{eqn:mder-quadratic}
\frac{ \msc(z)}{ z + 2 \msc (z) } = - \msc'(z),
\eeq
which follows from differentiating the quadratic equation satisfied by $\msc$. \qed

\bel \label{lem:a-cov-3}
We have, for $z = x +  \pm \i y$ and $w = u \pm \i v$, where $v, y >0$
\begin{align}
\left| \frac{ \msc'(z) \msc'(w)}{( 1- \msc (z) \msc (w) )^2} \right| \leq C \frac{|\msc'(z) \msc'(w) |}{ (x-u)^2 + y^2 + v^2 }
\end{align}
\eel
\proof From the fact that $| \msc(z) | \leq 1- c y$ we see that
\beq
| 1 - \msc(z) \msc(w) | \geq c (y + v).
\eeq
From \eqref{eqn:a-cov-lead} we see that
\beq
|1 - \msc(z) \msc (w) | \geq c |z-w| \geq c|x-u|.
\eeq
This completes the proof. \qed

\section{Hermite polynomial asymptotics} \label{a:hermite}

We collect here some standard results on the Hermite polynomials $h_n (x)$, all available from \cite{deift1999strong}.   
In the bulk we have that there exists a $\delta_a >0$, we have
\beq
h_n \left( \sqrt{\frac{n}{2}}x \right) \e^{- n \frac{x^2}{4}} =2 \sqrt{ \frac{1}{ \pi \sqrt{ 2 n } } } \frac{1}{ (4-x^2)^{1/4} } \left( \cos \left[ \frac{n}{2} F( x) - \frac{1}{2} \arcsin \left(\frac{x}{2} \right) \right] + \O (n^{-1} ) \right)
\eeq
uniformly in $-2+ \delta_a \leq x \leq 2 - \delta_a$ where
\beq
F(x) = \left| \int_x^2 \sqrt{|4 - y^2|} \d y \right|.
\eeq
In particular,
\beq
n^{1/4} \left| h_n \left( \sqrt{\frac{n}{2}}x \right) \e^{- n \frac{x^2}{4}}  \right| \leq C
\eeq
for $|x| \leq 2- \delta_a$.  Inside the spectrum, but near the edge we have 
\begin{align}
h_n \left( \sqrt{\frac{n}{2}}x \right) \e^{- n \frac{x^2}{4}} &= (2 n)^{-1/4} \notag\\
&\times \bigg\{ \left( \frac{x+2}{x-2} \right)^{1/4} \left[ \frac{3}{4} n F(x) \right]^{1/6} \Ai \left( - \left[\frac{3}{4}n F(x) \right]^{2/3} \right) ( 1 + \O (n^{-1} ) ) \notag\\
&- \left( \frac{x-2}{x+2} \right)^{1/4} \left[ \frac{3}{4} n F(x) \right]^{-1/6}  \Ai' \left( - \left[ \frac{3}{4}n F(x) \right]^{2/3} \right) ( 1 + \O (n^{-1} ) )\bigg\} ,
\end{align}
for $2- \delta_a < x < 2$. 
Above, $\Ai$ is the Airy function.  It is bounded on the real line and smooth. So for $r>1$, the following asymptotics are useful, and can be found in Chapter 11, Section 1 of \cite{olver1997asymptotics}, 
\begin{align}
\Ai (-r) &= \pi^{-1/2} r^{-1/4} \left\{ \cos \left[ \frac{2}{3} r^{3/2} - \frac{\pi}{4} \right] + \O (r^{-3/2} ) \right\} \notag\\
\Ai'(-r) &= \pi^{-1/2} r^{1/4} \left\{ \sin\left[ \frac{2}{3} r^{3/2} - \frac{\pi}{4} \right] + \O (r^{-3/2} ) \right\} .
\end{align}
Letting $x= 2 - y$ for $\delta_a  > y >0$ we see that 
\beq
n F(x) \asymp n y^{3/2}.
\eeq
So, using the above we see that 
\beq
\left| h_n \left( \sqrt{\frac{n}{2}}x \right) \e^{- n \frac{x^2}{4}}  \right| \leq \frac{C}{ y^{1/4} + n^{-1/6} }
\eeq
for $ x = 2 - y$, $0 < y < \delta_a$. 

Now, near the edge, outside the spectrum we have, for $ 2 < x < 2 + \delta_a$,
\begin{align}
h_n \left( \sqrt{\frac{n}{2}}x \right) \e^{- n \frac{x^2}{4}} &= (2 n)^{-1/4} \notag\\
&\times \bigg\{ \left( \frac{x+2}{x-2} \right)^{1/4} \left[ \frac{3}{4} n F(x) \right]^{1/6} \Ai \left(  \left[ \frac{3}{4} n F(x) \right]^{2/3} \right) ( 1 + \O (n^{-1} ) ) \notag\\
&- \left( \frac{x-2}{x+2} \right)^{1/4} \left[ \frac{3}{4} n F(x) \right]^{-1/6}  \Ai' \left(  \left[ \frac{3}{4} n F(x) \right]^{2/3} \right) ( 1 + \O (n^{-1} ) )\bigg\} .
\end{align}
In this region, we have from \cite{olver1997asymptotics} that,
\beq
| \Ai (r ) | \leq C \frac{ \e^{ - \frac{2}{3} r^{3/2} } }{r^{1/4}}, \qquad | \Ai'(r) | \leq C (1+r)^{1/4} \e^{ - \frac{2}{3} r^{3/2} } .
\eeq
Therefore, for $x = 2 + y$ and $0 < y < \delta_a$ we have
\beq
\left| h_n \left( \sqrt{\frac{n}{2}}x \right) \e^{- n \frac{x^2}{4}}  \right| \leq \frac{C\e^{ - c ny^{3/2}}}{ y^{1/4} + n^{-1/6} }
\eeq
The asymptotics in \cite{deift1999strong} show that the  above also extends to $ y > \delta_a$.

It is a straightforward calculation to see that the above upper bounds extend also to
\beq
y \to \left| n^{1/4} h_{n-1} \left(  \sqrt{\frac{n}{2}} (2 + y ) \right) \e^{ - n\frac{ (2+y)^2 }{4} } \right|.
\eeq

\section{Some integrals of Hermite polynomials}

Let $H_n (x)$ denote the physicist's Hermite polynomials,
\beq
H_n (x) = (-1)^n \e^{x^2} \left( \frac{ \d }{ \d x} \right)^n \e^{ - x^2} ,
\eeq
and the Hermite functions,
\beq
\varphi_n (x) = \frac{1}{ ( 2^n n! \sqrt{ \pi } )^{1/2} } \e^{ - x^2/2} H_n (x).
\eeq
In this appendix we will use several results from Gradshteyn and Ryzhik \cite{gradshteyn2014table}. There are many editions of this work, and all references are to the most recent 2014 version (8th edition) of \cite{gradshteyn2014table} unless noted.
\bel \label{lem:herm-int-1}
We have
\beq
\int \varphi_{2m } (x) \d x =  \frac{2^{1/2}}{m^{1/4}} \left( 1 + \O ( m^{-1} ) \right)
\eeq
\eel
\proof The integral equals,
\beq
\int \varphi_{2m} (t) \d t = \frac{1}{2^{m} ( (2m)!)^{1/2} \pi^{1/4} } \int_{\rr} \e^{ - t^2/2} H_{2m} (t )  \d t.
\eeq
Formula (7.376.1) of \cite{gradshteyn2014table} is,
\beq \label{eqn:a1}
\int_{\rr} \e^{ \i x y } \e^{ - x^2/2} H_n (x) \d  x = (2 \pi)^{1/2} H_n (y) \e^{-y^2/2} \i^{n}.
\eeq
We need to evaluate this at $n = 2m$ and $y=0$.
Formula (8.956.6) of \cite{gradshteyn2014table} is,
\beq \label{eqn:a2}
H_{2m} (0) = (-1)^m 2^m (2m-1)!! 
\eeq
For the double factorial we have,
\beq \label{eqn:a3}
(2m-1)!! = \frac{ (2m-1)!}{2^{m-1} (m-1)!}
\eeq
Therefore,
\begin{align}
\int_\rr \varphi_{2m} (t) \d t &=\left(  \frac{1}{2^{m} ((2m)!)^{1/2} \pi^{1/4} }  \right) \left( (2 \pi)^{1/2} (-1)^m \right) \left( (-1)^m 2^m \frac{ (2m-1)!}{2^{m-1} (m-1)!} \right) \notag\\
&=  2^{1/2} \pi^{1/4} \frac{ (2m-1)!}{2^{m-1} ( (2m)!)^{1/2} (m-1)!} \notag\\
&= 2^{1/2} \pi^{1/4} \frac{ ( (2m)!)^{1/2}}{2^m  (m)!}.
\end{align}
The first line is just plugging in the formulas \eqref{eqn:a1}, \eqref{eqn:a2} and \eqref{eqn:a3}. The second line is just a few cancellations. In the third line we simplified the factorial terms.  Stirling's approximation is,
\beq
n! = ( 2 \pi n)^{1/2} n^n \e^{-n} \left( 1 + \O (n^{-1} ) \right)
\eeq
Therefore,
\begin{align}
\int_\rr \varphi_{2m} (t) \d t &= 2^{1/2} \pi^{1/4} \frac{ ( (2m)!)^{1/2}}{2^m  (m)!} \notag\\
&=  2^{1/2} \pi^{1/4} \frac{ ( 2 \pi )^{1/4} ( 2 m )^{1/4} (2 m)^{m} \e^{-m}}{2^m ( 2 \pi )^{1/2} m^{1/2} m^m \e^{-m} } \left( 1 + \O ( m^{-1} ) \right) \notag\\
&= \frac{2^{1/2}}{m^{1/4}} \left( 1 + \O ( m^{-1} ) \right) \label{eqn:d1}
\end{align}
This yields the claim. \qed

\bel \label{lem:herm-int-2}
We have,
\begin{align}
\left( \frac{ 2 m + 1 }{2} \right)^{1/2} \int_0^\infty \varphi_{2m+1} (x) \d x &= \frac{ m^{1/4}}{ 2^{1/2}}  \left( 1 + \O( m^{-1/2} ) \right)
\end{align}
\eel
\proof 
Clearly $(\frac{2m+1}{2} )^{1/2} = m^{1/2} ( 1 + \O (m^{-1} ))$. By definition,
\begin{align}
m^{1/2} \int_0^\infty \varphi_{2m+1} (t) \d t =
&= \frac{m^{1/2}}{ (2^{2m+1})^{1/2} ((2m+1)!)^{1/2} \pi^{1/4} } \int_0^\infty \e^{ - t^2/2} H_{2m+1} (t) \d t.
\end{align}
We require the formula,
\beq \label{eqn:odd-int}
\int_0^\infty \e^{ -2 \alpha x^2} x^\nu H_{2n+1} (x) \d x = (-1)^n 2^{2n-\nu/2} \frac{ \Gamma ( \nu/2+1) \Gamma (n+3/2)}{ \pi^{1/2} \alpha^{\nu/2+1}} F ( -n, \nu/2+1 ; 3/2 ; (2 \alpha)^{-1} ),
\eeq
where $F$ is the (Gauss) hypergeometric function (also denoted $_2 F_1$).  It is given in the 2007 version (7th edition) of Gradshteyn and Ryzhik \cite{gradshteyn2007table} as formula (7.376.3). It was then erroneously changed in the 2014 version \cite{gradshteyn2014table}, and then corrected back to the 2007 (7th) edition formula in the erratum \cite{veestraeten2015some}$^1$.{\let\thefootnote\relax\footnotetext{$1$. The source of confusion seems to be that the original source \cite{buchholz2013confluent} itself contains an error; note that the formula (13.18b) in that work is derived via substitution of (13.3b) into (10.4b). The $\mu$ in (10.4b) equals $\mu=1/2$ and so the $\Gamma(n+1/2)$ in (13.18b) should be $\Gamma(n+\mu+1/2) = \Gamma(n+3/2)$. In any case, it is also not hard to directly check \eqref{eqn:odd-int} from the relation $H_{2n+1}(x) = (-1)^n (2n+1)! (n!)^{-1} (2x)_1 F_1 (-n, 3/2; x^2)$, where $_1 F_1$ is the confluent hypergeometric function of the first kind, and direct integration of the series term-by-term. }}

 Evaluating this with $n=m$ and $\nu  = 0$ and $\alpha = 1/4$ we see that, 
\begin{align} \label{eqn:b3}
& \frac{m^{1/2}}{ (2^{2m+1})^{1/2} ((2m+1)!)^{1/2} \pi^{1/4} } \int_0^\infty \e^{ - t^2/2} H_{2m+1} (t) \d t \notag\\
= & \frac{m^{1/2}}{ (2^{2m+1})^{1/2} ((2m+1)!)^{1/2} \pi^{1/4} } (-1)^m 2^{2m} \frac{ \Gamma(m+3/2)}{ \pi^{1/2} 4^{-1}} F(-m, 1 ; 3/2 ; 2),
\end{align}
where we simplified $\Gamma(1) = 1$.   From (8.339) of \cite{gradshteyn2014table} and \eqref{eqn:a3},
\beq
\Gamma \left( \frac{1}{2} + n \right) = \frac{ \sqrt{ \pi}}{2^n} ( 2n -1 )!! =  \frac{ (2n)!}{4^n n!} \sqrt{ \pi }.
\eeq
Evaluating this with $n = m+1$ we get,
\beq \label{eqn:b1}
\Gamma \left( \frac{3}{2} + m \right) = \frac{ (2 (m+1) )!}{4^{m+1} (m+1)!} \sqrt{ \pi }.
\eeq
We use the following integral representation of the hypergeometric function, (9.111) of \cite{gradshteyn2014table}
\beq
F (a, b ; c ; z) = \frac{1}{ B ( b, c-b) } \int_0^1 t^{b-1} (1-t)^{c-b-1}(1-tz)^{-a} \d t,
\eeq
where the beta function is,
\beq
B(x, y) = \frac{ \Gamma(x) \Gamma(y) }{ \Gamma(x+y) }.
\eeq
Therefore,
\begin{align} \label{eqn:b2}
F(-m, 1, 3/2, 2) &= \frac{ \Gamma (3/2) }{ \Gamma(1) \Gamma ( 1/2 ) } \int_0^1 (1 - t)^{-1/2} (1 - 2 t)^m \d t \notag\\
&= \frac{1}{2}   \int_0^1 (1 - t)^{-1/2} (1 - 2 t)^m \d t 
\end{align}
So,
\begin{align}
& \frac{m^{1/2}}{ (2^{2m+1})^{1/2} ((2m+1)!)^{1/2} \pi^{1/4} } \int_0^\infty \e^{ - t^2/2} H_{2m+1} (t) \d t \notag\\
= & \left( \frac{m^{1/2}}{ (2^{2m+1})^{1/2} ((2m+1)!)^{1/2} \pi^{1/4} } \right) \times \left( \frac{ (-1)^m 2^{2m}}{ \pi^{1/2} 4^{-1} }   \frac{ (2 (m+1) )!}{4^{m+1} (m+1)!} \sqrt{ \pi } \right) \notag\\
\times & \left( \frac{1}{2}   \int_0^1 (1 - t)^{-1/2} (1 - 2 t)^m \d t  \right) \notag\\
= & \pi^{-1/4} (-1)^m \frac{ m^{1/2} (  2 (m+1) )! }{ 2^m 2^{1/2} (m+1)! ( ( 2 m + 1 )! )^{1/2} } \left( \frac{1}{2}   \int_0^1 (1 - t)^{-1/2} (1 - 2 t)^m \d t  \right) \notag\\
= & \pi^{-1/4} \frac{ m^{1/2} (2m+2)(2m+1) ( (2m)!)^{1/2} }{2^{m} 2^{3/2} (m+1) m! ( 2m + 1)^{1/2}} \left( (-1)^m   \int_0^1 (1 - t)^{-1/2} (1 - 2 t)^m \d t\right) \notag\\
= & \pi^{-1/4} m \left( \frac{ (  ( 2m)!)^{1/2} }{ 2^m m!} \right)  \left( (-1)^m   \int_0^1 (1 - t)^{-1/2} (1 - 2 t)^m \d t\right) \left( 1 + \O ( m^{-1} ) \right) \label{eqn:c3}
\end{align}
The first equality is direct substitution of \eqref{eqn:b1} and \eqref{eqn:b2} into \eqref{eqn:b3}. In the second equality we cleaned up the powers of $2$ and $\pi$. In the third equality we simplified the factorials to only have $(2m)!$ and $m!$. The third asymptotic is using things like $2m +1 \sim 2m$, etc., and simplifying the powers of $2$.  Via Stirling's approximation,
\beq \label{eqn:c1}
\frac{ (  ( 2m)!)^{1/2} }{ 2^m m!} = \frac{ (2 \pi )^{1/4} ( 2m )^{1/4} ( 2m )^{m} \e^{-m} }{ 2^m (2 \pi )^{1/2} m^{1/2} m^m \e^{-m} }  \left( 1 + \O ( m^{-1} ) \right)  = \pi^{-1/4} m^{-1/4}  \left( 1 + \O ( m^{-1} ) \right)  .
\eeq
Now consider the integral,
\beq
 (-1)^m   \int_0^1 (1 - t)^{-1/2} (1 - 2 t)^m \d t = \int_0^1 (2t - 1)^m (1- t)^{-1/2} \d t.
\eeq
It is easy to see that,
\beq
\int_0^{1/2} |1 - 2 t|^m (1 - t)^{-1/2} \leq 2 \int_0^{1/2} (1 - 2t )^m \leq \frac{C}{m}
\eeq
by exact calculation. For the other half of the integral,
\begin{align}
\int_{1/2}^1  (2t - 1)^m (1- t)^{-1/2} \d t &= \int_0^{1/2} (1 - 2u )^{m} u^{-1/2} \d u \notag\\
&= \int_0^{1/2} \e^{ m \log (1 - 2 u ) } u^{-1/2} \d u \notag\\
&= \frac{1}{2^{1/2} m } \int_0^\infty \exp \left[ - x(1+m^{-1} ) \right] \frac{1}{ ( 1 - \e^{-x/m} )^{1/2} } \d x ,
\end{align}
via the change of variable $x = - m \log (1 - 2 u )$. 
It is straightforward to check that,
\beq
\left| \frac{1}{ ( 1 - \e^{ - s } )^{1/2} } - \frac{1}{s^{1/2} } \right| \leq C
\eeq
for all $s>0$.  Therefore,
\begin{align}
 & \frac{1}{2^{1/2} m } \int_0^\infty \exp \left[ - x(1+m^{-1} ) \right] \frac{1}{ ( 1 - \e^{-x/m} )^{1/2} } \d x \notag \\
 =&  \frac{1}{ 2^{1/2} m^{1/2} } \int_0^\infty \exp\left[ - x (1 + m^{-1} ) \right] x^{-1/2} \d x + \O ( m^{-1} ) \notag\\
 = &\frac{1}{ 2^{1/2} m^{1/2} } \int_0^\infty \exp\left[ - x \right] x^{-1/2} \d x + \O ( m^{-1} ) 
\end{align}
The final integral can be calculated exactly by relating it to the Gaussian integral,
\beq
\frac{1}{ 2^{1/2} m^{1/2} } \int_0^\infty \exp\left[ - x \right] x^{-1/2} \d x = \frac{ \pi^{1/2}}{ 2^{1/2} m^{1/2}}.
\eeq
Summarizing,
\beq \label{eqn:c2}
 (-1)^m   \int_0^1 (1 - t)^{-1/2} (1 - 2 t)^m \d t = \frac{\sqrt{ \pi}}{ \sqrt{ 2 m }} ( 1  + \O (m^{-1/2} ) ).
\eeq
Plugging \eqref{eqn:c1} and \eqref{eqn:c2} into \eqref{eqn:c3} we obtain,
\begin{align}
  &  \frac{m^{1/2}}{ (2^{2m+1})^{1/2} ((2m+1)!)^{1/2} \pi^{1/4} } \int_0^\infty \e^{ - t^2/2} H_{2m+1} (t) \d t \notag \\
 = & \pi^{-1/4} m \left( \pi^{-1/4} m^{-1/4} \right) ( \pi^{1/2} ( 2 m )^{-1/2} ) \left( 1 + \O( m^{-1/2} ) \right) \notag\\
 = & \frac{ m^{1/4}}{ 2^{1/2}}  \left( 1 + \O( m^{-1/2} ) \right).
\end{align}
This completes the proof. \qed

\section{Verification that the variance functional is positive}

\bel \label{lem:var-positive}
The functional $V ( \varphi )$ defined in \eqref{eqn:main-V-def} is positive for $N \geq 4$.
\eel
\proof With the functional $V( \varphi)$ as in \eqref{eqn:main-V-def} we have, (see, e.g., \cite[Theorem 2.4]{johansson1998fluctuations})
\beq
V ( \varphi ) = \frac{1}{2} \sum_{k=1}^\infty k c_k^2 + \frac{s_4}{2} c_2^2 - \frac{2 s_3}{N^{1/2}} c_1 c_2
\eeq
where $c_i$ is the $i$th coefficient of $\varphi$ in the Chebyshev polynomial basis,
\beq
c_k = \frac{1}{ \pi} \int_{-\pi}^{\pi} \varphi ( 2 \cos ( \theta ) ) \cos ( k \theta ) \d \theta = \frac{2}{ \pi}\int_{-1}^1 \varphi ( 2 x ) T_k (x) \frac{ \d x }{ \sqrt{ 1 - x^2}}
\eeq
and $T_k$ is the $k$th Chebyschev polynomial of the first kind, $T_k ( \cos ( \theta ) ) = \cos (k \theta )$.

Now for any centered random variable with variance $1$ we have the inequality, (see, for example, the remark after \cite[Theorem 2.2]{erdHos2012bulk})
\beq
s_4 \geq s_3^2 -2.
\eeq
Hence,
\begin{align}
V( \varphi ) &\geq \frac{c_1^2}{2} + c_2^2 + \frac{s_4}{2} c_2^2 - 2 \frac{ s_3}{N^{1/2}} c_1 c_2 \notag\\
& \geq \frac{c_1^2}{2} + \frac{s_3^2}{2} c_2^2 -  2 \frac{ s_3}{N^{1/2}} c_1 c_2 \notag\\
& \geq \frac{c_1^2}{2} + \frac{s_3^2}{2} c_2^2 - \frac{s_3^2}{2} c_2^2 - 2 \frac{c_1^2}{N}.
\end{align}
We therefore conclude the proof. \qed



\bibliography{mybib}{}

\begin{thebibliography}{10}

\bibitem{az}
G.~W. Anderson and O.~Zeitouni.
\newblock A {CLT} for a band matrix model.
\newblock {\em Probab. Theory Related Fields}, 134(2):283--338, 2006.

\bibitem{arkharov1971limit}
L.~Arkharov.
\newblock Limit theorems for the characteristic roots of a sample covariance
  matrix.
\newblock {\em Dokl. Akad. Nauk}, 199(5):994--997, 1971.

\bibitem{bahouri2011fourier}
H.~Bahouri, J.-Y. Chemin, and R.~Danchin.
\newblock {\em Fourier analysis and nonlinear partial differential equations},
  volume 343.
\newblock Springer, 2011.

\bibitem{bai2008clt}
Z.~D. Bai and J.~W. Silverstein.
\newblock {CLT} for linear spectral statistics of large-dimensional sample
  covariance matrices.
\newblock In {\em Ann. Probab.}, volume~32, pages 553--605. World Scientific,
  2004.

\bibitem{bao2021quantitative}
Z.~Bao and Y.~He.
\newblock Quantitative {CLT} for linear eigenvalue statistics of {W}igner
  matrices.
\newblock {\em Preprint, arXiv:2103.05402}, 2021.

\bibitem{bekerman2018mesoscopic}
F.~Bekerman and A.~Lodhia.
\newblock Mesoscopic central limit theorem for general $\beta$-ensembles.
\newblock {\em Ann. Inst. Henri Poincar{\'e} Probab. Stat.}, 54(4):1917--1938,
  2018.

\bibitem{benaych2019lectures}
F.~Benaych-Georges and A.~Knowles.
\newblock Lectures on the local semicircle law for wigner matrices.
\newblock 2019.

\bibitem{borodin2017gaussian}
A.~Borodin, V.~Gorin, and A.~Guionnet.
\newblock Gaussian asymptotics of discrete $\beta$-ensembles.
\newblock {\em Publ. Math. Inst. Hautes {\'E}tudes Sci.}, 125(1):1--78, 2017.

\bibitem{bourgade2018extreme}
P.~Bourgade.
\newblock Extreme gaps between eigenvalues of {W}igner matrices.
\newblock {\em J. Eur. Math. Soc., to appear}, 2021.

\bibitem{fixed-wig}
P.~Bourgade, L.~Erd{\H{o}}s, H.-T. Yau, and J.~Yin.
\newblock Fixed energy universality for generalized {W}igner matrices.
\newblock {\em Comm. Pure Appl. Math.}, 69(10):1815--1881, 2016.

\bibitem{bourgade2016fixed}
P.~Bourgade, L.~Erd{\H{o}}s, H.-T. Yau, and J.~Yin.
\newblock Fixed energy universality for generalized {W}igner matrices.
\newblock {\em Comm. Pure Appl. Math.}, 69(10):1815--1881, 2016.

\bibitem{blz}
P.~Bourgade, P.~Lopatto, and O.~Zeitouni.
\newblock {\em In preparation}, 2022.

\bibitem{bourgade2019gaussian}
P.~Bourgade and K.~Mody.
\newblock Gaussian fluctuations of the determinant of {W}igner matrices.
\newblock {\em Electron. J. Probab.}, 24:1--28, 2019.

\bibitem{buchholz2013confluent}
H.~Buchholz.
\newblock {\em The confluent hypergeometric function: with special emphasis on
  its applications}, volume~15.
\newblock Springer Science \& Business Media, 2013.

\bibitem{chatterjee2009fluctuations}
S.~Chatterjee.
\newblock Fluctuations of eigenvalues and second order poincar{\'e}
  inequalities.
\newblock {\em Probab. Theory Related Fields}, 143(1):1--40, 2009.

\bibitem{cipolloni2020functional}
G.~Cipolloni, L.~Erd{\H{o}}s, and D.~Schr{\"o}der.
\newblock Functional central limit theorems for {W}igner matrices.
\newblock {\em Preprint, arXiv:2012.13218}, 2020.

\bibitem{cipolloni2021spectral}
G.~Cipolloni, L.~Erd{\H{o}}s, and D.~Schr{\"o}der.
\newblock On the spectral form factor for random matrices.
\newblock {\em Preprint, arXiv:2109.06712}, 2021.

\bibitem{costin1995gaussian}
O.~Costin and J.~L. Lebowitz.
\newblock Gaussian fluctuation in random matrices.
\newblock {\em Phys. Rev. Lett.}, 75(1):69, 1995.

\bibitem{dallaporta2011note}
S.~Dallaporta and V.~Vu.
\newblock A note on the central limit theorem for the eigenvalue counting
  function of {W}igner matrices.
\newblock {\em Electron. Comm. Probab.}, 16:214--322, 2011.

\bibitem{deift1999strong}
P.~Deift, T.~Kriecherbauer, K.~T.-R. McLaughlin, S.~Venakides, and X.~Zhou.
\newblock Strong asymptotics of orthogonal polynomials with respect to
  exponential weights.
\newblock {\em Comm. Pure. Appl. Math.}, 52(12):1491--1552, 1999.

\bibitem{deleporte2021universality}
A.~Deleporte and G.~Lambert.
\newblock Universality for free fermions and the local weyl law for
  semiclassical schr{\"o}dinger operators.
\newblock {\em Preprint arXiv:2109.02121}, 2021.

\bibitem{diaconis2001linear}
P.~Diaconis and S.~Evans.
\newblock Linear functionals of eigenvalues of random matrices.
\newblock {\em Trans. Amer. Math. Soc.}, 353(7):2615--2633, 2001.

\bibitem{diaconis1994eigenvalues}
P.~Diaconis and M.~Shahshahani.
\newblock On the eigenvalues of random matrices.
\newblock {\em J. Appl. Probab.}, 31(A):49--62, 1994.

\bibitem{wegner}
L.~Erd{\H{o}}s, B.~Schlein, and H.-T. Yau.
\newblock Wegner estimate and level repulsion for {W}igner random matrices.
\newblock {\em Int. Math. Res. Not. IMRN}, 2010(3):436--479, 2010.

\bibitem{erdHos2010wegner}
L.~Erd{\H{o}}s, B.~Schlein, and H.-T. Yau.
\newblock Wegner estimate and level repulsion for {W}igner random matrices.
\newblock {\em Int. Math. Res. Not. IMRN}, 2010(3):436--479, 2010.

\bibitem{local-relax}
L.~Erd{\H{o}}s, B.~Schlein, and H.-T. Yau.
\newblock Universality of random matrices and local relaxation flow.
\newblock {\em Invent. Math.}, 185(1):75--119, 2011.

\bibitem{erdHos2011universality}
L.~Erd{\H{o}}s, B.~Schlein, and H.-T. Yau.
\newblock Universality of random matrices and local relaxation flow.
\newblock {\em Invent. Math.}, 185(1):75--119, 2011.

\bibitem{erdHos2017dynamical}
L.~Erd{\H{o}}s and H.-T. Yau.
\newblock {\em A dynamical approach to random matrix theory}, volume~28.
\newblock American Mathematical Soc., 2017.

\bibitem{erdHos2012bulk}
L.~Erd{\H{o}}s, H.-T. Yau, and J.~Yin.
\newblock Bulk universality for generalized {W}igner matrices.
\newblock {\em Probab. Theory Related Fields}, 154(1):341--407, 2012.

\bibitem{forrester2006asymptotic}
P.~Forrester, N.~Frankel, and T.~Garoni.
\newblock Asymptotic form of the density profile for {G}aussian and {L}aguerre
  random matrix ensembles with orthogonal and symplectic symmetry.
\newblock {\em J. Math. Phys.}, 47(2):023301, 2006.

\bibitem{garoni2005asymptotic}
T.~Garoni, P.~J. Forrester, and N.~Frankel.
\newblock Asymptotic corrections to the eigenvalue density of the {GUE} and
  {LUE}.
\newblock {\em J. Math. Phys.}, 46(10):103301, 2005.

\bibitem{girko2012theory}
V.~L. Girko.
\newblock {\em Theory of stochastic canonical equations}, volume 535.
\newblock Springer Science \& Business Media, 2012.

\bibitem{gradshteyn2007table}
I.~S. Gradshteyn and I.~Ryzhik.
\newblock Table of integrals, series, and products, 7th edition.
\newblock {\em Elsevier/Academic Press, Amsterdam}, 48:1171, 2007.

\bibitem{gradshteyn2014table}
I.~S. Gradshteyn and I.~M. Ryzhik.
\newblock {\em Table of integrals, series, and products}.
\newblock Academic press, 2014.

\bibitem{guionnet2002large}
A.~Guionnet.
\newblock Large deviations upper bounds and central limit theorems for
  non-commutative functionals of {G}aussian large random matrices.
\newblock {\em Ann. Inst. Henri Poincar{\'e} Probab. Stat.}, 38(3):341--384,
  2002.

\bibitem{gustavsson2005gaussian}
J.~Gustavsson.
\newblock Gaussian fluctuations of eigenvalues in the {GUE}.
\newblock {\em Ann. Inst. H. Poincar{\'e} Probab. Statist.}, 41(2):151--178,
  2005.

\bibitem{he2017mesoscopic}
Y.~He and A.~Knowles.
\newblock Mesoscopic eigenvalue statistics of {W}igner matrices.
\newblock {\em Ann. Appl. Probab.}, 27(3):1510--1550, 2017.

\bibitem{he2020mesoscopic}
Y.~He and A.~Knowles.
\newblock Mesoscopic eigenvalue density correlations of {W}igner matrices.
\newblock {\em Probab. Theory Related Fields}, 177(1):147--216, 2020.

\bibitem{huang2019rigidity}
J.~Huang and B.~Landon.
\newblock Rigidity and a mesoscopic central limit theorem for {D}yson
  {B}rownian motion for general $\beta$ and potentials.
\newblock {\em Probab. Theory Related Fields}, 175(1):209--253, 2019.

\bibitem{johansson1998fluctuations}
K.~Johansson.
\newblock On fluctuations of eigenvalues of random {H}ermitian matrices.
\newblock {\em Duke Math. J.}, 91(1):151--204, 1998.

\bibitem{jonsson1982some}
D.~Jonsson.
\newblock Some limit theorems for the eigenvalues of a sample covariance
  matrix.
\newblock {\em J. Multivariate Anal.}, 12(1):1--38, 1982.

\bibitem{kopel2015regularity}
P.~Kopel.
\newblock Regularity conditions for convergence of linear statistics of {GUE}.
\newblock {\em Preprint, arXiv:1510.02988}, 2015.

\bibitem{landon2021single}
B.~Landon, P.~Lopatto, and P.~Sosoe.
\newblock Single eigenvalue fluctuations of general {W}igner-type matrices.
\newblock {\em Preprint, arXiv:2105.01178}, 2021.

\bibitem{meso}
B.~Landon and P.~Sosoe.
\newblock Applications of mesoscopic {CLT}s in random matrix theory.
\newblock {\em Ann. Appl. Probab.}, 30(6):2769--2795, 2020.

\bibitem{lee2018local}
J.~O. Lee and K.~Schnelli.
\newblock Local law and {T}racy--{W}idom limit for sparse random matrices.
\newblock {\em Probab. Theory Related Fields}, 171(1):543--616, 2018.

\bibitem{li2021central}
Y.~Li, K.~Schnelli, and Y.~Xu.
\newblock Central limit theorem for mesoscopic eigenvalue statistics of
  deformed {W}igner matrices and sample covariance matrices.
\newblock {\em Ann. Inst. Henri Poincar{\'e} Probab. Stat.}, 57(1):506--546,
  2021.

\bibitem{lytova2009central}
A.~Lytova and L.~Pastur.
\newblock Central limit theorem for linear eigenvalue statistics of random
  matrices with independent entries.
\newblock {\em Ann. Probab.}, 37(5):1778--1840, 2009.

\bibitem{mehta2004random}
M.~L. Mehta.
\newblock {\em Random matrices}.
\newblock Elsevier, 2004.

\bibitem{olver1997asymptotics}
F.~Olver.
\newblock {\em Asymptotics and special functions}.
\newblock AK Peters/CRC Press, 1997.

\bibitem{orourke}
S.~O'Rourke.
\newblock Gaussian fluctuations of eigenvalues in {W}igner random matrices.
\newblock {\em J. Stat. Phys.}, 138(6):1045--1066, 2010.

\bibitem{pastur2011eigenvalue}
L.~A. Pastur and M.~Shcherbina.
\newblock {\em Eigenvalue distribution of large random matrices}.
\newblock Number 171. American Mathematical Soc., 2011.

\bibitem{shcherbina2011central}
M.~Shcherbina.
\newblock Central limit theorem for linear eigenvalue statistics of the
  {W}igner and sample covariance random matrices.
\newblock {\em J. Math. Phys., Analysis, Geometry}, (7):176--192, 2011.

\bibitem{soshnikov2002gaussian}
A.~Soshnikov.
\newblock Gaussian limit for determinantal random point fields.
\newblock {\em Ann. Probab.}, 30(1):171--187, 2002.

\bibitem{soshnikov2000gaussian}
A.~B. Soshnikov.
\newblock Gaussian fluctuation for the number of particles in {A}iry, {B}essel,
  sine, and other determinantal random point fields.
\newblock {\em J. Stat. Phys.}, 100(3):491--522, 2000.

\bibitem{sosoewong}
P.~Sosoe and P.~Wong.
\newblock Regularity conditions in the {CLT} for linear eigenvalue statistics
  of {W}igner matrices.
\newblock {\em Ad. Math.}, 249:37--87, 2013.

\bibitem{veestraeten2015some}
D.~Veestraeten et~al.
\newblock Some remarks, generalizations and misprints in the integrals in
  {G}radshteyn and {R}yzhik.
\newblock {\em Scientia. Series A. Mathematical Sciences}, 26:115--131, 2015.

\end{thebibliography}
\bibliographystyle{abbrv}

\end{document}